\input ./style/arxiv-general.cfg
\documentclass[MSNbibl,number,citesort,seceqn,dvips]{arxbj}
\makeatletter
   \@ifpackageloaded{graphicx}{}{\usepackage{graphicx}}
\makeatother
\volume{22}
\issue{4}
\pubyear{2016}
\firstpage{2325}
\lastpage{2371}
\doi{10.3150/15-BEJ731}% Updated by VTEXPTS2LaTeX.exe, 14.08.2015 14:26
\docsubty{FLA}

\makeatletter
\newcommand{\rrvert}{\vert}
\newcommand{\llvert}{\vert}
\renewcommand{\mid}{|}
\newtheorem{teo}{Theorem}[section]
\newremark{rem}{Remark}[section]
\newtheorem{theorem}{Theorem}
\newtheorem{corollary}{Corollary}
\newtheorem{lemma}{Lemma}
\newremark{remark}{Remark}

\newcommand{\eqref}[1]{(\ref{#1})}
\newcommand{\pibar}{\overline{\Pi}}
\newcommand{\lambar}{\overline{\Lambda}}
\newcommand{\rhobar}{\overline{\rho}}
\newcommand{\veps}{\varepsilon}
\newcommand{\rmd}{{\mathrm d}}
\newcommand{\rmi}{{\mathrm i}}
\newcommand{\FFF}{{\mathcal F}}
\newcommand{\EEEE}{\mathfrak{E}}
\newcommand{\bfeins}{\mathbf{1}}
\newcommand{\myRR}{\mathbb{R}}
\newcommand{\NN}{\mathbb{N}}
\newcommand{\VV}{\mathbb{V}}
\newcommand{\XX}{\mathbb{X}}
\newcommand{\YY}{\mathbb{Y}}
\newcommand{\SSSS}{\mathfrak{S}}
\newcommand{\UUUU}{\mathfrak{U}}
\makeatother

\begin{document}
\begin{frontmatter}

\title{Distributional representations and dominance of a L\'evy
process over its maximal jump processes}
\runtitle{Distributional representations of a L\'{e}vy process}

\begin{aug}
%%%% inicialai - be tarpu
% Corresponding author: Boris Buchmann - boris.buchmann@anu.edu.au% Updated by VTEXPTS2LaTeX.exe, 17.08.2015 08:35
%Updated by VTEXPTS2LaTeX.exe, 14.08.2015 14:26
\author[A]{\inits{B.}\fnms{Boris}~\snm{Buchmann}\corref{}\thanksref{A,e1}\ead[label=e1,mark]{Boris.Buchmann@anu.edu.au}},
\author[B]{\inits{Y.}\fnms{Yuguang}~\snm{Fan}\thanksref{B}\ead[label=e2]{Yuguang.Fan@unimelb.edu.au}}
\and
\author[A]{\inits{R.A.}\fnms{Ross A.}~\snm{Maller}\thanksref{A,e3}\ead[label=e3,mark]{Ross.Maller@anu.edu.au}\ead
[label=u1,url]{www.foo.com}}
%%\runauthor{} %% auto
%\dedicated{}
\address[A]{Research School of Finance, Actuarial Studies \& Statistics,
Mathematical Sciences Institute, Australian National University, Australia.
\printead{e1}; \printead*{e3}}
\address[B]{School of Mathematics \& Statistics, University of Melbourne,
ARC Centre of Excellence for Mathematics \& Statistical Frontiers, Australia.
\printead{e2}}
\end{aug}

% HISTORY:
%
\received{\smonth{10} \syear{2014}}% Updated by VTEXPTS2LaTeX.exe,
%14.08.2015 14:26
%
\revised{\smonth{3} \syear{2015}}% Updated by VTEXPTS2LaTeX.exe,
%14.08.2015 14:26

% ABSTRACT
%
\begin{abstract}
Distributional identities for a L\'evy process $X_t$, its quadratic
variation process $V_t$ and its
maximal jump processes, are derived, and used to make ``small time''
(as $t\downarrow0$) asymptotic
comparisons between them. The representations are constructed using
properties of the underlying
Poisson point process of the jumps of $X$. Apart from providing insight
into the connections between $X$, $V$, and their maximal jump
processes, they enable investigation of a great variety of limiting
behaviours. As an application, we study ``self-normalised'' versions of
$X_t$, that is, $X_t$ after division by $\sup_{0<s\le t}\Delta X_s$,
or by
$\sup_{0<s\le t}\llvert \Delta X_s\rrvert $.
Thus, we obtain necessary and sufficient conditions for $X_t/\sup
_{0<s\le t}\Delta X_s$ and
$X_t/\sup_{0<s\le t}\llvert \Delta X_s\rrvert $ to converge in
probability to~1, or
to $\infty$,
as $t\downarrow0$, so that $X$ is either comparable to, or dominates,
its largest jump.
The former situation tends to occur when the singularity at 0 of the L\'
evy measure
of $X$ is fairly mild (its tail is slowly varying at~0), while the
latter situation
is related to the relative stability or attraction to normality of $X$
at 0 (a steeper singularity at 0).
An important component in the analyses is the way the largest positive
and negative
jumps interact with each other. Analogous ``large time'' (as $t\to
\infty$) versions of the results can also be obtained.
\end{abstract}

% KEYWORDS
% visi is mazosios raides ir pagal abecele
%
\end{frontmatter}

%s1 #&#
\section{Introduction}\label{s1} We study relations between a L\'evy
process $X=(X_t)_{t\ge0}$, its
quadratic variation process $V=(V_t)_{t\ge0}$ and its maximal jump
processes, with
particular interest in how these processes, and how positive and
negative parts of the $X$ process, interact. Representations of
distributions related to these processes are calculated and used as a
basis for making asymptotic (small time) comparisons in their behaviours.

A convenient way of proceeding is to derive identities for the
distributions of $X_t$ modified by subtracting a number of its largest jumps,
%%namely, $\sup_{0<s\le t}\Delta X_s$,
or its jumps of largest modulus, up until time $t$,
joint with $V_t$, modified similarly.
These identities are obtained by considering the Poisson point process
of jumps of $X$, allowing for possible ties in the order statistics of
the jumps.

The distributions
%%of the modified $X_t$, $V_t$, and the large jump processes
thus obtained enable the study of a wide variety of small or large time
kinds of behaviour of $X$. As an application, we investigate
``self-normalised'' versions of $X_t$, giving a comprehensive analysis
of the behaviour of $X_t/\sup_{0<s\le t}\Delta X_s$ and $X_t/\sup
_{0<s\le t}\llvert \Delta X_s\rrvert $ as $t\downarrow0$, and
similarly with $X_t$
replaced by $\llvert X_t\rrvert $.
Two extreme situations are considered; first, when $X$ is of comparable
size to a maximal jump process, for example, $X_t/\sup_{0<s\le
t}\llvert \Delta X_s\rrvert \buildrel P \over\to1$ as $t\downarrow
0$; or,
alternatively, when $X$
dominates a maximal jump process, in the sense that
$X_t/\sup_{0<s\le t}\llvert \Delta X_s\rrvert \buildrel P \over\to
\infty$ as
$t\downarrow0$;
and similarly with $X_t$ replaced by $\llvert X_t\rrvert $, and/or
$\llvert \Delta X_s\rrvert $
replaced by $\Delta X_s$.
Complementary to these is the way the largest positive and negative
jumps interact with each other.

Such results can be seen as continuations in one way or another of a
growing literature in this area which has some classical antecedents.
The original developments occurred in the context of random walks,
where the concept of ``trimming'' by removing extremes from a sample
sum has been studied extensively in the past. Our particular emphasis
on the ratio of the process to its extremes goes back in the random
walk situation to results of Darling \cite{darl} and Arov and Bobrov
\cite{ab}.
Later, Maller and Resnick \cite{MR} gave conditions for a random walk
to be comparable in magnitude to its large values (a heavy-tailed
situation), while Kesten and Maller \cite{KM1992,KM1994}
studied the other end of the spectrum, when the sum dominates its large
values (see Table~1 of \cite{KM1994} for a convenient summary).

Subsequent to these papers there was much development in the general
area of trimmed sums, especially concerning heavy tailed distributions;
see, for example, Cs\"org\H{o}, Haeusler and Mason \cite{CHM}, Berkes
and Horv\'ath \cite{BH}, Berkes, Horv\'ath and Schauer \cite{BHS},
and Griffin and Pruitt \cite{GP1}.
We mention in this context also results of Silvestrov and Teugels \cite
{st} concerning sums and maxima of random walks and triangular arrays, and
Ladoucette and Teugels \cite{LT3} for an insurance application.
There are also recent results about the St. Petersburg game; Gut and
Martin-L\"of \cite{GML} give a ``maxtrimmed'' version of the game,
while Fukker, Gy\"orfi and Kevei \cite{FGK} determine the limit
distribution of the St. Petersburg sum conditioned on its maximum. Cs\"
org\H{o} and Simons \cite{CS} give a review of the later St.
Petersburg literature.

For almost sure versions of particular kinds of sum/max relationships,
see Feller \cite{feller1968}, Kesten and Maller \cite{KM1995} and
Pruitt \cite{pruitt}.

Studies of small time or local behaviour of L\'evy processes go back to
the work of L\'evy and Khintchine \cite{kh,kh2}, in the
1930s. More recent work, relevant to our topic, includes that of Doney
\cite{doney2004}, who gives conditions for a L\'evy process $X$ to
remain positive near 0 with probability approaching 1, and Andrew \cite
{and2}, who similarly analyses the behaviours of the positive and
negative jump processes near 0. There is a connection also with results
of Bertoin \cite{bereg}, who in studying regularity of a L\'evy
process $X$ at 0 was concerned with the dominance of the positive part
of $X$ over its negative part, when $X$ is of bounded variation.
For further background along these lines, we refer to Doney \cite{doney2007}.

Despite all this activity, there seems to have been little done so far
by way of relating the L\'evy process directly to its large jumps, as
we do herein. Of course, our methods rely substantially on previously
developed foundational work. Our representations of the trimmed L\'evy
process, for example, are inspired by those of LePage \cite{LePa,LePb}, LePage, Woodroofe and Zinn \cite{LWZ} and Mori \cite
{mori} for trimmed sums via order statistics,
%(LePage, Woodroofe and Zinn \cite{LWZ} restricted themselves to
%symmetric and continuous distributions.) Our representation
%%for trimming of positive jumps is based on
and Khintchine's \cite{kh} inverse L\'evy measure method. %% to
%represent $\XX_t^+$.
(The corresponding representations are incorporated in our Lemma \ref
{chain}.) In another direction, Rosi\'nski \cite{ros} collects a
number of alternative series representations for L\'evy processes,
especially with a view to simulation of the process.

%Compactness and subsequential versions of many of the results also
%hold.
%%%, and, furthermore, many of the (large time) results when
%transferred to random walks %%are new in that situation too.

Our paper is organised as follows.
The dominance results are in Sections~\ref{s3} and \ref{s4}.
Section~\ref{s+-} compares the positive and negative jump processes.
Before this, in Section~\ref{s2}, we set up notation and, in Theorem
\ref{randrep},
derive the distribution identities using the Poisson point process
structure of the jumps.
Section~\ref{s2} also recalls some basic facts concerning Poisson
point processes and
constructs the distribution of the relevant Poisson random measure from
the jumps of $X$. Particular attention is paid to the possibility of
tied jumps, related to atoms in the canonical measure of~$X$.
We make brief mention of some other possible applications of the
methodology in the final discussion Section~\ref{s8}.

%================================================================
%
% PRELIMINARIES
%==================================================================
%s2 #&#
\section{Distributional representations}\label{s2}
%%==================================================================
Our object of study will be
a real-valued L\'evy process $X=(X_t)_{t\ge0}$ with canonical triplet
$(\gamma,\sigma^2,\Pi)$,
thus having characteristic function $Ee^{\rmi\theta X_t}= e^{t\Psi
(\theta)}$, $t\ge0$,
$\theta\in\myRR$, with characteristic exponent
%
%e2.1 #&#
\begin{equation}
\label{ce} \Psi(\theta):= \rmi\theta\gamma- \frac{1} 2
\sigma^2 \theta^2 +\int_{\myRR_*}
\bigl(e^{\rmi\theta x}-1-\rmi\theta x\bfeins_{\{\llvert x\rrvert \le
1\}} \bigr) \Pi(\rmd x).
\end{equation}
Here, $\gamma\in\myRR$, $\sigma^2\ge0$ and $\Pi$ is a L\'evy
measure on $\myRR$, that is, a Borel measure on $\myRR_*: = \myRR
\setminus\{0\}$
such that $\int_{\myRR_*}(x^2\wedge1)\Pi(\rmd x)<\infty$.
Define measures $\Pi^{(+)}$, $\Pi^{(-)}$, and $\Pi^{\llvert \cdot
\rrvert }$ on
$(0,\infty)$ such that
$\Pi^{(+)}$ is $\Pi$ restricted to $(0,\infty)$,
$\Pi^{(-)}$ is $\Pi(-\cdot)$ restricted to $(0,\infty)$,
and $\Pi^{\llvert \cdot\rrvert }:= \Pi^{(+)} + \Pi^{(-)}$.
The positive, negative and two-sided tails of $\Pi$ are
\begin{eqnarray}\label{pidef}
\pibar^+(x) &:=& \Pi\bigl\{(x,\infty)\bigr\},\qquad
\pibar^-(x) := \Pi\bigl\{(-\infty,-x)\bigr\}\quad\mbox{and}
\nonumber\\[-8pt]\\[-8pt]\nonumber
\pibar(x)&:=&\pibar^+(x)+\pibar^-(x),\qquad x>0.\nonumber
\end{eqnarray}
We are only interested in small time behaviour of $X_t$, so we
eliminate trivial cases by assuming $\pibar(0+)=\infty$ or
$\pibar^+(0+)=\infty$, as appropriate.
Let $\Delta\Pi(y):=\Pi(\{y\})$, $y\in\myRR_*$, and $\Delta
\pibar(y):= \pibar(y-)-\pibar(y)$, $y>0$.
Denote the jump process of $X$ by $(\Delta X_t)_{t\ge0}$, where
$\Delta X_t = X_t-X_{t-}$, $t>0$, with $\Delta X_0 \equiv0$.
The quadratic variation process associated with $X$ is
\[
%%\label{defqv}
V_t:= \sigma^2t+\sum
_{0<s\le t} (\Delta X_s)^2,\qquad t>0,
\]
with $V_0\equiv0$.
%% We consider also the maximal jump process $\sup_{0<s\le t}\Delta X_s$
%%and associated variant processes $\sup_{0<s\le t}\left\vert\Delta X_s
%\right\vert$,
%%$\sup_{0<s\le t}\Delta X_s^+$, and $\sup_{0<s\le t}\Delta X_s^-$,
%where
%%$\Delta X_t^+=\max(0, \Delta X_t)$ and $\Delta X_t^-=\Delta X_t^+-
%\Delta X_t$, $t>0$.
Recall that $X$ is of bounded variation if $\sum_{0<s\le t}\llvert
\Delta X_s\rrvert <\infty$
a.s. for all $t>0$, equivalently, if $\sigma^2=0$ and
$\int_{\llvert x\rrvert \le1}\llvert x\rrvert \Pi(\rmd x)<\infty$.
If this is the case, \eqref{ce} takes the form
\[
\rmi\theta\,\rmd_X + \int_{\myRR}
\bigl(e^{\rmi\theta x}-1\bigr)\Pi(\rmd x),
\]
where $\rmd_X$ is the drift of $X$.

%=======================================================================
% Representation of the distribution
%
%%%%========================================================================

In deriving representations for the joint distributions of $X_t$, $V_t$
and the $r$th maximal jump processes, it is convenient to work with
the processes having the $r$ largest jumps, or the $r$ jumps largest in
modulus, subtracted. These ``trimmed'' processes are no longer L\'evy
processes, but we can give useful representations for their marginal
distributions.
The expressions are in terms
of a truncated L\'evy process, together with one or two Poisson
processes, and a
Gamma %%sum of $r$ exponential
random variable, all processes and random variables independent of one another.

%Introduce the notation $\Delta X_t^{(1)} = \sup_{0<s\le t}\Delta X_s$,
% and let $\wt{\Delta X}_t^{(1)}$
%be any term of maximum modulus among $(\Delta X_s)_{0<s\le t}$, that
%is,
%any $\Delta X_y$ with $y\in(0,t]$ such that $\sup_{0<s\le t}\left\vert
%\Delta
%X_s\left\vert\le\left\vert\Delta X_y\left\vert$.
For any integer $r=1,2,\ldots,$ let $\Delta X_t^{(r)}$ and $\widetilde
{\Delta X}_t^{(r)}$ be the
$r$th largest positive jump and the $r$th largest jump in modulus
up to time $t$, respectively.
Formal definitions of these, allowing for the possibility of tied
values (we choose the order uniformly among the ties), are given in
Section~\ref{sub21} below.
``One-sided'' and ``modulus'' trimmed versions of $X$ are then defined as
%
%
%e2.2 #&#
\begin{equation}
\label{trims} {}^{(r)}X_t:= X_t-\sum
_{i=1}^r {\Delta X}_t^{(i)}
\quad\mbox{and}\quad{}^{(r)}\widetilde X_t:=
X_t- \sum_{i=1}^r \widetilde{
\Delta X}_t^{(i)},
\end{equation}
with corresponding trimmed quadratic variation processes
\[
%%\label{vtl}
^{(r)}{ V}_t:=V_t- \sum
_{i=1}^r \bigl({\Delta X}_t^{(i)}
\bigr)^2 \quad\mbox{and}\quad{}^{(r)}{\widetilde
V}_t:=V_t-\sum_{i=1}^r
\bigl( \widetilde{\Delta X}_t^{(i)}\bigr)^2,\qquad
t>0.
\]

Recall the definitions of the tails of $\Pi$ in \eqref{pidef}. Let
\[
%%\label{linv}
\pibar^\leftarrow(x)=\inf\bigl\{y>0: \pibar(y) \le x\bigr\},\qquad x>0,
\]
be the right-continuous inverse of the nonincreasing function $\pibar$,
and similarly for $\pibar^{+, \leftarrow}$ and $\pibar^{-,\leftarrow}$.
By convention, the inf of the empty set is taken as $\infty$. The
following properties of the inverse function will be used frequently
(see Resnick \cite{res87}, Section~0.2).
For each $x, y > 0$, $\pibar^{\leftarrow}(x)\le y$ if and only if
$\pibar(y) \le x$;
%%$\pibar^{\leftarrow}(x) > y$ if and only if $\pibar(y) > x$;
$\pibar(\pibar^{\leftarrow}(x)) \le x \le\pibar(\pibar
^{\leftarrow}(x)-)$;
and $ \overline{\Pi}^{\leftarrow}(\pibar(x))\le x$;
%%$\pibar(\pibar^{\leftarrow}(\pibar(x))) = \pibar(x)$.
similarly, for $\pibar^{\pm}$.
We refer to Appendix A in Fan \cite{fanthesis} for more details.

%%a representation formula for a L\'evy process with the maximum jump
%or the maximum modulus %%jump removed. Its proof is based on random
%walk results presented in Section \ref{s1-3}.

We introduce four families of processes, indexed by $v>0$, truncating
jumps from sample paths of $X_t$ and $V_t$, respectively.
Let $v,t>0$. When $\pibar(0+)=\infty$, we set
%
%
%e2.3 #&#
\begin{eqnarray}\label{21b}
\widetilde X_t^v&:=&X_t-\sum
_{0<s\le t} \Delta X_s \bfeins_{\{\llvert \Delta X_s\rrvert \ge
\overline{\Pi}^{\leftarrow}(v)\}}
\quad\mbox{and}
\nonumber\\[-8pt]\\[-8pt]\nonumber
\widetilde V_t^v&:=&V_t-
\sum_{0<s\le t} (\Delta X_s)^2
\bfeins_{\{\llvert \Delta X_s\rrvert \ge\overline{\Pi}^{\leftarrow}(v)\}}.
\end{eqnarray}
When $\pibar^+(0+)=\infty$, we set
\[
%%\label{22b}
X_t^v:=X_t-\sum
_{0<s\le t} \Delta X_s \bfeins_{\{\Delta X_s\ge
\pibar^{+,\leftarrow}(v)\}}\quad
\mbox{and}\quad V_t^v:=V_t-\sum
_{0<s\le t} (\Delta X_s)^2
\bfeins_{\{\Delta X_s\ge
\pibar^{+,\leftarrow}(v)\}}.
\]

Under the assumptions $\pibar(0+)=\infty$ and $\pibar^+(0+)=\infty
$, $(\widetilde X_t^v)_{t\ge0}$ and $(X_t^v)_{t\ge0}$ are well-defined
L\'
evy processes with canonical triplets, respectively,
%
%
%e2.4 #&#
\begin{equation}
\label{trip2} \biggl(\gamma-\bfeins_{\{\overline{\Pi}^{\leftarrow}(v)\le
1\}
}\int_{\overline{\Pi}^{\leftarrow}(v)\le
\llvert x\rrvert \le1}x
\Pi(\rmd x), \sigma^2, \Pi(\rmd x)\bfeins_{\{
\llvert x\rrvert <\overline{\Pi}^{\leftarrow}(v)\}} \biggr)
\end{equation}
and
%
%e2.5 #&#
\begin{equation}
\label{trip1} \biggl(\gamma-\bfeins_{\{\overline{\Pi}^{+,\leftarrow
}(v)\le1\}
}\int_{\overline{\Pi}^{+,\leftarrow}
(v)\le x \le1}x
\Pi(\rmd x), \sigma^2, \Pi(\rmd x)\bfeins_{\{
x<\overline{\Pi}^{+,\leftarrow}(v)\}} \biggr).
\end{equation}

Our main result in this section gives very general
representations for the joint distributions of
$({}^{(r)}\widetilde X_t, {}^{(r)}\widetilde V_t, \llvert \widetilde
{\Delta X}_t^{(r)}\rrvert )$
and of
$({}^{(r)}X_t, ^{(r)}V_t, \Delta X_t^{(r)})$,
allowing for possible tied values in the large jumps.
We make the convention throughout that a Poisson random variable with
parameter $0$ is $0$. Note that then the expressions in \eqref{Ppar},
\eqref{cor1aa} and \eqref{cor1bb} below are zero when $\Pi$ has no
atoms. But we do not assume this.

%
%th2.1 #&#
\begin{teo}\label{randrep}
Let $r\in\NN=\{1,2,3,\ldots\}$ and $\SSSS_r$ be a $\operatorname{Gamma}(r,1)$
random variable.
Suppose $Y^\pm=(Y_t^\pm)_{t\ge0}$ and $Y=(Y_t)_{t\ge0}$
are independent Poisson processes with $EY_1^\pm=EY_1=1$.
Assume that $X$, $\SSSS_r$, $Y^+$, $Y^-$, and $Y$ are independent as
random elements.

\begin{longlist}[(ii)]
\item[(i)] Assume $\pibar(0+)=\infty$. For each $v>0$, let
%
%e2.6 #&#
\begin{equation}
\label{Ppar} \kappa^\pm(v):= \bigl(\pibar\bigl(\overline{
\Pi}^{\leftarrow
} (v )- \bigr)-v \bigr)\frac{\Delta\Pi( \pm\overline{\Pi}^{\leftarrow
} (v )
)}{\Delta\pibar( \overline{\Pi}^{\leftarrow}
(v )
)}\bfeins_{\{\Delta\pibar(\overline{\Pi}^{\leftarrow}
(v )
)\neq0\}}
\end{equation}
and for $v>0$, $t>0$, set
%
%e2.7 #&#
\begin{equation}
\label{cor1aa} \widetilde{G}_t^{v}:=\overline{
\Pi}^{\leftarrow}(v) \bigl(Y^+_{t\kappa
^+(v)} - Y^-_{t\kappa^-(v)}\bigr) \quad
\mbox{and}\quad\widetilde{H}_t^{v}:= \bigl(\overline{
\Pi}^{\leftarrow}(v)\bigr)^2\bigl(Y^+_{t\kappa^+(v)} +
Y^-_{t\kappa^-(v)}\bigr).
\end{equation}
Then, for each $t>0$, we have
%
%e2.8 #&#
\begin{equation}
\label{randdis1} \bigl({}^{(r)}\widetilde X_t,
{}^{(r)}\widetilde V_t, \bigl\llvert\widetilde{\Delta
X}_t^{(r)}\bigr\rrvert\bigr) \stackrel{\mathrm{D}} {=}
\bigl(\widetilde{ X}_t^{v} + \widetilde
{G}_t^{v}, \widetilde{ V}_t^{v} +
\widetilde{H}_t^{v}, \overline{\Pi}^{\leftarrow} (v )
\bigr) \mid_{v=\SSSS_r/t}.
\end{equation}

\item[(ii)] Assume $\pibar^+(0+)=\infty$.
%% and $(Y_t)_{t\ge0}$ is a Poisson process with $EY_1=1$, independent
%of $\SSSS_i$.
For each $v>0$, let $\kappa(v):=\pibar^+(\overline{\Pi
}^{+,\leftarrow}(v)-)-v$,
and for $v>0$, $t>0$, set
%
%e2.9 #&#
\begin{equation}
\label{cor1bb} G_t^{v}:=\overline{\Pi}^{+,\leftarrow}(v)Y_{t\kappa(v)}
\quad\mbox{and}\quad H_t^{v}:= \bigl(\overline{
\Pi}^{+,\leftarrow}(v)\bigr)^2Y_{t\kappa(v)}.
\end{equation}
Then, for each $t>0$, we have
%
%
%e2.10 #&#
\begin{equation}
\label{2rrep1} \bigl({}^{(r)}X_t, ^{(r)}V_t,
\Delta X_t^{(r)} \bigr) \stackrel{\mathrm{D}} {=}
\bigl(X_t^{v} + G_t^v,
{V}_t^{v} + H_t^v, \overline{
\Pi}^{+,\leftarrow} (v ) \bigr) \mid_{v=\SSSS_r/t}.
\end{equation}
\end{longlist}
\end{teo}

%
%re2.1 #&#
\begin{rem}\label{rem0}
Processes ${}^{(r)}\widetilde X_t$ and ${}^{(r)}X_t$ are not L\'
evy processes;
their increments are not independent, or homogeneous in distribution.
But the identities \eqref{randdis1} and \eqref{2rrep1} express
their marginal distributions in terms of distributions of L\'evy
processes, mixed in a sense according to their $r$th largest jumps,
with allowance made for ties.
This opens the possibility for results obtained from analyses of the
underlying L\'evy processes to be transferred to the trimmed processes.
We exemplify this procedure in a variety of ways in Sections~\ref{s3}
and \ref{s4}.
\end{rem}

As an immediate corollary of Theorem \ref{randrep}, the following identities
will be useful.

%
%co1 #&#
\begin{corollary}
Using the notation in Theorem \ref{randrep}, we have,
%%, for the special case when $r = 1$, we obtain the following results
%regarding the maximal jump (and in modulus) process.
for $x\in\myRR$, $y\ge0$, $t>0$, $r=1,2,\ldots:$

\begin{longlist}[(ii)]
\item[(i)] when $\pibar(0+)=\infty$,
%
%e2.11 #&#
\begin{eqnarray}
\label{cor1a} &&P \bigl({}^{(r)}\widetilde X_t\le x\bigl
\llvert\widetilde{\Delta X}_t^{(r)}\bigr\rrvert,
{}^{(r)}\widetilde V_t\le y\bigl\llvert\widetilde{\Delta
X}_t^{(r)}\bigr\rrvert^2 \bigr)
\nonumber\\[-8pt]\\[-8pt]\nonumber
&&\quad =\int_0^\infty P \bigl(\widetilde{
X}_t^{v} + \widetilde{G}_t^{v} \le
x \overline{\Pi}^{\leftarrow}(v), \widetilde{ V}_t^{v} +
\widetilde{H}_t^{v}\le y\bigl(\overline{\Pi
}^{\leftarrow}(v)\bigr)^2 \bigr)P(\SSSS_r\in t\,\rmd v);
\end{eqnarray}
\item[(ii)] when $\pibar^+(0+)=\infty$,
%
%e2.12 #&#
\begin{eqnarray}
\label{cor1b} &&P \bigl({}^{(r)} X_t\le x \Delta
X_t^{(r)}, {}^{(r)}V_t\le y \bigl(
\Delta X_t^{(r)}\bigr)^2 \bigr)
\nonumber\\[-8pt]\\[-8pt]\nonumber
&&\quad =\int_0^\infty P \bigl(X_t^{v}
+ G_t^{v} \le x \overline{\Pi}^{+,\leftarrow}(v),
V_t^{v} + H_t^{v}\le y\bigl(
\overline{\Pi}^{+,\leftarrow}(v)\bigr)^2 \bigr)P(\SSSS_r
\in t\,\rmd v).
\end{eqnarray}
%
%%where $P(\EEEE\in t\rmd v)=te^{-tv}\,\rmd v$, $v\ge%%0$.
\end{longlist}
\end{corollary}

%================REMARK==========================

%\begin{remark}\label{rem1} {\rm
%Although not explicit in the proof of Theorem \ref{randrep}, the
%extra Poisson processes in \eqref{randdis1} and \eqref{2rrep1}
%originally derive from the possibility of a binomial number of ties in
%the maximum of a random walk approximation to $X_t$, which are
%transmuted to Poisson random variables in the continuous time process.
%When $\Pi(\cdot)$ is a continuous function, there are no such ties
%a.s., and the Poisson random variables disappear from
%\eqref{randdis1} and \eqref{2rrep1}.
%%We expand on this further in Remark \ref{rem2} at the end of the
%proof of Theorem \ref{randrep}
%%and in Theorem \ref{tdel}.}
%}
%\end{remark}
%
%{\color{red}A basic tool to be used is a version of a compound
%Poisson approximation to $X$, which is reviewed in Subsection
%\ref{sub21}. The proof of Theorem \ref{randrep} then proceeds by
%calculating, in Subsection \ref{sub22}, corresponding representations
%for trimmed random walks. These are then transferred to the L\'evy
%process
%by taking limits. Subsection \ref{sub23} ties these ideas together.}

In proving Theorem \ref{randrep}, we make use of the underlying
Poisson point process (PPP) structure of the jumps of a L\'evy process.
We begin in Section~\ref{sub21} with a precise definition of the
order statistics of a PPP when tied values may be present. In
Section~\ref{sub22}, we review basic properties of standard PPPs and
in Section~\ref{sub23} construct the distribution of a Poisson random
measure (PRM) from the jumps of a L\'{e}vy process
through a series of marking and deterministic transformations. Also, in
Section~\ref{sub23}, we derive the joint distribution of the trimmed
point process using the point process order statistics. This machinery
allows us to complete the proof of Theorem \ref{randrep} in
Section~\ref{sub24}.

%===============================================================
% subsection 2.1: ordered statistics with ties
%===============================================================
%s2.1 #&#
\subsection{Order statistics with ties}
\label{sub21}
Introduce $\XX$ as the point measure associated with the jumps of $X$:
\[
%% \label{defX}
\XX=\sum_{s}\delta_{(s,\Delta X_s)}.
\]
$\XX$ is a Poisson point process\footnote{For necessary material on
point processes, we refer to Chapter~12 in Kallenberg \cite{Kall} or
Chapter~5 in Resnick \cite{res07}.} (PPP) on $[0,\infty)\times\myRR
_*$ with intensity measure $\rmd s\otimes\Pi(\rmd x)$.
Analogously, the PPPs of positive and negative jumps and jumps in
modulus associated with $\XX$ are
\begin{eqnarray*}
\XX^+&=&\sum_{s}\bfeins_{(0,\infty)}(\Delta
X_s)\delta_{(s,\Delta
X_s)},%%\label{defX+}
\qquad
\XX^- =\sum_{s} \bfeins_{(0,\infty)}(-\Delta
X_s)\delta_{(s,-\Delta X_s)},
\\
\XX^{\llvert \cdot\rrvert }&=&\XX^++\XX^-=\sum_{s}
\delta_{(s,\llvert \Delta X_s\rrvert )},
\end{eqnarray*}
having intensity measures $\rmd s\otimes\Pi^{\pm,\llvert \cdot\rrvert
}(\rmd x)$, respectively. For $t>0$,
we consider restrictions of these processes to the time interval $[0,t]$ by introducing
\[
\XX_t(\cdot):=\XX\bigl([0,t]\times\myRR_*\cap\cdot\bigr) \quad
\mbox{and}\quad\XX_t^{\pm,\llvert \cdot\rrvert }(\cdot)=\XX^{\pm
,\llvert \cdot\rrvert }
\bigl([0,t]\times(0,\infty)\cap\cdot\bigr).
\]
Assume $\pibar(0+)=\infty$ and $t>0$. Our first task is to specify
the points with maximum modulus in~$\XX_t$.

%%%%--------------------------------------------------------------------------------------------------------------------

Let $\widetilde T^{(1)}(\XX_t)$ be randomly chosen,
independently of $(X_t)_{t\ge0}$,
according to the discrete uniform distribution in the set $\{0\le s\le
t: \llvert {\Delta X}_s\rrvert =\sup_{0\le u\le t}\llvert {\Delta
X}_u\rrvert \}$, which is almost
surely finite. Then define
$\widetilde{\Delta X}_t^{(1)} = \widetilde{\Delta X}^{(1)}(\XX_t): =
\Delta X_{\widetilde
T^{(1)}(\XX_t)}$. Define the maximum modulus trimmed point process on
$[0,t] \times\myRR_*$ by
\[
{}^{(1)}\widetilde\XX_t: = \XX_t -
\delta_{(\widetilde T^{(1)}(\XX
_t), \widetilde
{\Delta X}_t^{(1)})}.
\]
Let $r = 2,3, \ldots.$ Iteratively, we define $ \widetilde T^{(r)}(\XX
_t): =
\widetilde T^{(1)}(^{(r-1)}\widetilde\XX_t)$ and $\widetilde{\Delta
X}^{(r)}_t: = \Delta
X_{\widetilde T^{(r)}(\XX_t)}$.
The $r$-\textit{fold modulus trimmed point process of modulus jumps} is then
defined by
\[
%\wt T^{(r)}(\XX_t): = \wt T^{(1)}(^{(r-1)}\XX_t),\quad
%\wt{\Delta X}^{(r)}_t: = \Delta X_{\wt T^{(r)}(\XX_t)} \quad
%\mbox{and}\quad
{}^{(r)}\widetilde\XX_t:=
\XX_t - \sum_{i=1}^r \delta
_{(\widetilde T^{(i)}(\XX
_t), \widetilde{\Delta X}_t^{(i)})}.
\]

In a similar way, under the assumption $\pibar^+(0+)=\infty$, we can
define the ordered pairs
\[
\bigl(T^{(1)}\bigl(\XX_t^+\bigr),\Delta
X^{(1)}_t\bigr), \bigl(T^{(2)}\bigl(
\XX_t^+\bigr),\Delta X^{(2)}_t\bigr),
\bigl(T^{(3)}\bigl(\XX_t^+\bigr),\Delta X^{(3)}_t
\bigr),\dots \in[0,t]\times(0,\infty),
\]
such that $\Delta X^{(1)}_t\ge\cdots\ge\Delta X^{(r)}_t$ are the
$r$th largest order statistics of positive jumps of $X$ sampled on
time interval $[0,t]$. By subtracting the points corresponding to large
jumps, analogously as we did for
${}^{(r)}\widetilde\XX_t$, we then define the $r$-\textit{positive trimmed
point process of positive jumps} by
\[
{}^{(r)}\XX^+_t:=\XX^+_t-\sum
_{1\le i\le r}\delta_{(T^{(i)}(\XX
_t^+), \Delta X_t^{(i)})}.
\]

%================================================================
%
% Subsection 2.2: Homogenous Poisson point process
%==================================================================
%s2.2 #&#
\subsection{Standard Poisson point process}\label{sub22}
In this section, we provide alternative constructions of $\XX
_t,{}^{(r)}\widetilde\XX_t,\XX^+_t,{}^{(r)}\XX^+_t$, this time starting
from homogeneous processes.

Let $(\UUUU_i)$, $(\UUUU_i')$ and $(\EEEE_i)$ be independent, where
$(\UUUU_i)$ and $(\UUUU_i')$ are i.i.d. sequences of uniformly distributed
random variables in $(0,1)$, and $(\EEEE_i)$
is an i.i.d. sequence of exponentially distributed
random variables with common parameter $E\EEEE_i=1$.
Then $\SSSS_r=\sum_{i=1}^r\EEEE_i$ is a Gamma$(r,1)$ random
variable, $r\in\NN$.

For $t>0$, we introduce
\[
%% \label{VVorder}
\VV_t:=\sum_{i\ge1}
\delta_{(t\UUUU_i, \SSSS_i/t)}\quad\mbox{and}\quad\VV'_t:=\sum
_{i\ge1}\delta_{(t\UUUU_i, \UUUU_i', \SSSS_i/t)}.
\]
Then $\VV_t$ and $\VV_t'$ are homogeneous PPPs on $[0,t]\times
(0,\infty)$ and
$[0,t]\times(0,1)\times(0,\infty)$
with intensity measures $\rmd s\otimes\rmd v$ and $\rmd s\otimes\rmd
u' \otimes\rmd v$, respectively.
For $r\in\NN:=\{0,1,2,\ldots\}$, we define their $r$-fold trimmed
counterparts by
\[
%%\label{ktrim-def}
{}^{(r)}\VV_t:=\sum
_{i>r}\delta_{(t\UUUU_i,\SSSS_i/t)}\quad\mbox{and}\quad
{}^{(r)}\VV'_t:=\sum
_{i>r}\delta_{(t\UUUU_i, \UUUU_i',\SSSS
_i/t)}.
\]

When $\pibar(0+)=\infty$, we consider the transformation
\[
\bigl(I,I,\overline{\Pi}^{\leftarrow}\bigr):[0,t]\times(0,1)\times
(0,\infty)
\to[0,t]\times(0,1)\times(0,\infty),\qquad\bigl(s,u',v\bigr)\mapsto
\bigl(s,u',\overline{\Pi}^{\leftarrow}(v)\bigr).
\]
Still assuming $\pibar(0+)=\infty$, by the Radon--Nikodym theorem,
there exist Borelian functions $g^{\pm}:(0,\infty)\to(0,\infty)$
with $g^++g^-\equiv1$ such that $\rmd\Pi^{\pm}= g^{\pm}\,\rmd\Pi
^{\llvert \cdot\rrvert }$ and, in particular,
%
%
%e2.13 #&#
\begin{equation}
\label{g+-} \pibar^{\pm}(x)=\int_{(x,\infty)}
g^{\pm}(y)\Pi^{\llvert \cdot\rrvert }(\rmd y),\qquad x>0.
\end{equation}
We use $g^+$ to return the sign to the process by a second
transformation $m:[0,t]\times(0,1)\times(0,\infty)\to[0,t]\times
\myRR_*$, defined by
%
%
%e2.14 #&#
\begin{equation}
\label{defm} m\bigl(s,u',x\bigr):= \cases{ (s,x),&\quad if
$u'<g^+(x)$,
\cr
(s,-x),&\quad if $u'\ge g^+(x)$.}
\end{equation}
In summary, let ${\VV'}_t^{m\circ(I,I,\overline{\Pi}^{\leftarrow
})}$ be the point
process on $[0,t]\times\myRR^*$, being the image of the composition of
the above transformations applied to $\VV_t'$:
\begin{eqnarray*}
%% \label{chain*}
\VV'_t&\stackrel{\bigl(I,I,\overline{
\Pi}^{\leftarrow}\bigr)}\longrightarrow& {\VV'}^{(I,I,\overline{\Pi
}^{\leftarrow})}_t:=
\sum_{i\ge1}\delta_{(t\UUUU_i, \UUUU_i', \overline{\Pi
}^{\leftarrow}(\SSSS
_i/t))}
\nonumber
\\
&\stackrel{m}\longrightarrow& {\VV'}^{m\circ(I,I,\overline{\Pi
}^{\leftarrow})}_t:=
\sum_{i\ge1}\delta_{m(t\UUUU_i,\UUUU_i',\overline{\Pi
}^{\leftarrow}(\SSSS_i/t))}.
\end{eqnarray*}
Their trimmed counterparts are similarly defined by setting, for $r\in
\NN$,
\begin{eqnarray*}
%% \label{trchain*}
{}^{(r)}\VV'_t&\stackrel{\bigl(I,I,
\overline{\Pi}^{\leftarrow
}\bigr)}\longrightarrow& {}^{(r)} {
\VV'}^{(I,I,\overline{\Pi}^{\leftarrow})}_t:= \sum
_{i>r}\delta_{(t\UUUU_i, \UUUU_i', \overline{\Pi
}^{\leftarrow}(\SSSS
_i/t))}
\nonumber
\\
&\stackrel{m}\longrightarrow& {}^{(r)}{\VV'}^{m\circ(I,I,\overline{\Pi
}^{\leftarrow})}_t:=
\sum_{i>r}\delta_{m(t\UUUU_i, \UUUU_i',\overline{\Pi
}^{\leftarrow}(\SSSS_i/t))}.
\end{eqnarray*}

When $\pibar^+(0+)=\infty$ we can contrive $\overline{\Pi
}^{+,\leftarrow}$ as a transformation of $(0,\infty)$ into
$(0,\infty)$ and we will consider the image measures of
$\VV_t$ and ${}^{(r)}\VV_t$ under $(I,\overline{\Pi}^{+,\leftarrow
}):[0,t]\times(0,\infty)\to[0,\infty)\times(0,\infty)$, defined by
\[
\label{mappedplusktrim-def} \VV_t^{(I,\overline{\Pi}^{+,\leftarrow})}:=
\sum
_{i\ge1}\delta_{(t\UUUU_i,\overline{\Pi}^{+,\leftarrow
}(\SSSS_i/t))}\quad\mbox{and}\quad
{}^{(r)}\VV_t^{(I,\overline{\Pi}^{+,\leftarrow})}:= \sum
_{i>r}\delta_{(t\UUUU_i,\overline{\Pi}^{+,\leftarrow}(\SSSS
_i/t))}.
\nonumber
\]

%===============================================================
% subsection: r-trimmed PPP
%===============================================================

%s2.3 #&#
\subsection{Representations for $r$-trimmed PPPs}\label{sub23}
In this section, the original point process $\XX$, its ordered jumps,
and the trimmed point process, is related to a corresponding standard
version $\VV$. %%We have the following results.
%
%le1 #&#
\begin{lemma}\label{chain} Let $t>0$ and $r\in\NN$.

\begin{longlist}[(ii)]
\item[(i)] If $\pibar(0+)=\infty$, we have the following
distributional equivalences:
%
%
%e2.15 #&#
%e2.16 #&#
\begin{eqnarray}
&& \XX_t \stackrel{\mathrm{D}} {=} {\VV'}^{m\circ(I,I,\overline{\Pi
}^{\leftarrow})}_t,\label{chbasic}
\\
&& \bigl(\widetilde T^{(i)}(\XX_t),\widetilde{\Delta
X}_t^{(i)} \bigr)_{i\ge1} \stackrel{\mathrm{D}} {=}
\bigl(m \bigl(t\UUUU_i,\UUUU_i',
\overline{\Pi}^{\leftarrow}(\SSSS_i/t) \bigr)
\bigr)_{i\ge1},\label{chpair01}
\\
&& \bigl\{ \bigl(\widetilde T^{(i)}(\XX_t),\widetilde{\Delta
X}_t^{(i)} \bigr)_{1\le i\le r},{}^{(r)}
\widetilde\XX_t \bigr\}
\nonumber\\[-8pt]\label{chpair1} \\[-8pt]\nonumber
&&\quad  \stackrel{\mathrm{D}} {=} \bigl\{ \bigl(m
\bigl(t\UUUU_i,\UUUU_i',\overline{
\Pi}^{\leftarrow
}(\SSSS_i/t)\bigr) \bigr)_{1\le i\le r},{}^{(r)}{
\VV'}^{m\circ(I,I,\overline{\Pi
}^{\leftarrow})}_t \bigr\}.
\nonumber
\end{eqnarray}

\item[(ii)] If $\pibar^+(0+)=\infty$, we have the following
distributional equivalences:
\begin{eqnarray*}
\XX_t^+&\stackrel{\mathrm{D}} {=}& \VV^{(I,\overline{\Pi
}^{+,\leftarrow})}_t,
%%\label{chbasic+}
\\
\bigl(T^{(i)}\bigl(\XX_t^+\bigr),{\Delta
X}_t^{(i)} \bigr)_{i\ge1}&\stackrel{\mathrm{D}}
{=}& \bigl(t\UUUU_i,\overline{\Pi}^{+,\leftarrow}(
\SSSS_i/t) \bigr)_{i\ge1},%%\label{chpair01+}
\\
\bigl\{ \bigl(T^{(i)}\bigl(\XX_t^+\bigr),\Delta
X_t^{(i)} \bigr)_{1\le i\le
r},{}^{(r)}
\XX^+_t \bigr\}&\stackrel{\mathrm{D}} {=}& \bigl\{ \bigl(t
\UUUU_i,\overline{\Pi}^{+,\leftarrow}(\SSSS_i/t)
\bigr)_{1\le i\le r},{}^{(r)}\VV^{(I,\overline{\Pi
}^{+,\leftarrow})}_t \bigr\}.
%%label{ch:pair1+}
\end{eqnarray*}
\end{longlist}
\end{lemma}

%%\verb\begin{proof} [of Lemma \rm{\ref{chain}}] \verb*

\begin{pf} %{Proof of Lemma {\ref{chain}}}
(i)~Assume $\pibar
(0+)=\infty$, and introduce
\[
%%\label{deftildem}
\widetilde m:(0,1)\times(0,\infty)\to\myRR_*,\qquad\widetilde m
\bigl(u',x\bigr):= x\bfeins_{u'<g^+(x)} -x
\bfeins_{u'\ge g^+(x)}.
\]
(The\vspace*{1pt} mapping $\widetilde m$ is the same as the $m$ in \eqref{defm} without
the time component.)

Let $\mu^T:=\mu\circ T^{-1}$ denote the image measure of a measure
$\mu$ under a transformation $T$. Using this notation, and in view
of~\eqref{g+-}, we get
from $(\rmd v)^{\overline{\Pi}^{\leftarrow}}= \rmd\Pi^{\llvert \cdot
\rrvert }$ that
%
%
%e2.17 #&#
\begin{eqnarray}\label{UUU}
\bigl(\rmd u'\otimes\rmd v\bigr)^{\widetilde m\circ(I,\overline{\Pi
}^{\leftarrow})}\bigl((x,\infty)
\bigr) &=&\bigl(\rmd u'\otimes\rmd\Pi^{\llvert \cdot\rrvert }\bigr)
\bigl(
\widetilde m^{-1} \bigl((x,\infty) \bigr) \bigr)
\nonumber
\\
&=& \int_{(x,\infty)} g^+(v) \Pi^{\llvert \cdot\rrvert }(\rmd v)
\\
&=&\pibar^+(x),\qquad x>0,\nonumber
\end{eqnarray}
and similarly with $(x,\infty)$ replaced by
$(-\infty, -x)$, and $g^+$, $\pibar^+$, replaced by
$g^-$, $\pibar^-$.
With $m$ as in \eqref{defm}, and since the tail functions determine
the corresponding measures, \eqref{UUU} extends to
%
%
%e2.18 #&#
\begin{equation}
\label{imageM2} \bigl(\rmd s\otimes\rmd u' \otimes\rmd v
\bigr)^{m\circ(I, I, \overline{\Pi
}^{\leftarrow}
)}=\rmd s\otimes\rmd\Pi.
\end{equation}
Let $h:= m\circ(I, I, \overline{\Pi}^{\leftarrow})$.
It follows from \eqref{imageM2} that $\XX_t$ and ${\VV'}^{m\circ
(I, I, \overline{\Pi}^{\leftarrow})}_t={\VV'}^h_t$ share a common
intensity measure
$\rmd s\otimes\rmd\Pi$. Since both $\XX$ and ${\VV'}^h$ are simple
PPPs, this completes the proof of~\eqref{chbasic}.

In order to show~\eqref{chpair01}, introduce record times, defined
recursively by
\[
R_n:=\min\bigl\{i>R_{n-1}:\overline{\Pi}^{\leftarrow}(
\SSSS_i/t)>\overline{\Pi}^{\leftarrow}(\SSSS_{R_{n-1}}/t)
\bigr\},\qquad R_1:=1, n= 2, 3, 4,\ldots.
\]
Observe that $(R_n)_{n\ge1}$ is independent of $(\UUUU_i)$ and
$(\UUUU'_i)$.

Construct the sequence $(\widetilde T^{(i)}({\VV'}^h_t),\widetilde
{\Delta
X}^{(i)}({\VV'}^h_t))_{i\ge1}$ associated with trimming the process
${\VV'}_t^h$ by choosing a sequence of independent permutations
$(\sigma_n)_{n\ge1}$, where
\[
\sigma_n:\{R_{n-1},\dots,R_n-1\}\stackrel{1:1}
\longrightarrow\{ R_{n-1},\dots,R_n-1\},\qquad n=2,3, 4,
\ldots,
\]
are chosen according to the discrete
uniform distribution amongst the finitely many candidates,
independently of $(X_t)_{t\ge0}$. By our construction of trimming,
the pairs $(\{R_n\},\{\sigma_n\})$ and $(\{\UUUU_i\},\{\UUUU'_i\})$
are also independent. Consequently,
\begin{eqnarray*}
&& \bigl\{\widetilde T^{(i)}\bigl({\VV'}^h_t
\bigr),\widetilde{\Delta X}^{(i)}\bigl({\VV'}^h_t
\bigr) \bigr\}_{i\ge1}
\\
&&\quad= \bigl\{ \bigl( \widetilde T^{(i)}\bigl({\VV'}^h_t
\bigr),\widetilde{\Delta X}^{(i)}\bigl({\VV'}^h_t
\bigr) \bigr)_{R_{n-1}\le i<R_n} \bigr\}_{n\ge2}
\\
&&\quad= \bigl\{m \bigl(t\UUUU_{\sigma_n(i)},\UUUU'_{\sigma
_n(i)},
\overline{\Pi}^{\leftarrow}(\SSSS_{R_{n-1}}/t) \bigr)_{R_{n-1}\le i<R_n}
\bigr\}_{n\ge2}
\\
&&\quad\stackrel{\mathrm{D}} {=} \bigl\{m \bigl(t\UUUU_i,\UUUU
'_i,\overline{\Pi}^{\leftarrow}(\SSSS
_{R_{n-1}}/t) \bigr)_{R_{n-1}\le i<R_n} \bigr\}_{n\ge2}
\\
&&\quad= \bigl\{m\bigl(t\UUUU_i, \UUUU'_i,
\overline{\Pi}^{\leftarrow}(\SSSS_{i}/t)\bigr) \bigr
\}_{i\ge
1}.
\end{eqnarray*}
In view of~\eqref{chbasic}, this completes the proof of~\eqref{chpair01}.
Note that~\eqref{chpair1} follows from~\eqref{chpair01}. Part~(ii)
is shown analogously.
\end{pf}

Next is our main theorem giving the representation for trimmed PPPs.
For $x>0$, write $\XX_t^{+\cdot<x}$
and $\XX_t^{\llvert \cdot\rrvert <x}$ for point processes generated
by deleting all
points in $\XX^+_t$ and $\XX_t$ not lying in the regions $[0,t]\times
(0,x)$ and
$[0,t]\times(-x,x)_*$, respectively:
\[
\XX_t^{+\cdot< x}(\cdot):=\XX^+\bigl([0,t]\times(0,x)\cap\cdot
\bigr)
\]
and
\[
\XX_t^{\llvert \cdot\rrvert <x}(\cdot):=\XX
\bigl([0,t]
\times(-x,x)_*\cap\cdot\bigr).
\]
%
%
%th1 #&#
\begin{theorem}\label{orderMod}
Assume that $\XX$, $(\UUUU_i)$, $(\UUUU'_{i})$, $\SSSS_r$,
$Y^{\pm}=(Y^{\pm}(t))_{t\ge0}$, $Y=(Y(t))_{t\ge0}$, are independent
processes, with $Y^\pm$ and $Y$ being standard Poisson processes.

\begin{longlist}[(ii)]
\item[(i)] Assume $\pibar(0+)=\infty$. Then, for all $t>0$, $r\in\NN$,
%
%
%e2.19 #&#
\begin{eqnarray}
\label{mainmod} && \bigl(\bigl\llvert\widetilde{\Delta X}^{(r)}_t
\bigr\rrvert,{}^{(r)}\widetilde\XX_t \bigr)
\nonumber\\[-8pt]\\[-8pt]\nonumber
&&\quad \stackrel{\mathrm{D}} {=} \Biggl(\overline{\Pi}^{\leftarrow}(v),
\XX^{\llvert \cdot\rrvert <\overline
{\Pi}^{\leftarrow}(v)}_t+ \sum_{i=1}^{Y^+(t\kappa^+(v))}
\delta_{(t\UUUU_i, \overline{\Pi
}^{\leftarrow}(v))} + \sum_{i=1}^{Y^-(t\kappa^-(v))}
\delta_{(t\UUUU_i', -\overline
{\Pi}^{\leftarrow}(v))} \Biggr)_{v=\SSSS_r/t},
\end{eqnarray}
where $\kappa^{\pm}(v)$ are the quantities in~\eqref{Ppar}.

\item[(ii)] Assume $\pibar^+(0+)=\infty$. Then for all $t>0$, $r\in\NN$,
\begin{eqnarray*}
%%\label{main+}
\bigl(\Delta X^{(r)}_t,{}^{(r)}
\XX_t^+ \bigr) \stackrel{\mathrm{D}} {=} \Biggl(\pibar^{+,\leftarrow}(v),
\XX^{+\cdot<\pibar
^{+,\leftarrow}(v)}_t+ \sum_{i=1}^{Y(t\kappa(v))}
\delta_{(t\UUUU
_i, \pibar^{+,\leftarrow}(v))} \Biggr)_{v=\SSSS_r/t},
\end{eqnarray*}
where $\kappa(v)=\pibar^+(\pibar^{+,\leftarrow}(v)-)-v$.
\end{longlist}
\end{theorem}

\begin{pf} %{Proof of Theorem~\ref{orderMod}}
Let $t>0$,
$r\in\NN$, and introduce a point measure $\widetilde\VV_t'$ as follows:
\[
\widetilde\VV_t':=\sum_{i\ge1}
\delta_{(t\UUUU_{i+r},\UUUU
'_{i+r},
(\SSSS_{i+r} - \SSSS_r)/t)}.
\]
Then $\widetilde\VV_t'$ is independent of $\mathfrak{V}:=\SSSS_r/t$
with $\widetilde\VV
_t'\stackrel{\mathrm{D}}{=}\VV'_t$. Observe that
%
%
%e2.20 #&#
\begin{eqnarray}\label{trunceqshift}
&&E\exp\biggl\{-\lambda\mathfrak{V}-\int f \,\rmd\bigl\{\delta
_{(0,0,\mathfrak{V}
)}\star\widetilde\VV'_t\bigr\} \biggr\}
\nonumber
\\
&&\quad =E\exp\biggl\{-\lambda\mathfrak{V}-\int_0^t
\int_0^1\int_{\mathfrak{V}
}^\infty
\bigl(1-e^{-f(s,u',v)} \bigr) \,\rmd s\,\rmd u' \,\rmd v \biggr\}
\\
&&\quad =E\exp\biggl\{-\lambda\mathfrak{V}-\int f\,\rmd\widetilde\VV
'_t{}^{\cdot\ge
\mathfrak{V}} \biggr\},\nonumber
\end{eqnarray}
for all nonnegative Borelian $f$ and $\lambda\ge0$. Here $\widetilde
\VV
'{}_t^{\cdot\ge v}(\cdot):=\widetilde\VV'{}_t([0,t]\times
(0,1)\times
[v,\infty)\cap\cdot)$.

Assume $\pibar(0+)=\infty$. Combining~\eqref{chpair1} and~\eqref
{trunceqshift} yields
%
%
%e2.21 #&#
\begin{eqnarray}
\label{triplelaw} \bigl(\bigl\llvert\widetilde{\Delta X}^{(r)}_t
\bigr\rrvert,{}^{(r)}\widetilde\XX_t \bigr) &\stackrel{
\mathrm{D}} {=}& \bigl(\overline{\Pi}^{\leftarrow}(\mathfrak{V}), \bigl
\{\delta
_{(0,0,\mathfrak{V})}\star\widetilde\VV_t'\bigr\}
^{m\circ(I,I,\overline{\Pi}^{\leftarrow})} \bigr)
\nonumber\\[-8pt]\\[-8pt]\nonumber
&\stackrel{\mathrm{D}} {=}& \bigl(\overline{\Pi}^{\leftarrow
}(\mathfrak{V}),
\bigl\{\widetilde\VV_t'{}^{\cdot\ge\mathfrak{V}}\bigr\}
^{m\circ(I,I,\overline{\Pi}^{\leftarrow})} \bigr).
\end{eqnarray}
Next, set $\YY_t:=\{\widetilde\VV_t'{}^{\cdot\ge\mathfrak{V}}\}
^{m\circ
(I,I,\overline{\Pi}^{\leftarrow})}$, and let
$\YY_t^{\llvert \cdot\rrvert <y}$, $\YY_t^{\cdot\ge y}$, $\YY_t^{\cdot
\le y}$
be the point processes obtained from $\YY_t$ by removal of points
not lying in the regions $[0,t]\times(-y,y)_*$, $[0,t]\times[y,\infty
)$, $[0,t]\times(-\infty,y]$, respectively.

Let $\lambda\ge0$ and $f_{-1},f_0,f_1:[0,t]\times\myRR_*\to
[0,\infty]$ be Borel functions.
Define $\Phi:[0,t]\times\myRR_*\times(0,\infty)\to[0,\infty]$ by setting
\[
\Phi[s,x,y]:= f_0(s,x) \bfeins_{(0,y)}\bigl(\llvert x
\rrvert\bigr)+f_1(s,x) \bfeins_{[y,\infty)}(x)+f_{-1}(s,x)
\bfeins_{(-\infty,-y]}(x).
\]
Observe that
%
%
%e2.22 #&#
\begin{eqnarray}\label{decompsumPhi}
&& \int_0^t\int_{\myRR_*}
\bigl(f_0\,\rmd\YY_t^{\llvert \cdot\rrvert <\overline{\Pi
}^{\leftarrow}
(\mathfrak{V})}+ f_1\,\rmd
\YY_t^{\cdot\ge\overline{\Pi}^{\leftarrow}(\mathfrak{V})}+ f_{-1}\,\rmd\YY
_t^{\cdot\le-\overline{\Pi}^{\leftarrow}(\mathfrak
{V})}
\bigr)
\nonumber\\[-8pt]\\[-8pt]\nonumber
&&\quad  =\int_0^t\int_{\myRR_*}
\Phi\bigl[\cdot,\cdot,\overline{\Pi}^{\leftarrow}(\mathfrak{V})\bigr]
\,\rmd
\YY_t
\end{eqnarray}
%
%\begin{eqnarray}\nonumber
%\lefteqn{\hspace*{-4cm}\int_0^t \int_{\myRR_*} \left(f_0(s,x)
%\YY_t^{\left\vert\cdot\right\vert<\pibarinv(\VVVV)}(\rmd s,\rmd x)+
%f_1(s,x) \YY_t^{\cdot\ge\pibarinv(\VVVV)}(\rmd s,\rmd x)+
%f_{-1}(s,x) \YY_t^{\cdot\le-\pibarinv(\VVVV)}(\rmd s,\rmd x)
%\right)}\\
%&=&\int_0^t\int_{\myRR_*}\Phi[s,x,\pibarinv(\VVVV)] \YY_t(\rmd s,\rmd
%x)\label{decompsumPhi}
%\end{eqnarray}
and
%
%
%e2.23 #&#
\begin{eqnarray}
\label{laplace1}
\hspace*{-10pt}&&E\exp\biggl\{-\lambda\mathfrak{V}- \int_0^t
\int_{\myRR_*}\Phi\bigl[s,x,\overline{\Pi}^{\leftarrow
}(\mathfrak{V})\bigr] \YY_t(\rmd s,\rmd x)\biggr\}\nonumber
\\[-2pt]
&&\quad =E\exp\biggl\{-\lambda\mathfrak{V}- \int_0^t
\int_0^1\int_{\mathfrak{V}}^\infty
\Phi\bigl[m\bigl(s,u',\overline{\Pi}^{\leftarrow} (v)\bigr),
\overline{\Pi}^{\leftarrow}(\mathfrak{V})\bigr] \widetilde\VV
'_t\bigl(\rmd s,\rmd u',\rmd v\bigr)
\biggr\}
\\[-2pt]
&&\quad =E\exp\biggl\{-\lambda\mathfrak{V} -\int_0^t
\int_0^1\int_{\mathfrak{V}}^\infty
\bigl(1-\exp\bigl\{ -\Phi\bigl[m\bigl(s,u',\overline{
\Pi}^{\leftarrow}(v)\bigr),\overline{\Pi}^{\leftarrow
}(\mathfrak{V})\bigr]
\bigr\} \bigr)\,\rmd s\,\rmd u'\,\rmd v\biggr\}.
\nonumber
\end{eqnarray}
As $\{v>0:\overline{\Pi}^{\leftarrow}(v)<\overline{\Pi
}^{\leftarrow}(\mathfrak{V})\}\subseteq(\mathfrak{V},\infty)$ and
$\overline{\Pi}^{\leftarrow}(v)= \overline{\Pi}^{\leftarrow
}(\mathfrak{V})$ for $v\in[\mathfrak{V},\pibar(\overline{\Pi
}^{\leftarrow}
(\mathfrak{V})-)]$,
the last integral in the exponent equals
\begin{eqnarray}
\label{integrfnu} &&\int_0^t\int
_0^1\int_{0}^\infty
\bfeins_{(0,\overline{\Pi}^{\leftarrow}(\mathfrak{V}))}\bigl(\overline
{\Pi}^{\leftarrow}(v)\bigr)
\bigl(1-e^{-f_0(m(s,u',\overline{\Pi}^{\leftarrow}(v)))} \bigr) \,\rmd
s\,\rmd u'\,\rmd v
\nonumber\\[-8pt]\\[-8pt]\nonumber
&&\quad {}+\kappa^+(\mathfrak{V}) \int_0^t
\bigl(1-e^{-f_1(s,\overline{\Pi}^{\leftarrow}(\mathfrak
{V}))} \bigr) \,\rmd s +\kappa^-(\mathfrak{V}) \int
_0^t \bigl(1-e^{-f_{-1}(s,-\overline{\Pi}^{\leftarrow
}(\mathfrak{V}))} \bigr) \,\rmd s,
\nonumber
\end{eqnarray}
with $\kappa^\pm(v)$ as in \eqref{Ppar}.
It follows from \eqref{imageM2} and a change of variables that
\begin{eqnarray*}
&&\int_0^t\int_0^1
\int_{0}^\infty\bfeins_{(0,\overline{\Pi}^{\leftarrow}(\mathfrak
{V}))}\bigl(
\overline{\Pi}^{\leftarrow}(v)\bigr) \bigl(1-e^{-f_0(m(s,u',\overline
{\Pi}^{\leftarrow}(v)))} \bigr) \,\rmd s\,\rmd
u'\,\rmd v
\nonumber
\\[-2pt]
&&\quad= \int_0^t\int_{(-\overline{\Pi}^{\leftarrow}(\mathfrak
{V}),\overline{\Pi}^{\leftarrow}(\mathfrak{V}))}
\bigl(1-e^{-f_0(s,x)} \bigr) \,\rmd s\Pi(\rmd x).
\end{eqnarray*}
We get from \eqref{decompsumPhi}, \eqref{laplace1} and \eqref{integrfnu}
%
%
%e2.24 #&#
\begin{eqnarray}
\label{quadruplelaw}
&& E\exp\biggl\{-\lambda\mathfrak{V}-\int
_0^t\int_{\myRR
_*}
\bigl(f_0\,\rmd\YY_t^{\llvert \cdot\rrvert <\overline{\Pi}^{\leftarrow
}(\mathfrak{V})}+ f_1\,\rmd
\YY_t^{\cdot\ge\overline{\Pi}^{\leftarrow}(\mathfrak{V})}+ f_{-1}\,\rmd\YY
_t^{\cdot\le-\overline{\Pi}^{\leftarrow}(\mathfrak
{V})}
\bigr) \biggr\}\nonumber
\\[-2pt]
&&\quad=E\exp\Biggl\{-\lambda\mathfrak{V}-\XX'_t{}^{\llvert \cdot
\rrvert <\overline
{\Pi}^{\leftarrow}(\mathfrak{V}
)}(f_0)
\\[-2pt]
&&\qquad{}-\sum_{i=1}^{Y^+(t\kappa^+(\mathfrak{V}))}f_1\bigl(t
\UUUU_i, \overline{\Pi}^{\leftarrow} (\mathfrak{V})\bigr) -\sum
_{i=1}^{Y^-(t\kappa^-(\mathfrak{V}))}f_{-1}\bigl(t
\UUUU_i', -\overline{\Pi}^{\leftarrow} (\mathfrak{V})\bigr) \Biggr\},
\nonumber
\end{eqnarray}
completing the proof of the following identity in law:
\begin{eqnarray*}
&& \bigl(\mathfrak{V}, \YY_t^{\llvert \cdot\rrvert <\overline{\Pi
}^{\leftarrow
}(\mathfrak{V})}, \YY_t^{\cdot\ge\overline{\Pi}^{\leftarrow}(\mathfrak{V})},
\YY_t^{\cdot\le-\overline{\Pi}^{\leftarrow}
(\mathfrak{V})} \bigr)
\\[-2pt]
&&\quad \stackrel{\mathrm{D}} {=} \Biggl(\mathfrak{V}, \XX'_t{}^{\llvert
\cdot
\rrvert <\overline{\Pi}^{\leftarrow}(\mathfrak{V})}, \sum
_{i=1}^{Y^+(t\kappa^+(\mathfrak{V}))}\delta_{(t\UUUU_i,
\overline{\Pi}^{\leftarrow}
(\mathfrak{V}))},
\sum_{i=1}^{Y^-(t\kappa^-(\mathfrak{V}))}\delta_{(t\UUUU_i',
-\overline{\Pi}^{\leftarrow}
(\mathfrak{V}))}
\Biggr),
\nonumber
\end{eqnarray*}
where
$\mathfrak{V},\XX'_t,Y^+,Y^-,(\UUUU_i),(\UUUU'_i)$ are independent with
$\XX_t'\stackrel{\mathrm{D}}{=}\XX_t$.
The proof of part~(i) is completed by combining \eqref{triplelaw} and
\eqref{quadruplelaw}. The proof of part~(ii) is similar.
\end{pf}

%%%%==============subSECTION------------------------------------------------

%============================proof of representaiton
%thm==============================
%-----------------------------------------------------------------
%s2.4 #&#
\subsection{Representations for the $r$-trimmed L\'evy
processes}\label{sub24}

By the L\'evy--It\^o decomposition (Sato \cite{sato1}, Theorem 19.2, page~120),
we can decompose a real-valued L\'evy process $X_t$, defined on the
probability space $(\Omega, \FFF, P)$, as
%
%
%e2.25 #&#
\begin{equation}
\label{decomp3} X_t=\gamma t+\sigma Z_t+X_t^{(J)},\qquad
t\ge0,
\end{equation}
where $\gamma\in\myRR$, $\sigma\ge0$, $(Z_t)_{t\ge0}$ is a
standard Brownian motion, and
$(X_t^{(J)})_{t\ge0}$, the jump process of $X$, is independent of
$(Z_t)_{t\ge0}$.
It satisfies, locally uniform in $t\ge0$,
%
%e2.26 #&#
\begin{equation}
\label{j1a} X_t^{(J)} =\mbox{a.s. }\lim
_{\varepsilon\downarrow0} \biggl(\sum_{0<s\le t}\Delta
X_s \bfeins_{\{\llvert \Delta X_s \rrvert >\veps\}} -t\int_{\veps
<\llvert x\rrvert \le1} x\Pi(\mathrm{d}x) \biggr).
\end{equation}
Now we can complete the proof of Theorem \ref{randrep}.

\begin{pf*}{Proof of Theorem \ref{randrep}}
We will prove part~(i), the identity for the $r$-fold modulus trimmed
L\'evy process. Trimming of positive jumps\vspace*{1pt} as in part~(ii) follows similarly.
Let $t > 0$, $r \in\NN$ be fixed. By \eqref{decomp3} and the
definition of $^{(r)}\widetilde X_t$, the $r$-fold modulus trimmed L\'evy
process is
\[
%%\label{ktrimlevy1}
^{(r)}\widetilde X_{t} = \gamma t+\sigma
Z_{t} + X_{t}^{(J)} - \sum
_{i=1}^r \widetilde{\Delta X}_t^{(i)},\qquad
t>0. %%\bfeins_{[\wt T^{(i)}_t,t]}(t), t>0.
\]
Note that the jump process of $^{(r)}\widetilde X_t$ and its quadratic
variation are obtained by applying the summing functional to the
$r$-fold modulus trimmed point process $^{(r)} \widetilde\XX$ and to the
squared jumps of $^{(r)}\widetilde\XX$. Using \eqref{j1a}, we can write
%% uniformly in $0\le u\le1$
%
%e2.27 #&#
\begin{eqnarray}
\label{ktrimlevy2}
&& X_{t}^{(J)} - \sum
_{i=1}^r \widetilde{\Delta X}_{t}^{(i)}
%%%%
%\bfeins_{[\wt
%T^{(i)}_t,t]}(t)
= \mbox{a.s. } \lim_{\veps\downarrow0}
\biggl(\int_{[0,t] \times\{
\llvert x\rrvert >
\veps\}} x {}^{(r)}\widetilde\XX(\rmd s,
\rmd x) -t\int_{\veps
<\llvert x\rrvert \le1} x \Pi(\rmd x) \biggr).
\end{eqnarray}
The corresponding $r$-trimmed quadratic variation is simply
\[
%%\label{trimQua}
^{(r)}\widetilde V_t = \int_{[0,t]\times\myRR_*}
x^2 {}^{(r)}\widetilde\XX(\rmd s, \rmd x).
\]
Recall from Lemma~\ref{chain} and Theorem \ref{orderMod} that the
distribution of $^{(r)}\widetilde\XX_t$ %, and $ \left\vert\wt{\Delta
%X}_t^{(r)}\left\vert$
can be decomposed as the superposition of three independent point
measures, as in \eqref{mainmod}. Splitting the integral in \eqref
{ktrimlevy2} into these components gives
%%, uniformly in $0\le u\le1$,
%
\begin{eqnarray*}
%%\label{ktrimLevy2a}
&&\mbox{a.s. } \lim_{\veps\downarrow0} \biggl(\int
_{[0,t] \times\{
\llvert x\rrvert >
\veps\}} x {}^{(r)}\widetilde\XX(\rmd s, \rmd x) - t
\int_{\veps
<\llvert x\rrvert \le1} x \Pi(\rmd x) \biggr)
\nonumber
\\
&&\quad \stackrel{\mathrm{D}} {=} \mbox{a.s. } \lim_{\veps\downarrow0} \biggl(
\int_{[0,t] \times\{\llvert x\rrvert >
\veps\}} x \XX^{\llvert \cdot\rrvert <\overline{\Pi}^{\leftarrow
}(\SSSS
_r/t)}(\rmd s, \rmd x)- t\int
_{\veps<\llvert x\rrvert \le1} x \Pi(\rmd x) \biggr)
\nonumber
\\
&&\qquad{}+ \overline{\Pi}^{\leftarrow}(\SSSS_r/t) \bigl( Y^+\bigl(t
\kappa^+(\SSSS_r/t)\bigr)-Y^-\bigl(t\kappa^-(\SSSS_r/t)
\bigr) \bigr).
\end{eqnarray*}
A similar expression holds for $^{(r)} \widetilde V_t$.
%%Letting $u=1$ gives
Thus, we conclude
\begin{eqnarray*}
%%\label{ktrimlevy4}
\bigl({}^{(r)}\widetilde X_t, {}^{(r)}
\widetilde V_t, \bigl\llvert\widetilde{\Delta X}_t^{(r)}
\bigr\rrvert\bigr) &\stackrel{\mathrm{D}} {=}& \bigl\{ \widetilde{X}_t^{v}
+ \overline{\Pi}^{\leftarrow}(v) \bigl(Y^+_{t\kappa
^+(v)} -
Y^-_{t\kappa
^-(v)}\bigr),
\\
&& \widetilde{V}_t^{v} + \overline{
\Pi}^{\leftarrow}(v)^2 \bigl(Y^+_{t\kappa^+(v)} +
Y^-_{t\kappa
^-(v)}\bigr), \overline{\Pi}^{\leftarrow}(v) \bigr
\}_{v=\SSSS_r/t}.
\end{eqnarray*}
This is \eqref{randdis1} and completes the proof of part~(i).
\end{pf*}

This completes our derivation of the trimming identities. In the next
sections, we turn to applications of them.

%%%%%%SECTION 3%%%%%%%%%%%%%%%%%%%%%%%%%%%%%%%%%%%%

%s3 #&#
\section{$X$ comparable with its large jump processes}\label{s3}
In this section, we apply Theorem \ref{randrep}
to complete a result of Maller and Mason \cite{MM10} concerning the
ratio of the process to its jump of largest magnitude. Note that when
$\pibar(0+)=\infty$, we have
$\llvert \widetilde{\Delta X}_t^{(1)}\rrvert =\sup_{0<s\le t}\llvert
\Delta X_s\rrvert >0$ a.s.
for all
$t>0$; similarly,
when $\pibar^+(0+)=\infty$, ${\Delta X}_t^{(1)}=\sup_{0<s\le
t}\Delta X_s>0$ a.s. for all $t>0$. Recall that $\pibar(x)$ is said to
be slowly varying (SV) as $x\downarrow0$ if
$\lim_{x\downarrow0}\pibar(ux)/\pibar(x)=1$ for all $u>0$ (e.g., Bingham,
Goldie and Teugels \cite{BGT87}).

%
%th2 #&#
\begin{theorem}\label{tX}
Suppose $\sigma^2=0$ and $\pibar(0+)=\infty$. Then
%
%
%e3.1 #&#
\begin{equation}
\label{rlim} \frac{X_t}{\widetilde{\Delta X}_t^{(1)}} \buildrel P \over
\to 1,\qquad\mbox{as }t\downarrow0,
\end{equation}
iff $\pibar(x)\in SV$ at 0
(so that $X$ is of bounded variation) and $X$ has drift 0. These imply
%
%
%e3.2 #&#
\begin{equation}
\label{slim2} \frac{\llvert \widetilde{\Delta X}_t^{(2)}\rrvert
}{\llvert \widetilde{\Delta X}_t^{(1)}\rrvert } \buildrel P \over\to
0,\qquad\mbox{as }t\downarrow0;
\end{equation}
and conversely \eqref{slim2} implies $\pibar(x)\in SV$ at 0.
\end{theorem}

For the proof, we need two preliminary lemmas.
The first calculates a distribution related to the large jumps, and the
second applies Theorem \ref{randrep} to derive a useful inequality.
%-------
%
%le2 #&#
\begin{lemma}\label{lemtsum2}
Assume $\pibar(0+)=\infty$.
Then for $t>0$, $0<u<1$,
%
%
%e3.3 #&#
\begin{equation}
\label{t5bb} P \bigl(\bigl\llvert\widetilde{\Delta X}_t^{(2)}
\bigr\rrvert\le u\bigl\llvert\widetilde{\Delta X}_t^{(1)}
\bigr\rrvert\bigr) =t\int_{(0,\infty)} e^{-t\pibar(u\overline{\Pi
}^{\leftarrow
}(v))}\,\rmd v.
%%%=t\int_{(0,\infty)} e^{-t\pibar(uy)}\left\vert\rmd\pibar(y)\right
%\vert,
\end{equation}
%
%and (ii) for $t>0$, $x\ge y>0$,
%\begin{equation}\label{eqid0}
% P\left(\left\vert\wt{\Delta X}_t^{(2)}\right\vert\le y, \left\vert
%\wt{\Delta X}_t^{(1)}\right\vert\le x
%\right)
%=\left(t(\pibar(y)-\pibar(x))+1\right)e^{-t\pibar(y)}.
%\end{equation}
A similar expression to \eqref{t5bb} %%and \eqref{eqid0} are
is true when $\pibar^+(0+)=\infty$,
with $\llvert \widetilde{\Delta X}_t^{(1)}\rrvert $ and $\llvert
\widetilde{\Delta X}_t^{(2)}\rrvert $
replaced by ${\Delta X}^{(1)}_t$ and ${\Delta X}^{(2)}_t$, and $\pibar
$ and $\overline{\Pi}^{\leftarrow}$ replaced by $\pibar^+$ and
$\overline{\Pi}^{+,\leftarrow}$.
\end{lemma}

\begin{pf} %{Proof of Lemma \ref{lemtsum2}}
Assume $\pibar(0+)=\infty$ and
take $t>0$. We get from~\eqref{chpair01} that
%
%
%e3.4 #&#
\begin{equation}
\label{eqdinproofofLemlemtsum2} \bigl(\bigl\llvert\widetilde{\Delta X}^{(1)}_t
\bigr\rrvert,\bigl\llvert\widetilde{\Delta X}^{(2)}_t\bigr
\rrvert\bigr)\stackrel{\mathrm{D}} {=}\bigl(\overline{
\Pi}^{\leftarrow
}(\EEEE_1/t), \overline{\Pi}^{\leftarrow}\bigl((
\EEEE_1 + \EEEE_2)/t\bigr) \bigr),
\end{equation}
where $\EEEE_1$ and $\EEEE_2$ are independent unit exponential random
variables. %%with common parameter $E[\EEEE_1]=E[\EEEE_2]=1$.
%%For Part (i),
Take $0<u<1$ and $v>0$ and let $y_{t,u}(v):=t\pibar(u\overline{\Pi
}^{\leftarrow}(v/t))$.
%%Observe that $y_{t,u}(v)>v$ for all $v>0$.
Then, in view of~\eqref{eqdinproofofLemlemtsum2},
\begin{eqnarray*}
P \bigl(\bigl\llvert\widetilde{\Delta X}^{(2)}_t\bigr
\rrvert\le u \bigl\llvert\widetilde{\Delta X}^{(1)}_t\bigr
\rrvert\bigr) &=& P \bigl(\overline{\Pi}^{\leftarrow}\bigl((\EEEE_1
+ \EEEE_2)/t\bigr)\le u \overline{\Pi}^{\leftarrow}(\EEEE
_1/t) \bigr)
\\
&=&P \bigl(\EEEE_1 + \EEEE_2\ge y_{t,u}(
\EEEE_1) \bigr)
\\
&=&\int_{(0,\infty)}e^{-(y_{t,u}(v)-v)}e^{-v} \,\rmd v
\\
&=&\int_{(0,\infty)}\exp\bigl\{-t \pibar\bigl(u\overline{
\Pi}^{\leftarrow
}(v/t)\bigr)\bigr\} \,\rmd v.
\end{eqnarray*}
Changing the variable from $v/t$ to $v$ gives \eqref{t5bb}.
%%Changing variable from $v$ to
%$\pibar(y)$ then gives the second, by recalling the identities $(\rmd
%v)^{\pibarinv}=\rmd\Pi^{\left\vert\cdot\right\vert}=\left\vert\rmd
%\pibar\right\vert$.
%For Part (ii), take $x\ge y>0$. Then \eqref{eqid0} follows from
%\eqref{eqdinproofofLemlemtsum2} and
%\begin{eqnarray*}
% P\left(\left\vert\widetilde{\Delta X}^{(2)}_t\left\vert\le y, \left
%\vert\widetilde{\Delta
%X}^{(1)}_t\left\vert\le x\right)
%&=& P\left(\pibarinv((\EEEE_1 + \EEEE_2)/t)\le y, \pibarinv(
%\EEEE_1/t)\le x\right)
%\\
%&=& P\left(\EEEE_1 + \EEEE_2\ge t\pibar(y), \EEEE_1\ge t\pibar(x)
%\right)\\
%&=& P\left(\EEEE_2\ge t\pibar(y) - \EEEE_1, \EEEE_1\ge t\pibar(y)
%\right)\\
%{}&&{}+P\left(\EEEE_2\ge t\pibar(y) - \EEEE_1, t\pibar(y)>\EEEE_1
%\ge t\pibar(x)\right)\\
%&=& (1+t(\pibar(y)-\pibar(x))) e^{-t \pibar(y)}.
%\end{eqnarray*}
The version for large jumps, rather than jumps large in modulus, is
proved similarly.
\end{pf}

%%---symmetrisation lemma----
%%---proposition on lower bound of the trimmed inequality with
%centering.
%
%le3 #&#
\begin{lemma}\label{lemtsum4}
Assume $\pibar(0+)=\infty$, and let $a_t$ be any nonstochastic
function in $\myRR$. Then for $t>0$ and $0<u<1/4$,
%
%
%e3.5 #&#
\begin{equation}
\label{L1} 4P \bigl(\bigl\llvert{}^{(1)}\widetilde
X_t-a_t\bigr\rrvert> u\bigl\llvert\widetilde{\Delta
X}^{(1)}_t\bigr\rrvert\bigr) \ge P \bigl(\bigl\llvert
\widetilde{\Delta X}^{(2)}_t\bigr\rrvert>4u\bigl\llvert
\widetilde{\Delta X}^{(1)}_t\bigr\rrvert\bigr).
\end{equation}
%
%%at points $u$ of continuity of the righthand side.
Assuming $\pibar^+(0+)=\infty$, the same inequality \eqref{L1} holds with
${}^{(1)}X_t$, ${\Delta X}_t^{(1)}$ and ${\Delta X}^{(2)}_t$
in place of ${}^{(1)}\widetilde X_t$, $\llvert \widetilde{\Delta
X}_t^{(1)}\rrvert $
and $\llvert \widetilde
{\Delta X}^{(2)}_t\rrvert $.
\end{lemma}
%
%\begin{remark}\label{reSum4}
%Although \eqref{L1} is proved only for continuity points of the right
%hand side, it holds for all $u>0$. This fact is not needed for our
%present purposes so we omit the details.
%\end{remark}

\begin{pf} %{Proof of Lemma \ref{lemtsum4}}
Let $\EEEE$ be an exponential random variable with $E\EEEE=1$, thus,
$\EEEE\stackrel{\mathrm{D}}{=}\SSSS_1$.
Using the identity in \eqref{cor1a} with $r=1$, the left-hand side of
\eqref{L1} is, for $u>0$,
%
%e3.6 #&#
\begin{equation}
\label{L2} \int_0^\infty4P \bigl(\bigl\llvert
\widetilde X_t^v+\widetilde G_t^v-a_t
\bigr\rrvert> uy_v \bigr) P(\EEEE\in t\,\rmd v),
\end{equation}
where
% term from \eqref{cor1a} is absent since $\Pi$ is continuous, and
we abbreviate $y_v:=\overline{\Pi}^{\leftarrow}(v)$, $v>0$. For each
$v>0$, let
$(\overline X_t^v)_{t\ge0}$ and $(\overline G_t^v)_{t\ge0}$
be independent\vspace*{1pt} copies of $(\widetilde X_t^v)_{t\ge0}$ and $(\widetilde
G_t^v)_{t\ge0}$,
with $(\overline G_t^v)_{t\ge0}$ also independent of $(\overline
X_t^v)_{t\ge0}$. Define the symmetrised process $(\widehat
Y_t^v)_{t\ge
0}$ by
\[
\widehat Y_t^v= \bigl(\widetilde X_t^v+
\widetilde G_t^v \bigr)- \bigl(\overline
X_t^v+\overline G_t^v \bigr),\qquad
t>0,
\]
with jump process $\Delta\widehat Y_t^v=\widehat Y_t^v- \widehat Y_{t-}^v$,
$t>0$.
Then the integrand in \eqref{L2} satisfies
%
%e3.7 #&#
\begin{eqnarray}
\label{L3} 4P \bigl(\bigl\llvert\widetilde X_t^v+
\widetilde G_t^v-a_t\bigr\rrvert>
uy_v \bigr) &=& 2P \bigl(\bigl\llvert\widetilde X_t^v+
\widetilde G_t^v-a_t\bigr\rrvert>
uy_v \bigr)+ 2P \bigl(\bigl\llvert\overline X_t^v
+\overline G_t^v-a_t\bigr\rrvert>
uy_v \bigr)\nonumber
\\
&\ge& 2 P \bigl(\bigl\llvert\bigl(\widetilde
X_t^v+\widetilde G_t^v-a_t
\bigr)-\bigl(\overline X_t^v+\overline
G_t^v-a_t\bigr)\bigr\rrvert>
2uy_v \bigr)
\\
&=& 2P\bigl(\bigl\llvert\widehat
Y_t^v\bigr\rrvert>2uy_v\bigr).\nonumber
\end{eqnarray}
Substitute the inequality \eqref{L3} in \eqref{L2} and equate to the
left-hand side of \eqref{L1} to get
%
%e3.8 #&#
\begin{equation}
\label{L4} 4P \bigl(\bigl\llvert{}^{(1)}\widetilde
X_t-a_t\bigr\rrvert> u\bigl\llvert\widetilde{\Delta
X}^{(1)}_t\bigr\rrvert\bigr)\ge2\int_0^\infty
P \bigl(\bigl\llvert\widehat Y_t^v\bigr\rrvert>
2uy_v \bigr) P(\EEEE\in t\,\rmd v),\qquad u>0.
\end{equation}

Take $u \in(0, 1/4)$. Applying L\'evy's maximal inequality for random
walks (Feller (\cite{feller1971}, Lemma 2, page~147), we have
%
%e3.9 #&#
\begin{eqnarray}
\label{L5} 2 P \bigl(\bigl\llvert\widehat Y_t^v\bigr
\rrvert>2uy_v \bigr)&=& 2\lim_{n\to\infty} P \Biggl(
\Biggl\llvert\sum_{i=1}^{2^n}\bigl(\widehat
Y_{it/2^n}^v-\widehat Y_{(i-1)t/2^n}^v\bigr)
\Biggr\rrvert>2uy_v \Biggr)\nonumber
\\
&\ge& \lim_{n\to\infty}P
\Bigl(\max_{1\le j\le2^n}\bigl\llvert\bigl(\widehat
Y_{jt/2^n}^v-\widehat Y_{(j-1)t/2^n}^v\bigr)
\bigr\rrvert>2uy_v \Bigr)
\\
&\ge& %% P\left(\sup_{0<s\le t} \left\vert\widehat X_s^v\right\vert
%>2uy_v\right) \cr&\ge&
P \Bigl(\sup
_{0<s\le t} \bigl\llvert\Delta\widehat Y_s^v
\bigr\rrvert>4uy_v \Bigr).\nonumber
\end{eqnarray}
The L\'evy measure of $\widetilde X_t^v$ is
$\Pi(\rmd x)\bfeins_{\{\llvert x\rrvert <y_v\}}$, $x\in\myRR$, having
tail function
\[
\bigl(\pibar(x)-\pibar(y_v-)\bigr)\bfeins_{\{x<y_v\}},\qquad x>0.
\]
Suppose at first that $\Delta\pibar(y_v)>0$. %%$\pibar(y_v-)>v$.
Then $\widetilde G_t^v$ is nonzero. Its L\'evy measure consists of point
masses at $\pm y_v$ with magnitudes $\kappa^\pm(v)$, given by \eqref{Ppar}.
Hence, it has tail
\[
\bigl(\pibar(y_v- )-v \bigr) \frac{ (\Delta\Pi(y_v)+\Delta\Pi(-y_v)
)}{\Delta
\pibar( y_v )}\bfeins_{\{x<y_v\}}=
\bigl(\pibar(y_v- )-v \bigr) \bfeins_{\{x<y_v\}},\qquad x>0.
\]
Adding the two tails gives the tail of $\widetilde X_t^v+\widetilde
G_t^v$ as
$(\pibar(x)-v)\bfeins_{\{x<y_v\}}$, $x>0$.
The symmetrisation $\widehat Y_t^v$ has L\'evy tail being twice the
magnitude of this.
This result remains true when $\Delta\pibar(y_v)=0$, as $\widetilde
G_t^v\equiv0$ and $\pibar(y_v-)=v$ then.

We can now calculate the right-hand side of \eqref{L5} and deduce from
it that
%
%e3.10 #&#
\begin{eqnarray}
\label{L6} 2 P \bigl(\bigl\llvert\widehat Y_t^v\bigr
\rrvert>2uy_v \bigr) &\ge& 1-e^{-2t(\pibar(4uy_v)-v)}
\nonumber\\[-8pt]\\[-8pt]\nonumber
& \ge&
1-e^{-t(\pibar(4uy_v)-v)}.
\end{eqnarray}
Finally, \eqref{L4}, \eqref{L6} and Lemma \ref{lemtsum2} give
\begin{eqnarray*}
%%\label{L8}
4P \bigl(\bigl\llvert{}^{(1)}\widetilde X_t-a_t
\bigr\rrvert> u\bigl\llvert\widetilde{\Delta X}^{(1)}_t
\bigr\rrvert\bigr) &\ge& t\int_0^\infty
\bigl(e^{-tv}- e^{-t\pibar(4uy_v)} \bigr)\,\rmd v
\\
&=& P \bigl(\bigl\llvert
\widetilde{\Delta X}_t^{(2)}\bigr\rrvert>4u\bigl\llvert
\widetilde{\Delta X}_t^{(1)}\bigr\rrvert\bigr).
\end{eqnarray*}

This proves \eqref{L1}. To derive the version for $^{(1)}X_t$, define
the one-sided L\'evy process $X_t^*$ having triplet
$(\gamma,\sigma^2, \Pi^*(\rmd x)= \Pi(\rmd x)\bfeins_{(x>0)}$),
and let $\widetilde{\Delta X}_t^{*,(r)}$ be the jump of $r$th largest
modulus up until time $t$ for $(X^*_t)_{t\ge0}$, $r\in\NN$. Then
$\widetilde{\Delta X}_t^{*,(r)}= \Delta X_t^{(r)}$ and
${}^{(1)}X_t= {}^{(1)}\widetilde X_t^*= X_t^*- \widetilde{\Delta X}_t^{*,(1)}$.
Assuming $\pibar^+(0+)=\infty$, inequality \eqref{L1} with
${}^{(1)}X_t$, ${\Delta X}_t^{(1)}$ and ${\Delta X}^{(2)}_t$
replacing ${}^{(1)}\widetilde X_t$, $\llvert \widetilde{\Delta
X}_t^{(1)}\rrvert $
and $\llvert \widetilde
{\Delta X}^{(2)}_t\rrvert $ then follows from \eqref{L1} itself,
applied to $X_t^*$.
\end{pf}

%
%le4 #&#
\begin{lemma}\label{lemtsum5}
Assume $\pibar(0+)=\infty$. Then
%
%
%e3.11 #&#
\begin{equation}
\label{B1} \frac{\llvert \widetilde{\Delta X}_t^{(2)}\rrvert
}{\llvert \widetilde{\Delta X}_t^{(1)}\rrvert } \buildrel P \over\to
0,\qquad\mbox{as }t\downarrow0,
\end{equation}
implies $\pibar(x)$ is SV at 0.
\end{lemma}

\begin{pf} %{Proof of Lemma \ref{lemtsum5}}
From~\eqref{t5bb}, for $0<u<1$,
with $y_v:=\overline{\Pi
}^{\leftarrow}(v)$,
%
%e3.12 #&#
\begin{eqnarray}
\label{B2} P \bigl(\bigl\llvert\widetilde{\Delta X}_t^{(2)}
\bigr\rrvert>u\bigl\llvert\widetilde{\Delta X}_t^{(1)}\bigr
\rrvert\bigr) &=& %1- t\int_{(0,\infty)} e^{-t\pibar(uy)}\left\vert
%\rmd\pibar(y)\right\vert\cr
%&=&
t\int_0^\infty
\bigl(e^{-tv}- e^{-t\pibar(uy_v)} \bigr)\,\rmd v
\nonumber\\[-8pt]\\[-8pt]\nonumber
&=& \int
_0^\infty\bigl(e^{-v}-
e^{-t\pibar(uy_{v/t})} \bigr)\,\rmd v.
\end{eqnarray}
Assume \eqref{B1}, so the integral on the right-hand side of \eqref
{B2} tends to 0 as $t\downarrow0$. %% Recall that $y_v=\pibarinv(v)$,
%$v>0$.
Take any sequence $t_k\downarrow0$ and by Helly's theorem select for each
$u>0$ a subsequence $t_{k'}=t_{k'}(u)\downarrow0$ such that
$t_{k'}\pibar
(uy_{v/t_{k'}})$ converges vaguely to $g_u(v)$, as $k'\to\infty$,
where $g_u(v)$ is a monotone function of~$v$.
Since $t\pibar(uy_{v/t})\ge t\pibar(y_{v/t}-)\ge v$, we have
$g_u(v)\ge v$.
Fatou's lemma applied to \eqref{B2} shows then that $g_u(v)=v$ for
$v>0$, thus
$t_{k'}\pibar(uy_{v/t_{k'}})\to v$, and since this is true for all
subsequences we deduce
\[
\lim_{t\downarrow0} t\pibar(uy_{v/t})=v,\qquad v>0, 0<u<1.
\]
Given $x>0$, $v>0$, let $t(x)=v/\pibar(x)$. Then
%% $\pibar(y_{v/t(x)})=\pibar(\pibarinv(v/t(x)))=x$, and
$ y_{v/t(x)}= \overline{\Pi}^{\leftarrow}(\pibar(x))\le x$, implying
$ \pibar(uy_{v/t(x)})\ge\pibar(ux)$. So we get, for
$0<u<1$,
\[
1\le\frac{\pibar(ux)}{\pibar(x)}\le\frac{t(x)\pibar(uy_{v/t(x)})}{v} \to
\frac{v}{v}= 1,\qquad\mbox{as }x
\downarrow0,
\]
and $\pibar\in SV$ at 0.
\end{pf}

\begin{pf*}{Proof of Theorem \ref{tX}}
Observe that \eqref
{rlim} is equivalent to
\[
%%\label{slim}
\frac{\llvert {}^{(1)}\widetilde X_t\rrvert }{\llvert \widetilde
{\Delta X}_t^{(1)}\rrvert } \buildrel P \over\to 0,\qquad\mbox{as
}t\downarrow0,
\]
and this implies \eqref{B1} by Lemma \ref{lemtsum4}. Thus,
by Lemma \ref{lemtsum5}, $\pibar\in SV$ at 0. Hence, $\int_0^1\pibar
(x)\,\rmd x<\infty$ and $X$ is
of bounded variation, with drift $\rmd_X$. By, for example, Bertoin
(\cite{b}, Proposition 11, page~167), $X_t/t \buildrel P \over\to
\rmd_X$ as $t\downarrow0$,
while, for any $\delta>0$,
\[
P\Bigl(\sup_{0<s\le t}\llvert\Delta X_s\rrvert>
\delta t\Bigr)= 1-e^{-t\pibar(\delta t)} \to0,
\]
thus $\widetilde{\Delta X}^{(1)}_t/t \buildrel P \over\to 0$ as
$t\downarrow0$.
But
\[
\frac{\llvert X_t\rrvert }{t} = \frac{\llvert X_t\rrvert }{\llvert
\widetilde{\Delta X}^{(1)}_t\rrvert }\cdot\frac{\llvert \widetilde
{\Delta X}^{(1)}
_t\rrvert }{t}
\buildrel P \over\to (1) (0)=0,
\]
showing that $\rmd_X=0$.

Conversely, \eqref{rlim} holds when $\pibar\in SV$ at 0 and $\rmd
_X=0$, as shown
in Lemma 5.1 of Maller and Mason \cite{MM10}.
\end{pf*}

The next result follows by applying Theorem \ref{tX} to the L\'evy process
$(\sum_{0<s\le t}\llvert \Delta X_s\rrvert )_{t>0}$, when $X_t$ is of
bounded variation.

%
%co2 #&#
\begin{corollary}\label{cortX}
Suppose $\sigma^2=0$ and $\pibar(0+)=\infty$. $X_t$ is of bounded
variation and
\[
%%\label{rlimcor}
\frac{\sum_{0<s\le t}\llvert \Delta X_s\rrvert }{\sup_{0<s\le
t}\llvert \Delta X_s\rrvert } \buildrel P \over\to 1,\qquad\mbox{as
}t\downarrow0,
\]
iff $\pibar(x)\in SV$ at 0
(so that $X$ is of bounded variation) and $X$ has drift 0.
\end{corollary}

%
%re1 #&#
\begin{remark}\label{rem4}
As another corollary of Theorem \ref{tX}, it is not hard to show that
$\pibar(x)\in SV$ at 0 implies $t\pibar(\llvert X_t\rrvert )\stackrel
{\mathrm
{D}}{\longrightarrow}\EEEE$ as $t\downarrow0$.
The variance gamma model, widely used in
financial modelling, has L\'evy measure whose tail is slowly varying at
0 (Madan and Seneta (\cite{ms}, page~519)).
Our results for such processes provide useful intuition and, more specifically,
may be of immediate use in applications, such as for estimation of $\Pi
$ or simulation, and so forth.
\end{remark}

The next theorem gives a one-sided version of Theorem \ref{tX}.
Condition \eqref{pmrlim} reflects a kind of dominance of the positive
part of $X$ over its negative part. We defer the proof of Theorem \ref
{pmtX} to the following section, where we study such dominance ideas in detail.

%
%th3 #&#
\begin{theorem}\label{pmtX}
Suppose %%$\sigma^2=0$ and
$\pibar^+(0+)=\infty$. Then
%
%
%e3.13 #&#
\begin{equation}
\label{pmrlim} \frac{X_t}{{\Delta X}_t^{(1)}} \buildrel P \over\to
1,\qquad\mbox{as }t\downarrow0
\end{equation}
iff $\pibar^+(x)\in SV$ at 0, $X$ is of bounded variation with drift
0, and $\lim_{x\downarrow0} \pibar^-(x)/\pibar^+(x)=0$.
\end{theorem}

%%====================================================
%% Positive and negative jumps
%%
%%============================================
%s4 #&#
\section{Comparing positive and negative jumps}\label{s+-}
In this section, we
%%only deal with an $X_t$ of bounded variation such that
%%$\pibar^+(0+)=\pibar^-(0+)=\infty$. We
are concerned with comparing magnitudes of positive and negative jumps
of $X$, in various ways.
Define $\Delta X_t^+=\max(\Delta X_t,0)$, $\Delta X_t^-=\max(-\Delta
X_t,0)$, and
\[
\bigl(\Delta X^+\bigr)_t^{(1)}=\sup_{0<s\le t}
\Delta X_s^+ \quad\mbox{and}\quad\bigl(\Delta X^-
\bigr)_t^{(1)}=\sup_{0<s\le t}\Delta
X_s^-,\qquad t>0.
\]

In the Poisson point process of jumps $(\Delta X_t)_{t>0}$, the numbers
of jumps and their magnitudes in disjoint regions are independent.
Thus, the positive and negative jump processes are independent. When
the integrals are finite, define
\[
A_\pm(x):= \int_0^x
\pibar^\pm(y)\,\rmd y= x\int_0^1
\pibar^\pm(xy)\,\rmd y.
\]
We obtain the following.

%, for $\veps>0$,
%\[
%N^+_t(\veps):= \#\{{\rm number of jumps with} \Delta X_t>\veps\},
%\]
% and the magnitudes of these jumps, are independent of $N^-_t(\veps)$,
%the
% number of jumps with $\Delta X_t<-\veps$, and the magnitudes of these
%jumps.
% With $N_t(\veps)$ and $J_i(\veps)$ as in \eqref{j4ab}, let $J_i^+(
%\veps)=\max(J_i(\veps),0)$
% and $J_i^-(\veps):= J_i^+(\veps)-J_i(\veps)$.
% Notice that the sum of the magnitudes of the positive jumps is
% \[
% \sum_{i=1}^{N^+_t(\veps)} J_i(\veps)\bfeins_{\{J_i(\veps)>0\}}
% = \sum_{i=1}^{N_t(\veps)} J_i^+(\veps),
% \]
% as long as $N^+_t(\veps)\ge1$, and similarly for the negative jumps;
% also,
%\[
% \max_{1\le i\le N^+_t(\veps)} J_i(\veps)\bfeins_{\{J_i(\veps)>0\}}
% = \max_{1\le i\le N_t(\veps)} J_i^+(\veps),
% \]
%if $N^+_t(\veps)\ge1$, and similarly for the negative jumps.
%Thus in particular the 2-vectors
%\begin{equation}\label{+1}
%\left( \sum_{i=1}^{N_t(\veps)} J_i^+(\veps), \max_{1\le i\le N_t(
%\veps)} J_i^+(\veps)\right) \quad{\rm and}\quad
%\left( \sum_{i=1}^{N_t(\veps)} J_i^-(\veps), \max_{1\le i\le N_t(
%\veps)} J_i^-(\veps)\right)
%\end{equation}
%are independent.
%

%
%th4 #&#
\begin{theorem}\label{t+-}
Suppose $\pibar^\pm(0+)=\infty$. %% and $X$ is of bounded variation.
%Then
For \eqref{pm} assume $\sum_{0<s\le t}\Delta X_s^-$ is finite a.s.,
and for \eqref{Rs1} assume $\sum_{0<s\le t}\Delta X_s^+$ is finite a.s.
For \eqref{Rrs3}, assume both are finite a.s.
Then
%
%
%e4.1 #&#
\begin{equation}
\label{pm} \frac{\sum_{0<s\le t}\Delta X_s^-}{\sup_{0<s\le t}\Delta
X_s^+} \buildrel P \over\to 0,\qquad\mbox{as }t\downarrow0
\quad\mbox{if and only if}\quad\lim_{x\downarrow0}\frac{\int_0^x\pibar
^-(y)\,\rmd y}{x\pibar
^+(x)}= 0;
\end{equation}
also
%
%
%e4.2 #&#
\begin{equation}
\label{Rs1} \frac{\sup_{0<s\le t}\Delta X_s^-}{\sum_{0<s\le t}\Delta
X_s^+} \buildrel P \over\to 0,\qquad\mbox{as }t
\downarrow0\quad\mbox{if and only if}\quad\lim_{x\downarrow0}\frac
{x\pibar^-(x)}{\int_0^x\pibar^+(y)\,\rmd
y}=
0;
\end{equation}
and
%
%
%e4.3 #&#
\begin{equation}
\label{Rrs3} \frac{\sum_{0<s\le t}\Delta X_s^-}{\sum_{0<s\le t}\Delta
X_s^+} \buildrel P \over\to 0,\qquad\mbox{as }t
\downarrow0\quad\mbox{if and only if}\quad\lim_{x\downarrow0}\frac{\int
_0^x\pibar^-(y)\,\rmd y}{\int_0^x\pibar
^+(y)\,\rmd y}=
0.
\end{equation}
Finally,
%
%
%e4.4 #&#
\begin{equation}
\label{spm} \frac{\sup_{0<s\le t}\Delta X_s^-}{\sup_{0<s\le t}\Delta
X_s^+} \buildrel P \over\to 0,\qquad\mbox{as }t
\downarrow0\quad\mbox{if and only if}\quad\lim_{x\downarrow0}\frac
{\pibar^-(\veps x)}{\pibar^+(x)}=
0\qquad\mbox{for all } \veps>0.
\end{equation}
\end{theorem}

\begin{pf}%{Proof of Theorem \ref{t+-}}
To prove the equivalence in \eqref{pm},
note that, for any $\lambda> 0$,
%
%e4.5 #&#
\begin{eqnarray}
\label{L0spl} E\exp\biggl( -\lambda\frac{\sum_{0<s\le t}\Delta
X_s^-}{\sup_{0<s\le t} \Delta X_s^+} \biggr) &=& E E \biggl[
\exp\biggl( -\frac{\lambda}{
(\Delta X^+)_t^{(1)}} \sum_{0<s\le t}\Delta
X_s^- \biggr) \Big| \bigl(\Delta X^+\bigr)_t^{(1)}
\biggr]
\nonumber
\\
&=& E \biggl[ \exp\biggl( -t\int_{(0,\infty)}
\bigl(1-e^{-\lambda x/(\Delta
X^+)_t^{(1)}}\bigr) \Pi^{(-)}(\rmd x) \biggr) \biggr]
\\
&=& \int_{(0,\infty)}\exp\biggl(-t\int_{(0,\infty)}
\bigl(1- e^{-\lambda x/y} \bigr)\Pi^{(-)}(\rmd x) \biggr)\lambda
^+_t(\rmd y),\nonumber
\end{eqnarray}
where
\[
\lambda_t^+(x):=P \Bigl(\sup_{0<s\le t}\Delta
X_s\le x \Bigr)=e^{-t\pibar^+(x)},\qquad x>0, t>0.
\]
%
%%%$\lambda^+_t(\rmd y)=\rmd(e^{-t\pibar^+(y)})$, $y>0$.
By \eqref{L0spl}, the left-hand relation in \eqref{pm} holds if and
only if, for all $\lambda>0$,
%
%
%e4.6 #&#
\begin{equation}
\label{L21} \lim_{t\downarrow0} \int_{(0,\infty)}
\bigl(1- e^{-t\int_{(0,\infty)} (1- e^{-\lambda x/y} )\Pi
^{(-)}(\rmd x)} \bigr) \lambda_t^+(\rmd y)=0.
\end{equation}
Use the lower bound in the inequalities
(cf. Bertoin \cite{b}, Proposition 1, page~74)
%
%e4.7 #&#
\begin{equation}
\label{uy-} (\lambda/3y)A_-(y/\lambda) \le\int_{(0,\infty)}
\bigl(1- e^{-\lambda x/y} \bigr)\Pi^{(-)}(\rmd x) \le(\lambda/y)A_-(y/
\lambda),\qquad y>0, \lambda>0,
\end{equation}
with $\lambda=1$ to get a lower bound for the integral in \eqref{L21} of
%
%
%e4.8 #&#
\begin{eqnarray}
\label{u++} \int_{(0,\infty)} \bigl(1-e^{-tA_-(y)/3y} \bigr)
\lambda^+_t(\rmd y) \ge\int_{(0,z]}
\bigl(1-e^{-tA_-(y)/3y} \bigr)\lambda^+_t(\rmd y)
\end{eqnarray}
for any $z>0$.
It is easily checked that $A_-(z)/z$ is nonincreasing for $z>0$, so the
last integral in \eqref{u++} is not smaller than
\[
\bigl(1-e^{-tA_-(z)/3z} \bigr)\lambda^+_t(z)= \bigl(1-e^{-tA_-(z)/3z}
\bigr) e^{-t\pibar^+(z)}.
\]
Now choose $t=1/\pibar^+(z)$ and let $t\downarrow0$ (so $z\downarrow
0$) to get
the righthand relation in \eqref{pm}.

Conversely, assume the right-hand relation in \eqref{pm}. Then the
upper bound in \eqref{uy-} shows that the integral in \eqref{L21} is
no larger than
\[
\int_{(0,\infty)} \bigl(1-e^{-\lambda tA_-(y/\lambda)/y} \bigr)
\lambda^+_t(\rmd y).
\]
This is a nondecreasing function of $\lambda$ so it suffices to show
that it tends to 0 as $t\downarrow0$ for $\lambda>1$. Then since
$A_-$ is
nondecreasing, for any $z>0$ the integral is bounded above by
\[
\int_{[z,\infty)} \bigl(1-e^{-\lambda tA_-(y)/y} \bigr)
\lambda^+_t(\rmd y)+\lambda_t^+(z-) \le1-e^{-\lambda
tA_-(z)/z}+e^{-t\pibar^+(z-)}.
\]
Take $t>0$ and $a>0$ and let $z=\pibar^{+,\leftarrow}(a/t)$. Then the
last expression is no larger than
\[
1-e^{-a\lambda A_-(z)/z\pibar^+(z)}+e^{-a}.
\]
Letting $t\downarrow0$, so $z\downarrow0$, then $a\to\infty$, this
tends to 0
by the right-hand relation in \eqref{pm}.

The equivalence in \eqref{Rs1} is proved similarly to that in \eqref
{pm}, by reversing the numerator and denominator and interchanging
$+/-$ and noting that the left-hand relation in \eqref{Rs1} holds if
and only if the Laplace transform of the ratio on the left of \eqref
{Rs1} tends to 1 as $t\downarrow0$.

The equivalence in \eqref{Rrs3} can be inferred from that in \eqref
{Rs1} with the following device. The left-hand relation in \eqref
{Rrs3} holds if and only if
%
%e4.9 #&#
\begin{equation}
\label{RS3} \lim_{t\downarrow0} E\exp\biggl(-\lambda\sum
_{0<s\le t}\Delta X_s^-\Big/\sum
_{0<s\le
t}\Delta X_s^+ \biggr) =1\qquad\mbox{for all }
\lambda>0.
\end{equation}
The Laplace transform on the left-hand side of \eqref{RS3} equals
%
%e4.10 #&#
\begin{equation}
\label{Lapx} \int_{(0,\infty)}\exp\biggl(-t\int
_{(0,\infty)} \bigl(1- e^{-x/y} \bigr)\Pi^{(-)}(
\rmd x) \biggr) P \biggl(\sum_{0<s\le t}\Delta
X_s^+\in\lambda\,\rmd y \biggr).
\end{equation}
Define a measure $\rho(\cdot)$ on $(0,\infty)$ in terms of its tail:
\[
\rhobar(y)= \int_{(0,\infty)} \bigl(1- e^{-x/y} \bigr)
\Pi^{(-)}(\rmd x),\qquad y>0.
\]
Then $\rhobar(y)$ is nonincreasing, $\rhobar(0+)=\infty$, $\rhobar
(+\infty)=0$, and $\int_0^1y\rhobar(y)\,\rmd y<\infty$. So $\rho$ is
a L\'evy measure and we can define a L\'evy process $(U_t)_{t\ge0}$,
independent of $(X_t)_{t\ge0}$, having L\'evy characteristics
$(0,0,\rho)$ and jump process $\Delta U_t:=U_t-U_{t-}$, $t>0$. Then
\[
P \Bigl(\sup_{0<s\le t} \Delta U_s\le y
\Bigr)=e^{-t\rhobar(y)},\qquad y>0,
\]
and the right-hand side of \eqref{Lapx} is
\begin{eqnarray*}
\int_{(0,\infty)} e^{-t\rhobar(y)}P \biggl(\sum
_{0<s\le t}\Delta X_s^+\in\lambda\,\rmd y \biggr) &=&
\int_{(0,\infty)} P \Bigl(\sup_{0<s\le t} \Delta
U_s\le y \Bigr) P \biggl(\sum_{0<s\le t}
\Delta X_s^+\in\lambda\,\rmd y \biggr)
\\
&=& P \biggl(\sup
_{0<s\le t} \Delta U_s\Big /\sum
_{0<s\le t}\Delta X_s^+\le\lambda^{-1}
\biggr).
\end{eqnarray*}
Thus \eqref{RS3} holds if and only if
\[
\frac{\sup_{0<s\le t}\Delta U_s}{\sum_{0<s\le t}\Delta X_s^+}
\buildrel P \over\to 0,\qquad\mbox{as }t\downarrow0.
\]
Applying \eqref{Rs1}, with $U_t$ in the role of the negative jump
process, this is so if and only if
%
%e4.11 #&#
\begin{equation}
\label{RS5} \lim_{y\downarrow0}\frac{y\rhobar(y)}{A_+(y)}=0.
\end{equation}
The estimates in \eqref{uy-} give
\[
\frac{A_-(y)}{3y}\le\rhobar(y) \le\frac{A_-(y)}{y},
\]
so the equivalence in \eqref{Rrs3} follows from \eqref{RS5}.

Finally, for the equivalence in \eqref{spm}, use
\[
P \Bigl(\sup_{0<s\le t}\Delta X_s^->u\sup
_{0<s\le t}\Delta X_s^+ \Bigr) =\int
_{(0,\infty)} \bigl(1-e^{-t\pibar^-(uy)} \bigr)\lambda
_t^+(\rmd y),
\]
and similar calculations as above.
\end{pf}

%
%\begin{remark}\label{rem5} {\rm
%The equivalence in \eqref{Rrs3} can be compared with an a.s. version
%in Bertoin (1997): when $X$ is of bounded variation,
%\begin{equation}\label{Rrs4}
%\lim_{t\dto0} \frac{\sum_{0<s\le t}\Delta X_s^-}{\sum_{0<s\le t}
%\Delta X_s^+}= 0, %{\rm a.s.,}
%\quad{\rm if and only if} \quad
%\int_0^1 \frac{\pibar^-(x)\rmd x}{\int_0^x\pibar^+(y)\rmd y}<\infty.
%\end{equation}
%We hope to explore a.s versions of the other relations in Theorem
%\ref{+-} elsewhere.
%} \end{remark}

To complete this section, we give the deferred proof of Theorem \ref{pmtX}.

\begin{pf*}{Proof of Theorem \ref{pmtX}} Assume $\pibar^+(0+)=\infty$ and
suppose first that \eqref{pmrlim}
holds. Then
%%$^{(1)}X_t/{\Delta X}_t^{(1)} \topr0$, as $t\dto0$, so by
%\eqref{t4a},
%%\begin{equation}\label{pmslim2}
%%\frac{\right\vert{\Delta X}_t^{(2)}\left\vert}{\right\vert{\Delta
%X}_t^{(1)}\left\vert}\topr0,\qquad\mbox{as }t
%\dto0.
%%\end{equation}
the same proof as used for showing that \eqref{rlim} is equivalent to
$\pibar(x)\in SV$ at 0, shows here that $\pibar^+(x)\in SV$ at 0.
This implies that the L\'evy process $ \sum_{0<s\le t} \Delta X^+_s$
is of bounded variation, and so
%
%e4.12 #&#
\begin{equation}
\label{+su} \frac{\sum_{0<s\le t} \Delta X^+_s}{\sup_{0<s\le t} \Delta
X^+_s} \buildrel P \over\to 1,
\end{equation}
by Theorem \ref{tX} applied to $ \sum_{0<s\le t} \Delta X^+_s$.
Now, $\pibar^+(x) \in SV$ at $0$ implies $\int_0^1\pibar^+(y)\,\rmd
y<\infty$ hence
$\lim_{x\downarrow0}x\pibar^+(x)=0$. This means
\[
P \Bigl(\sup_{0<s\le t} \Delta X^+_s>\delta t \Bigr)
\le1-e^{-
t\pibar^+(\delta t)}\to0,\qquad\mbox{as }t\downarrow0 \mbox{ for }
\delta>0,
\]
thus $\sup_{0<s\le t} \Delta X^+_s/t \buildrel P \over\to 0$.
So by \eqref{pmrlim}
%
%e4.13 #&#
\begin{equation}
\label{Xpm} \frac{X_t}{t}= \frac{X_t}{{\Delta X}_t^{(1)}} \cdot\frac
{{\Delta X}_t^{(1)}}{t}
\buildrel P \over\to 0.
\end{equation}
Then $\sigma^2=0$ and $x\pibar(x)\to0$ as $x\downarrow0$, by Doney and
Maller (\cite{doneymallera}, Theorem 2.1).

Use the L\'evy--It\^o decomposition \eqref{decomp3} (with $\sigma
^2=0$) to write $X_t$ as
%
%e4.14 #&#
\begin{eqnarray}
\label{Xpp} X_t&=&\gamma t+ \mbox{a.s. } \lim
_{\varepsilon\downarrow0} \biggl(\sum_{0<s\le t}\Delta
X_s \bfeins_{\{\veps<\llvert \Delta X_s \rrvert \le
1\}} -t\int_{\veps<\llvert x\rrvert \le1} x\Pi(
\rmd x) \biggr)+X_t^{(B,1)}
\nonumber\\[-8pt]\\[-8pt]\nonumber
&=& \gamma t+\sum
_{0<s\le t}\Delta X_s^+-t\int_{0<x<1}
x\Pi(\mathrm{d}x) -X_t^{(-)}+o_P(t).
\end{eqnarray}
Here,
\[
X_t^{(B,1)} = \sum_{0<s\le t }\Delta
X_s \bfeins_{\{\llvert \Delta X_s\rrvert > 1\}
} = o_P(t),\qquad\mbox{as }t
\downarrow0,
\]
because $
P (\llvert X_t^{(B,1)}\rrvert >\delta t ) \le1-e^{- t\pibar(1)}\to0$, as
$ t\downarrow0$, for $\delta>0$, and
\[
X_t^{(-)}:= \mbox{a.s. } \lim_{\varepsilon\downarrow0}
\biggl(\sum_{0<s\le t}\Delta X_s^-
\bfeins_{\{\veps<\Delta X_s^-\le
1\}} -t\int_{-1\le x<-\veps} x\Pi(\mathrm{d}x)
\biggr).
\]
In view of \eqref{Xpm} and \eqref{Xpp}, we see that $X_t^{(-)}/t$ has
a finite limit in probability as $t\downarrow0$,
and so by Doney and Maller \cite{doneymallera}, Theorem 2.1 (see also
Doney \cite{doney2004}), the integral
$\int_{(x,1)}y \Pi^{(-)}(\rmd y)$ has a finite limit as $x\downarrow0$.
This means that $X_t^{(-)}$, and hence $X_t$ are of bounded variation,
with drift $\rmd_X=0$ by \eqref{Xpm} and Lemma 4.1 of Doney and
Maller \cite{doneymallera}.

So we can write
%
%e4.15 #&#
\begin{eqnarray}\label{sX}
1+o_P(1) &=& \frac{X_t}{\sup_{0<s\le t} \Delta X^+_s}= \frac
{\sum_{0<s\le t} \Delta X^+_s-\sum_{0<s\le t} \Delta
X^-_s}{\sup_{0<s\le t} \Delta X^+_s}
\nonumber\\[-8pt]\\[-8pt]\nonumber
&=& 1-\frac{\sum_{0<s\le t} \Delta X^-_s}{\sup_{0<s\le t} \Delta X^+_s}+o_P(1).
\end{eqnarray}
From this, we see that
%
%e4.16 #&#
\begin{equation}
\label{+-} \frac{\sum_{0<s\le t} \Delta X^-_s}{\sup_{0<s\le t} \Delta
X^+_s} \buildrel P \over\to 0,
\end{equation}
thus by \eqref{pm}
\[
\lim_{x\downarrow0}\frac{\int_0^x\pibar^-(y)\,\rmd y}{x\pibar
^+(x)}= 0.
\]
Since $\int_0^x \pibar^-(y) \,\rmd y \ge x\pibar^-(x)$, we have $\lim
_{x\downarrow0} \pibar^-(x)/\pibar^+(x)=0$, so we have proved the forward
part of Theorem \ref{pmtX}.

For the converse, assume $\pibar^+(x)\in SV$ at 0, $X$ is of bounded
variation with drift $\rmd_X=0$, and $\lim_{x\downarrow0} \pibar
^-(x)/\pibar^+(x)=0$.
Now $\pibar^+(x)\in SV$ at 0 implies \eqref{+su} by Theorem \ref
{tX}, and also
$\int_0^x\pibar^+(y)\,\rmd y \sim x\pibar^+(x)$ as $x\downarrow0$.
In addition, $\pibar^-(x)= o(\pibar^+(x))$ implies
\[
\int_0^x\pibar^-(y)\,\rmd y =o \biggl(\int
_0^x\pibar^+(y)\,\rmd y \biggr)= o\bigl(x
\pibar^+(x)\bigr),\qquad\mbox{as }x\downarrow0,
\]
and then \eqref{+-} follows as in \eqref{pm}. Thus, we get \eqref
{pmrlim} from \eqref{sX}.
\end{pf*}

%s5 #&#
\section{$X$ dominating its large jump processes}\label{s4}
In this section, we characterise divergences
like\footnote{Recall that $\pibar(0+)=\infty$ implies $\sup_{0<s\le
t}\llvert \Delta X_s\rrvert >0$ a.s. for all $t>0$ when writing
ratios like that in
\eqref{01}, and similarly for one-sided versions.}
%
%e5.1 #&#
\begin{equation}
\label{01} \frac{X_t}{\sup_{0<s\le t}\llvert \Delta X_s\rrvert
} \buildrel P \over\to
\infty,\qquad\mbox{as }t
\downarrow0;
\end{equation}
and similarly with $\llvert \Delta X_s\rrvert $ replaced by $\Delta X_s$.
We think of these kinds of conditions as expressing the ``dominance''
of $X$ over its largest jump processes, at small times.

These conditions will be shown to be related to the \textit{relative
stability} of the process $X$, and to its \textit{attraction to
normality}, as $t\downarrow0$. Relative stability is the convergence
of the
normed process to a finite nonzero constant which, by rescaling of the
norming function, can be taken as $\pm1$. Thus, we are concerned with
the property
%
%e5.2 #&#
\begin{equation}
\label{00} \frac{X_t}{b_t} \buildrel P \over\to \pm1,\qquad\mbox{as }t
\downarrow0,
\end{equation}
where $b_t>0$ is a nonstochastic function.
The concept is important in a variety of contexts, in particular, with
reference to the stability at 0 of certain passage times for the
process, as we discuss in more detail later. When $X_t$ is replaced by
$\llvert X_t\rrvert $ in \eqref{01}, we also bring into play the idea
of $X$ being
in \textit{the domain of attraction of the normal distribution}, as
$t\downarrow
0$; that is, when there are nonstochastic functions $a_t\in\myRR$,
$b_t>0$, such that $(X_t-a_t)/b_t\stackrel{\mathrm
{D}}{\longrightarrow}N(0,1)$, a standard normal
random variable, as $t\downarrow0$.

Before proceeding, we quote some preliminary results, including in the
next subsection a theorem originally due to Doney \cite{doney2004}
giving necessary and sufficient conditions for $X_t$ to stay positive
with probability approaching 1 as $t\downarrow0$.
The main result concerning relative stability is in
Section~\ref{sub42}, while Section~\ref{sub43} deals with 2-sided
versions. The domain of attraction of the normal is needed here.
Subsequential versions of the results are in Sections~\ref{sub44} and
\ref{sub45}.

%s5.1 #&#
\subsection{$X$ staying positive near 0, in probability}\label{sub41}
%Recall that $\Pi^{(-)}$ denotes the restriction of the measure $\Pi$
%to $(-\infty,0)$.
Versions of truncated first and second moment functions, we will use are
%
%
%e5.3 #&#
\begin{equation}
\label{nudef} \nu(x)=\gamma-\int_{x<\llvert y\rrvert \le1}y\Pi(\mathrm{d}y)
\quad\mbox{and}\quad V(x)=\sigma^{2}+\int_{0<\llvert y\rrvert \le x}y^{2}
\Pi(\mathrm{d}y),\qquad x>0.
\end{equation}
%
%These are right continuous functions on $(0,\infty)$. The function
%$\nu(x)$ may be unbounded as $x\to0+$, but $V(x)$ is always finite.
%It is nondecreasing with $V(0+)=\sigma^2\in[0,\infty)$, and $V(+
%\infty)<\infty$
%if and only if $EX^2<\infty$.
%We will also need the functions
% \begin{equation}\label{nupmdef}
%\nu_+(x)=\int_{(x,1]} y \Pi(\rmd y)
%\quad{\rm and}\quad
%\nu_-(x)=\int_{[-1,-x)} \left\vert y\right\vert\Pi(\rmd y), x>0.
%\end{equation}
%Both $\nu_+(x)$ and $\nu_-(x)$ are non-negative nonincreasing
%right continuous functions
%on $(0,1]$, and we have
%\[
%\nu(x)=\gamma-\nu_+(x)+\nu_-(x).
%\]
%Also set
%\begin{equation}\label{Vpmdef}
%V_+(x)=\int_{0<y\le x} y^2 \Pi(\rmd y)
%\quad{\rm and}\quad
%V_-(x)=\int_{-x\le y<0} y^2 \Pi(\rmd y), x>0,
%\end{equation}
%so that
%\[
%V(x)=\sigma^2+V_+(x)+V_-(x).
%\]
Variants of $\nu(x)$ and $V(x)$ are Winsorised first and second
moment functions defined by
%
%
%e5.4 #&#
\begin{equation}
\label{Adef} A(x)=\gamma+\overline{\Pi}^+(1) -\overline{\Pi}^-(1) -\int
_x^1\bigl(\overline{\Pi}^+(y) -\overline{
\Pi}^-(y)\bigr)\,\rmd y
\end{equation}
and
%
%
%e5.5 #&#
\begin{equation}
\label{Udef} U(x)=\sigma^2 +2\int_0^x
y\overline{\Pi}(y)\,\rmd y\qquad\mbox{for } x>0.
\end{equation}
$A(x)$ and $U(x)$ are continuous for $x>0$. Using Fubini's theorem, we
can show that
%
%
%e5.6 #&#
\begin{equation}
\label{Apart} A(x)=\nu(x)+x\bigl(\pibar^+(x)-\pibar^-(x)\bigr)
\end{equation}
and
%
%
%e5.7 #&#
\begin{equation}
\label{Upart} U(x)=V(x)+x^2\bigl(\pibar^+(x)+\pibar^-(x)\bigr) =
V(x)+x^2\pibar(x).
\end{equation}
These functions are finite for all $x>0$ by virtue of
property $\int_{0<\llvert y\rrvert \le1}y^2\Pi(\rmd y)<\infty$ of the
L\'evy
measure $\Pi$,
which further implies that
$\lim_{x\downarrow0}x^2\overline{\Pi}(x)=0$,
and, as is easily verified,
%
%
%e5.8 #&#
\begin{equation}
\label{xnu} \lim_{x\downarrow0}x\nu(x)=\lim_{x\downarrow0}xA(x)=0.
%%\lim_{x\downarrow0}x\nu_\pm(x)=0.
\end{equation}
Also, $\lim_{x\to\infty}A(x)/x=\lim_{x\to\infty}U(x)/x^2=0$.
We have the obvious inequality
\[
%%\label{Uin}
U(x)\ge\sigma^2+x^2\pibar(x)\ge
x^2\pibar(x),\qquad x\ge0.
\]
This can be amplified to
%
%e5.9 #&#
\begin{equation}
\label{Upart1} U(x)\ge\sigma^2+x^2\pibar(x-)\ge
x^2\pibar(x-),\qquad x>0.
\end{equation}
%
%%Note also that
%%\bea\label{Upart1}U(x)&=&
%%\sigma^2 +\int_{0<\left\vert y\right\vert<x}y^2\Pi(\rmd y)+x^2(\Delta
%\Pi(x) + \Delta
%\Pi(-%%x))+x^2\pibar(x)\nonumber\\
%%&=&\sigma^2 +\int_{0<\left\vert y\right\vert<x}y^2\Pi(\rmd y)+x^2
%\pibar(x-)
%%\ge x^2\pibar(x-).\eea
Another calculation gives (recall~$\Delta\Pi(x)=\Pi\{x\}$)
%
%e5.10 #&#
\begin{eqnarray}
\label{Upart2} \nu(x)-x\bigl(\Delta\Pi(x)-\Delta\Pi(-x)\bigr) = A(x)-
x\bigl(
\pibar^+(x-)-\pibar^-(x-)\bigr).
\end{eqnarray}

%
%le5 #&#
\begin{lemma}\label{sigma2}
Suppose $\sigma^2>0$. Then $X_t/\sqrt{t} \stackrel{\mathrm
{D}}{\longrightarrow}N(0,\sigma^2)$ and
$P(X_t> 0)\to1/2$, as $t\downarrow0$.
\end{lemma}

\begin{pf}%{Proof of Lemma \ref{sigma2}}
The asymptotic normality of
$X_t/\sqrt{t}$ when $\sigma^2>0$ is
proved in
Doney and Maller (\cite{doneymallera}, Theorem 2.5 and its corollary, page~760), and then
$\lim_{t\downarrow0}P(X_t> 0)=1/2$ is immediate.
\end{pf}

Next, we quote the (slightly modified) theorem originally due to Doney
\cite{doney2004}. It shows that $X$ remains positive with probability
approaching 1 iff $X$ dominates its large negative jumps, and explicit
equivalences for this are given in terms of the functions $A(x)$,
$U(x)$ and the negative tail of $\Pi$.
The latter conditions reflect the positivity of $X$ in that the
function $A(x)$ remains positive for small values of $x$; and $A(x)$
dominates $U(x)$ and the negative tail of $\Pi$ in certain ways.
Recall the notation $\Delta X_t^+=\max(\Delta X_t,0)$, $\Delta
X_t^-=\max(-\Delta X_t,0)$, and
$(\Delta X^+)_t^{(1)}=\sup_{0<s\le t}\Delta X_s^+$, $(\Delta
X^-)_t^{(1)}=\sup_{0<s\le t}\Delta X_s^-$.

%
%th5 #&#
\begin{theorem}\label{th33}
Suppose $\pibar^+(0+)=\infty$.

\begin{longlist}[(ii)]
\item[(i)] Suppose also that $\pibar^-(0+)>0$.
Then the following are equivalent:
%
%
%e5.11 #&#
\begin{eqnarray}
\label{37} \lim_{t\downarrow0}P(X_t> 0) &=&1;
\\
%e5.12 #&#
\label{Xdel} \frac{X_t}{(\Delta X^-)_t^{(1)}}& \buildrel P \over\to
&\infty,\qquad\mbox{as } t \downarrow 0;
\\
%
%
%e5.13 #&#
\label{36} \sigma^2&=&0\quad\mbox{and}\quad\lim
_{x\downarrow0} \frac{A(x)}{x\pibar^-(x)}=\infty;
\\
%
%
%e5.14 #&#
\label{36a} \lim_{x\downarrow0} \frac{A(x)}{\sqrt{U(x)\pibar
^-(x)}}&=&\infty;
\end{eqnarray}
there is a nonstochastic nondecreasing function $\ell(x)>0$, which is
slowly varying at 0, such that
%
%
%e5.15 #&#
\begin{equation}
\label{39} \frac{X_t}{t\ell(t)} \buildrel P \over\to \infty,\qquad
\mbox{as } t
\downarrow0.
\end{equation}
\item[(ii)] Suppose $X$ is spectrally positive, so $\pibar^-(x)=0$ for
$x>0$. Then (\ref{37}) is equivalent to
%
%
%e5.16 #&#
\begin{equation}
\sigma^{2}=0\quad\mbox{and}\quad A(x)\ge0\qquad\mbox{for all small }x,
\label{38}
\end{equation}
and this happens if and only if $X$ is a subordinator.
Furthermore, we then have $A(x)\geq0$, not only for small
$x$, but for all $x>0$.
\end{longlist}
\end{theorem}

%
%re2 #&#
\begin{remark}\label{rem2}
We adopt the convention that \eqref{Xdel} is taken to hold when
\eqref{37} holds but
$\sup_{0<s\le t} \Delta X_s^-=0$ a.s. for all small $t>0$.
This is the case when $\pibar^-(0+)<\infty$.
\end{remark}

%
%le6 #&#
\begin{lemma}\label{sig1}
If $\pibar^-(0+)>0$, then\vspace*{-2pt}
%
%
%e5.17 #&#
\begin{equation}
\label{41x} \limsup_{x\downarrow0}\frac{A(x)}{\sqrt{\pibar^-(x)}}<\infty.
\end{equation}
If $\pibar^-(0+)=0$ and $\pibar^+(0+)>0$, then
%
%
%e5.18 #&#
\begin{equation}
\label{41y} \limsup_{x\downarrow0}\frac{A(x)}{\sqrt{\pibar^+(x)}}<\infty.
\end{equation}
\end{lemma}

\begin{pf*}{Proof of Lemma \ref{sig1}}
(i)~Assume~$\pibar^-(0+)>0$ and, by way of contradiction, that there
is a nonstochastic sequence $x_k\downarrow0$ as $k\to\infty$ such that
\[
\frac{A(x_k)}{\sqrt{\pibar^-(x_k)}}= \frac{\gamma+\pibar^+(1)-\pibar
^-(1)-\int_{x_k}^1\pibar^+(y)\,\rmd y
+\int_{x_k}^1\pibar^-(y)\,\rmd y}{\sqrt{\pibar^-(x_k)}} \to\infty.
\]
Since $\pibar^-(0+)>0$, we deduce from this that
\[
\frac{-\pibar^-(1)+\int_{x_k}^1\pibar^-(y) \,\rmd y}{
\sqrt{\pibar^-(x_k)}} \to\infty.
\]
Thus, integrating by parts,
\[
\frac{-x_k\pibar^-(x_k)+
\int_{x_k<y\le1}y \Pi^{(-)}(\rmd y)}{\sqrt{\pibar^-(x_k)}} \to\infty.
\]
But
%%$x_k\pibar^-(x_k)/\sqrt{\pibar^-(x_k)}\le
%%\sqrt{x_k^2\pibar^-(x_k)}\to0$. Also,
by the Cauchy--Schwarz inequality,
\begin{eqnarray*}
\frac{ (\int_{x_k<y\le1} y\Pi^{(-)}(\rmd y) )^2}{\pibar^-(x_k)} &\le&
\frac{\int_{x_k<y\le1}
y^2 \Pi^{(-)}(\rmd y)\int_{x_k<y\le1}\Pi^{(-)}(\rmd y)}{\pibar^-(x_k)}
\nonumber
\\[-2pt]
&\le& \frac{\int_{x_k<\llvert y\rrvert \le1} y^2 \Pi(\rmd y)
(\pibar^-(x_k)-\pibar^-(1) )}{\pibar^-(x_k)}
\nonumber
\le \int_{0<\llvert y\rrvert \le1} y^2 \Pi(\rmd y)<\infty,
\end{eqnarray*}
giving a contradiction. Thus,~\eqref{41x} holds.

(ii) Alternatively, suppose $\pibar^-(0+)=0$ and $\pibar
^+(0+)>0$. Then, for $0<x<1$,
\[
\frac{A(x)}{\sqrt{\pibar^+(x)}} = \frac{\gamma+\pibar^+(1)-\int
_{x}^1\pibar^+(y)\,\rmd y}{\sqrt{\pibar^+(x)}} \le\frac{\gamma+\pibar
(1)}{\sqrt{\pibar^+(x)}},
%%\quad\dto\quad\frac{\gamma+\pibar(1)}{\sqrt{\pibar^+(0+)}}<\infty
%, %%\quad\mbox{as } x \dto0
\]
and since $\pibar^+(0+)>0$ the RHS is finite as $x\downarrow0$, so
\eqref
{41y} is proved.
\end{pf*}

\begin{pf}%{Proof of Theorem \ref{th33}}
Theorem \ref{th33} only differs from Theorem 1 in Doney \cite
{doney2004} (and his remark following the theorem, regarding part~(ii)
of our Theorem \ref{th33}) in that he assumes {a priori} that
$\sigma^2=0$. Clearly, \eqref{37}, \eqref{Xdel} and \eqref{39}
imply this by Lemma \ref{sigma2}.
\eqref{36a} also implies $\sigma^2=0$. To see this, suppose on the
contrary that $\sigma^2 > 0$. Then $U(x)\ge\sigma^2$ for all $x\ge
0$ and by Lemma \ref{sig1}, \eqref{41x} contradicts \eqref{36a}.
\end{pf}

We have the following subsequential version of Theorem \ref{th33}.
We omit the proof which is along the lines of Doney's proof, together
with similar ideas as in Theorem \ref{SPRS}.

%
%th6 #&#
\begin{theorem}\label{th33s}
Suppose $\pibar^+(0+)=\infty$.

\begin{longlist}[(ii)]
\item[(i)] Suppose also that $\pibar^-(0+)>0$. Then the following are equivalent:
there is a nonstochastic sequence $t_k\downarrow0$ such that
%
%
%e5.19 #&#
\begin{equation}
\label{37s} P(X_{t_k}>0)\to1; %% \mathrm{as} k\to\infty;
\end{equation}
there is a nonstochastic sequence $t_k\downarrow0$ such that
%
%
%e5.20 #&#
\begin{eqnarray}
\label{Xdels}
\frac{X_{t_k}}{ (\Delta X^-)_{t_k}^{(1)}}& \buildrel P
\over\to &\infty,\qquad\mbox{as } k
\to\infty;
\\
%
%
%e5.21 #&#
\label{36as} \limsup_{x\downarrow0} \frac{A(x)}{\sqrt{U(x)\pibar
^-(x)}}&=&\infty.
\end{eqnarray}

\item[(ii)] Suppose $X$ is spectrally positive, that is, $\pibar
^-(x)=0$ for
all $x>0$. Then (\ref{37s}) is equivalent to
$\lim_{t\downarrow0} P(X_t>0)\to1$, thus to \eqref{38},
equivalently, $X_t$ is a subordinator, and $A(x)\geq0$ for all $x>0$.
\end{longlist}
%
%(iii) Suppose $\pibar^-(0+)>0$. Then
%$X_{t_k}/t_k\topr\infty$ for a nonstochastic sequence $t_k\dto0$
%if and only if
%\begin{equation} \label{36as+}
%\limsup_{x\dto0}
% \frac{A(x)}{1+\sqrt{U(x)\pibar^-(x)}}=\infty.
%%\end{equation}
\end{theorem}

%
%re3 #&#
\begin{remark}\label{rem3}
We get equivalences for
\[
\frac{X_t}{(\Delta X^+)_t^{(1)}} \buildrel P \over\to -\infty
\]
(or the subsequential version)
by applying Theorem \ref{th33} (or Theorem \ref{th33s}) with $X$
replaced by $-X$.
\end{remark}

In the next two subsections, we characterise when $X$ dominates its
large positive jumps
and its jumps large in modulus, while remaining positive in
probability, and
when $\llvert X\rrvert $ dominates its jumps large in modulus.
These kinds of behaviour require more stringent conditions on $X$,
namely, relative stability or attraction to normality, in the
respective cases.

%s5.2 #&#
\subsection{Relative stability and dominance}\label{sub42}
Recall that $X$ is said to be \textit{relatively stable} (RS) at 0 if
\eqref{00} holds.
$X$ is \textit{positively relatively stable} (PRS) at 0 if \eqref{00}
holds with a ``$+$'' sign, and
\textit{negatively relatively stable} (NRS)
at 0 if \eqref{00} holds with a ``$-$'' sign.
In either case, the function $b_t>0$ is regularly varying at 0 with
index 1.
In Griffin and Maller (\cite{GM2}, Proposition 2.1) it is shown (when
$\pibar(0+)=\infty$) that there is a measurable nonstochastic
function $b_t>0$ such that
%
%e5.22 #&#
\begin{equation}
\label{rs00} \frac{\llvert X_t\rrvert }{b_t} \buildrel P \over\to
1,\qquad\mbox{as }t
\downarrow0,
\end{equation}
iff $X\in RS$ at 0, equivalently, iff
%
%e5.23 #&#
\begin{equation}
\label{prs2a} \sigma^{2}=0\quad\mbox{and}\quad\lim_{x\downarrow0}
\frac{\llvert A(x)\rrvert }{x\pibar(x)} =\infty.
\end{equation}
The following conditions characterise the convergence in \eqref{00}
(Kallenberg \cite{Kall}, Theorem~15.14): for all $x>0$,
%
%e5.24 #&#
\begin{equation}
\label{degcon} \lim_{t\downarrow0} t\pibar(xb_t)=0,\qquad
\lim_{t\downarrow0} \frac{tA(xb_t)}{b_t}=\pm1, \qquad
\lim_{t\downarrow0}
\frac{tU(xb_t)}{b_t^2}=0.
\end{equation}
Obvious modifications of these characterise convergence through a
subsequence $t_k$ in \eqref{00}.

Next is our main result relating ``one-sided'' dominance to positive
relative stability.
The identity \eqref{cor1a} supplies a key step in the proof.

%
%th7 #&#
\begin{theorem}\label{PRS}
Assume $\pibar^+(0+)=\infty$. Then the following are equivalent:
%
%e5.25 #&#
\begin{eqnarray}
\label{prs4a}
\frac{X_t}{(\Delta X^+)_t^{(1)}}& \buildrel P \over\to
&\infty,\qquad\mbox{as }t
\downarrow0;
\\
%
%e5.26 #&#
\label{prs4} \frac{X_t}{\llvert \widetilde{\Delta X}_t^{(1)}\rrvert
}& \buildrel P \over\to
&\infty,\qquad\mbox{as }t
\downarrow0;
\\
%
%e5.27 #&#
\label{prs2} \sigma^{2}&=&0\quad\mbox{and}\quad\lim_{x\downarrow0}
\frac{A(x)}{x\pibar(x)} =\infty;
\\
%
%e5.28 #&#
\label{prs1} X&\in&\mathrm{PRS}\qquad\mbox{at } 0;
\\
%
%e5.29 #&#
\label{prs3} \lim_{x\downarrow0} \frac{A(x)}{\sqrt{U(x)\pibar
(x)}}&=&\infty;
\\
%
%e5.30 #&#
\label{prs3a} \lim_{x\downarrow0} \frac{xA(x)}{U(x)}&=&
\infty.
\end{eqnarray}
%
%\begin{enumerate}[\rm(a)]
%
%\item\label{prs4a}
%$\displaystyle\frac{X_t}{(\Delta X^+)_t^{(1)}}\topr\infty, {\rm as}
% t\dto0;$
%
%\item\label{prs4}
%$\displaystyle
%\frac{X_t}{\left\vert\wt{\Delta X}_t^{(1)}\right\vert}\topr\infty, {
%\rm as} t\dto0; $
%
%\item\label{prs2}
%$\displaystyle\sigma^{2}=0 \mathrm{and} \lim_{x\dto0}
%\frac{A(x)}{x\pibar(x)} =\infty;$
%
%\item\label{prs1}
%$\displaystyle X\in{\rm PRS at} 0;$
%
%\item\label{prs3}
%$\displaystyle\lim_{x\dto0} \frac{A(x)}{\sqrt{U(x)\pibar(x)}}=
%\infty; $
%
%\item\label{prs3a}
%$\displaystyle\lim_{x\dto0} \frac{xA(x)}{U(x)}=\infty. $
%
%\end{enumerate}
\end{theorem}

Before proving the theorem, we record the following moment formulae.
Recall that $\widetilde{ X}_t^{v}$
is defined in \eqref{21b}.
%
%le7 #&#
\begin{lemma}\label{mom}
When $\overline{\Pi}^{\leftarrow}(v)<1$ and $t>0$:
%
%e5.31 #&#
\begin{eqnarray}
\label{mom1} t^{-1}E\widetilde{ X}_t^{v} &=&
\nu\bigl(\overline{\Pi}^{\leftarrow}(v)\bigr)-\overline{\Pi}^{\leftarrow
}(v)
\bigl( \Delta\Pi\bigl(\overline{\Pi}^{\leftarrow}(v)\bigr)-\Delta\Pi
\bigl(-
\overline{\Pi}^{\leftarrow}(v)\bigr) \bigr)
\nonumber\\[-8pt]\\[-8pt]\nonumber
&=& A\bigl(\overline{\Pi}^{\leftarrow}(v)\bigr)-\overline{
\Pi}^{\leftarrow
}(v) \bigl(\pibar^+\bigl(\overline{\Pi}^{\leftarrow}(v)-\bigr)
-\pibar^-\bigl(\overline{\Pi}^{\leftarrow}(v)-\bigr) \bigr).
\end{eqnarray}
%
%and when $\pibarinv(v)> 1$:
%\be\label{mom2}
%E\wt{ X}_t^{v}
%=t\left(\gamma+\int_{1<\left\vert x\right\vert< \pibarinv(v)}x\Pi(
%\rmd x)\right).
%\ee
For all $t>0$, $v>0$,
%
%e5.32 #&#
\begin{equation}
\label{mom3} E\bigl(\widetilde{ X}_t^{v}
\bigr)^2=t \biggl(\sigma^2+ \int_{\llvert x\rrvert <\overline{\Pi
}^{\leftarrow}(v)}x^2
\Pi(\rmd x) \biggr) + \bigl(E\widetilde{ X}_t^{v}
\bigr)^2.
\end{equation}
\end{lemma}

\begin{pf}%{Proof of Lemma \ref{mom}}
Let $(U_t)_{t\ge0}$ be a L\'evy process with triplet
$(\gamma_U,\sigma^2_U, \Pi_U)$.
Provided the participating integrals are finite (see Example~25.11 in
Sato \cite{sato1}), for instance, we have
\[
EU_t=t \biggl(\gamma_U+\int_{\llvert y\rrvert >1}y
\Pi_U(\rmd y) \biggr) \quad\mbox{and}\quad E(U_t)^2=t
\biggl(\sigma^2_U+\int_{\myRR_*}y^2
\Pi_U(\rmd y) \biggr)+(EU_t)^2.
\]
Apply these to $\widetilde{ X}_t^{v}$ with triplet as in \eqref
{trip2} to
get, when $\overline{\Pi}^{\leftarrow}(v)<1$ and $t>0$,
\begin{eqnarray*}
t^{-1}E\widetilde{ X}_t^{v} &=& \gamma-\int
_{\overline{\Pi}^{\leftarrow}(v)\le\llvert x\rrvert \le1}x\Pi(\rmd x)
\nonumber
\\
&=& \gamma-\int_{\overline{\Pi}^{\leftarrow}(v)<\llvert x\rrvert \le
1}x\Pi(\rmd x) -\overline{
\Pi}^{\leftarrow}(v) \bigl(\Delta\Pi\bigl(\overline{\Pi}^{\leftarrow}(v)
\bigr)-\Delta\Pi\bigl(-\overline{\Pi}^{\leftarrow} (v)\bigr) \bigr),
\end{eqnarray*}
which gives the first equation in \eqref{mom1}.
For the second equation in \eqref{mom1}, use \eqref{Upart2}.
\eqref{mom3} is proved similarly.
\end{pf}

\begin{pf*}{Proof of Theorem \ref{PRS}}
Assume $\pibar^+(0+)=\infty$ throughout.
%
%First note that any of~\eqref{prs4}--\eqref{prs1}
%
% If $\sigma^2>0$ \eqref{prs3} cannot hold by Lemma~\ref{sig1},
%while \eqref{prs3a} would imply $\lim_{x\downarrow0} x A(x)=\infty$,
%contradicting \eqref{xnu}.
%So we can assume $\sigma^2=0$ throughout.
%
%We prove \eqref{prs4} equivalent to \eqref{prs1} via \eqref{prs2}.

\begin{longlist}
\item[\textit{Case} (i).] Suppose $\pibar^-(0+)>0$.

\eqref{prs4a} $\Longrightarrow$ \eqref{prs4}: Assume
\eqref{prs4a}.
This implies $\lim_{t\downarrow0} P(X_t > 0)=1$, so by Theorem \ref
{th33}, \eqref{Xdel} holds.
\eqref{Xdel} together with \eqref{prs4a} implies \eqref{prs4}, because
$\llvert \widetilde{\Delta X}_t^{(1)}\rrvert =\max((\Delta
X^+)_t^{(1)}, (\Delta
X^-)_t^{(1)})$.

\eqref{prs4} $\Longrightarrow$ \eqref{prs2}: Assume \eqref{prs4}. Then
$\lim_{t\downarrow0} P(X_t > 0)=1$, so
\eqref{36} holds.
Since $\pibar^-(0+)>0$, \eqref{36} implies $\lim_{x\downarrow
0}A(x)/x=\infty$; in particular, $A(x)>0$ for all small $x$.
Since  $\lim_{t\downarrow0}P(X_t> 0)=1$,
Lemma 5 in Doney \cite{doney2004} gives
%
%e5.33 #&#
\begin{equation}
\label{prs10} U(x)\le3xA(x)\qquad\mbox{for all small } x, x\le x_0,\mbox
{ say}.
\end{equation}
Without loss of generality, assume $x_0<1$.

Note that \eqref{prs4} also implies
\[
\frac{^{(1)}\widetilde X_t}{\llvert \widetilde{\Delta
X}_t^{(1)}\rrvert }\buildrel
P \over\to \infty,\qquad\mbox{as }x\downarrow0
\]
(recall \eqref{trims}), so we have
%
%e5.34 #&#
\begin{equation}
\label{prs11} \lim_{t\downarrow0}P \bigl({}^{(1)}\widetilde
X_t\le a\bigl\llvert\widetilde{\Delta X}_t^{(1)}
\bigr\rrvert\bigr)=0\qquad\mbox{for some } a>0.
\end{equation}
(In fact, this holds for all $a>0$. But it will be enough to assume
\eqref{prs11}.)
Without loss of generality, take $a\le1$.

We will abbreviate $\overline{\Pi}^{\leftarrow}(v)$ to $y_v$
throughout this proof.
Then by \eqref{cor1a}, we can write
%
%e5.35 #&#
\begin{equation}
\label{prs12} P \bigl({}^{(1)}\widetilde X_t\le a\bigl
\llvert\widetilde{\Delta X}_t^{(1)}\bigr\rrvert\bigr) =\int
_0^\infty P \bigl(\widetilde{
X}_t^{v} + \widetilde{G}_t^{v} \le
ay_v \bigr)P(\EEEE\in t\,\rmd v),
\end{equation}
where $\EEEE=\SSSS_1$ is a unit exponential r.v.
By \eqref{Ppar} and \eqref{cor1aa},
%
%e5.36 #&#
\begin{eqnarray}
\label{ega} \bigl\llvert E\widetilde{G}_t^{v}\bigr
\rrvert&=& y_v\bigl\llvert EY^+_{t\kappa^+(v)}-EY^-_{t\kappa^-(v)}
\bigr\rrvert
\nonumber
\\
&\le& ty_v\bigl(\pibar(y_v-)-v\bigr) \frac{\Delta\Pi(y_v)+\Delta\Pi
(-y_v)}{\Delta\pibar(y_v)}
\bfeins_{\{\Delta\pibar(y_v)\ne0\}}
\nonumber\\[-8pt]\\[-8pt]\nonumber
&\le& ty_v\pibar(y_v-) \le tU(y_v)/y_v
\quad\bigl(\mbox{by \eqref{Upart1}}\bigr)
\nonumber
\\
&\le& 3tA(y_v)\qquad\bigl(\mbox{by \eqref{prs10}}\bigr)\nonumber
\end{eqnarray}
and similarly
%
%e5.37 #&#
\begin{equation}
\label{vgt} \operatorname{Var}\bigl(\widetilde{G}_t^{v}\bigr)
\le ty_v^2 \pibar(y_v-).
\end{equation}
With $x_0$ as in \eqref{prs10}, keep $v\ge\pibar(x_0)$, so $y_v\le x_0<1$.
%%Now $y_v=\pibarinv(v)\to0$ as $v\to0$, so we can assume so small
%that $y_v<1$.
Then
%
%e5.38 #&#
\begin{eqnarray}
\label{exa} E\widetilde{ X}_t^{v} &=& t
\bigl(A(y_v)-y_v \bigl(\pibar^+(y_v-) -
\pibar^-(y_v-) \bigr) \bigr) \qquad\bigl(\mbox{by \eqref{mom1}}\bigr)
\nonumber
\\
&\le& t \bigl(A(y_v)+y_v\pibar^-(y_v-) \bigr)
\\
&\le& 4tA(y_v)\qquad\bigl(\mbox{by \eqref{Upart1} and \eqref{prs10}}\bigr).\nonumber
\end{eqnarray}
Apply \eqref{ega} and \eqref{exa} to obtain from \eqref{prs12}
%
%e5.39 #&#
\begin{eqnarray}
\label{prs13} &&P \bigl({}^{(1)}\widetilde X_t\le a\bigl
\llvert\widetilde{\Delta X}_t^{(1)}\bigr\rrvert\bigr)
\nonumber\\[-8pt]\\[-8pt]\nonumber
&&\quad\ge\int_{\pibar(x_0)}^\infty P \bigl(\widetilde{
X}_t^{v}-E\widetilde{ X}_t^{v} +
\widetilde{G}_t^{v}-E\widetilde{ G}_t^{v}
\le ay_v-7tA(y_v) \bigr)P(\EEEE\in t\,\rmd v).
\end{eqnarray}

For $t>0$ and $a$ as in \eqref{prs11} define
%
%e5.40 #&#
\begin{equation}
\label{bdef} b_t:= \sup\biggl\{x>0: \frac{A(x)}{x} >
\frac{a^2}{56t} \biggr\},
\end{equation}
with $b_0:=0$. Recall that $\lim_{x\downarrow0}A(x)/x=\infty$, $\lim
_{x\to\infty}A(x)/x=0$, and $A(x)$ is continuous.
So $0<b_t<\infty$, $b_t$ is strictly increasing, $b_t\downarrow0$ as
$t\downarrow0$,
and
%
%e5.41 #&#
\begin{equation}
\label{prs14} \frac{tA(b_t)}{b_t}= \frac{a^2}{56}.
\end{equation}
Assume $t$ is small enough for $b_t\le x_0$ and keep $v< \pibar(b_t)$.
Then $y_v\ge b_t$, and so
$7tA(y_v)\le a^2y_v/8$ by definition of $b_t$.
This implies $7tA(y_v)\le ay_v/2$.
Thus, by Chebyshev's inequality and~\eqref{prs13}
%
%e5.42 #&#
\begin{eqnarray}
\label{prs14a} P \bigl({}^{(1)}\widetilde X_t\le a\bigl
\llvert\widetilde{\Delta X}_t^{(1)}\bigr\rrvert\bigr) &\ge&
\int_{\pibar(x_0)}^{\pibar(b_t)} P \bigl(\widetilde{
X}_t^{v}-E\widetilde{ X}_t^{v} +
\widetilde{G}_t^{v}-E\widetilde{ G}_t^{v}
\le ay_v/2 \bigr)P(\EEEE\in t\,\rmd v)
\nonumber\\[-8pt]\\[-8pt]\nonumber
&\ge& \int_{\pibar(x_0)}^{\pibar(b_t)} \biggl(1-
\frac{4 (\operatorname{Var}(\widetilde{X}_t^{v})+\operatorname
{Var}(\widetilde
{G}_t^{v}) )}{
a^2y_v^2} \biggr)P(\EEEE\in t\,\rmd v).
\end{eqnarray}
Also
%
%e5.43 #&#
\begin{eqnarray} \label{prs15}
\operatorname{Var}\bigl(\widetilde{X}_t^{v}\bigr)+\operatorname
{Var}\bigl(
\widetilde{G}_t^{v}\bigr) &\le& t \bigl(V(y_v)+y_v^2
\pibar(y_v-) \bigr) \qquad\bigl(\mbox{by \eqref{mom3} and \eqref{vgt}}\bigr)
\nonumber
\\
&\le& 2tU(y_v)\qquad\bigl(\mbox{see \eqref{Upart} and \eqref{Upart1}}\bigr)
\nonumber\\[-8pt]\\[-8pt]\nonumber
&\le& 6ty_vA(y_v)\qquad\bigl(\mbox{by \eqref{prs10}, since
}y_v\le x_0\bigr)
\nonumber
\\
&\le& a^2y_v^2/8.\nonumber
\end{eqnarray}
The last inequality holds since $y_v\ge b_t$. Hence, from \eqref{prs14a},
\begin{eqnarray*}
P \bigl({}^{(1)}\widetilde X_t\le a\bigl\llvert\widetilde{
\Delta X}_t^{(1)}\bigr\rrvert\bigr) &\ge& t\int
_{\pibar(x_0)}^{\pibar(b_t)}e^{-tv}\,\rmd v/2
\nonumber
\\
&=& e^{-t\pibar(x_0)} \bigl(1- e^{-t(\pibar(b_t)-\pibar(x_0))} \bigr)/2.
\end{eqnarray*}
Since the left-hand side tends to 0 as $t\downarrow0$ by \eqref
{prs11}, we
see from \eqref{prs14} that
%
%e5.44 #&#
\begin{equation}
\label{prs16} t\pibar(b_t) =\frac{a^2b_t\pibar(b_t)}{56A(b_t)}\to
0,\qquad\mbox{as }t
\downarrow0.
\end{equation}

Now take $\lambda>1$ and write, by \eqref{prs14},
\begin{eqnarray*}
\frac{b_{\lambda t}}{b_t} &=& \frac{56\lambda tA(b_{\lambda
t})}{a^2b_t} = \frac{56\lambda tA(b_{t})}{a^2b_t} +
\frac{56\lambda t (A(b_{\lambda t})-A(b_t) )}{a^2b_t}
\nonumber
\\
&=& \lambda+ \frac{56\lambda t\int_{b_t}^{b_{\lambda t}} (\pibar
^+(y)-\pibar^-(y) )\,\rmd y}{a^2b_t}
\nonumber
\\
&=& \lambda+ O \bigl(t\pibar(b_t) \bigr) \biggl( \frac{b_{\lambda t}-b_t}{b_t}
\biggr)
\nonumber
\\
&=& \lambda+ o \biggl( \frac{b_{\lambda t}}{b_t} \biggr).
\end{eqnarray*}

Thus, $b_t$ is regularly varying with index 1 as $t\downarrow0$.
Also,~\eqref{prs16} implies
$A(b_t)/b_t \pibar(b_t)\to\infty$ as $t\downarrow0$. From those we
obtain~\eqref{prs2} as follows. Given
$x>0$ choose $t=t(x)$ so that $b_{t-}\le x\le b_{t+}$. Then, for any
$\varepsilon\in(0,1)$, $b_{t(1-\varepsilon)}\le x\le
b_{t(1+\varepsilon)}$, while $b_{t(1+\varepsilon)}\sim(1+\varepsilon
) b_{t}\sim(1+\varepsilon) b_{(1-\varepsilon)t}/(1-\varepsilon)$ as
$t\downarrow0$. So
\begin{eqnarray*}
A(x)&=&A(b_{t(1-\varepsilon)})+\int_{b_{t(1-\varepsilon)}}^x\bigl(
\pibar^+(y)-\pibar^-(y)\bigr) \,\rmd y
\\
&\ge&A(b_{t(1-\varepsilon)})-b_{t(1+\varepsilon)} \pibar
(b_{t(1-\varepsilon)})
\\
&\ge&\bigl(1+o(1)\bigr)A(b_{t(1-\varepsilon)})\qquad\bigl(\mbox{by~\eqref{prs16}}\bigr).
\end{eqnarray*}
Hence, as $x\downarrow0$,
%
%
%e5.45 #&#
\begin{equation}
\label{3611} \frac{A(x)}{x\pibar(x)}\ge\frac{(1+o(1))A(b_{t(1-\varepsilon
)})}{b_{t(1-\varepsilon)}\pibar(b_{t(1-\varepsilon)})}\times\frac
{b_{t(1-\varepsilon)}}{b_{t(1+\varepsilon)}}\to
\infty,
\end{equation}
and \eqref{prs2} is proved.

\eqref{prs2} $\iff$ \eqref{prs1} is in Theorem 2.2 of Doney and
Maller \cite{doneymallera}.

\eqref{prs1} $\Longrightarrow$ \eqref{prs3}: \eqref{prs1} implies
$A(b_t)/\sqrt{U(b_t)\pibar(b_t)}\to\infty$ by \eqref{degcon} and then
\eqref{prs3} follows from the regular variation of $b_t$ (noted prior
to \eqref{rs00}),
by similar arguments as we used in proving \eqref{3611} from \eqref{prs16}.

\eqref{prs3} $\iff$ \eqref{prs3a}: \eqref{prs3} implies $X\in \mathrm{PRS}$ at
0, so
$b_tA(b_t)/U(b_t)\to\infty$ by \eqref{degcon}, and $b_t$ is
regularly varying with index 1 at 0. Then \eqref{prs3a} follows by
similar arguments as we used in proving \eqref{3611} from \eqref{prs16}.
Conversely, \eqref{prs3a} implies \eqref{prs3} because $U(x)\ge
x^2\pibar(x)$.

In the reverse direction, we will show that \eqref{prs3}
$\Longrightarrow$
\eqref{prs2} %%$\iff$ \eqref{prs1}
$\Longrightarrow$ \eqref{prs4} $\Longrightarrow$ \eqref{prs4a}.

\eqref{prs3} $\Longrightarrow$ \eqref{prs2}: \eqref{prs3} implies
\eqref
{36a}, hence $\sigma^2=0$ by Theorem \ref{th33}. Then \eqref{prs2}
follows from \eqref{prs3} since $U(x)\ge x^2\pibar(x)$.

\eqref{prs2} $\Longrightarrow$ \eqref{prs4}: Assume \eqref{prs2}.
This implies
$X\in \mathrm{PRS}$, so $X_t/b_t \buildrel P \over\to +1$ as $t\downarrow0$
for some $b_t>0$.
By \eqref{degcon}, $ \lim_{t\downarrow0} t\pibar(\veps b_t)=0$ for all
$\veps>0$. This implies
\[
%%\label{supn}
P \Bigl(\sup_{0<s\le t}\llvert\Delta X_s
\rrvert>\veps b_t \Bigr) =1- e^{-t\pibar(\veps b_t)} \to0,
\]
thus $\sup_{0<s\le t}\llvert \Delta X_s\rrvert /b_t \buildrel P \over
\to 0$ as
$t\downarrow0$.
So we get \eqref{prs4}.

\eqref{prs4} $\Longrightarrow$ \eqref{prs4a} is true since
$\llvert \widetilde{\Delta
X}_t^{(1)}\rrvert \ge(\Delta X^+)_t^{(1)}$.
So we have shown the equivalence of \eqref{prs4a}--\eqref{prs3a} for
case~(i).

\item[\textit{Case} (ii).] Suppose $\pibar^-(0+)=0$. By
part~(ii) of Theorem \ref{th33}, each of
\eqref{prs4a}--\eqref{prs1} implies $X$ is a subordinator (with
drift) and $A(x)\ge0$ for all $x\ge0$.
\eqref{prs4a} and \eqref{prs4} are the same thing in this case.

\eqref{prs4} $\Longrightarrow$ \eqref{prs2}:
Assume \eqref{prs4}. Since $X$ is a subordinator,
%%and $A(x)\ge0$ for all $x\ge0$,
we can write
\[
A(x)=\rmd_X+\int_0^x
\pibar^+(y)\,\rmd y,\qquad x\ge0,
\]
where $\rmd_X\ge0$ is the drift of $X$ and $\int_0^x\pibar^+(y)\,\rmd
y<\infty$.
The latter implies $\lim_{x\downarrow0}x\pibar^+(x)=0$. Of course
$\sigma
^2=0$ and if $\rmd_X>0$ then \eqref{prs2} clearly holds.
% \ben
% \frac{A(x)}{x\pibar(x)} \to\infty,
% \een
% so \eqref{prs2} holds.
So suppose $\rmd_X=0$.
As in \eqref{exa}, we get $E\widetilde{ X}_t^{v}\le tA(y_v)$ and
\eqref
{ega} and \eqref{vgt} remain true.
Since $\pibar^+(0+)=\infty$,
\[
\lim_{x\downarrow0}\frac{A(x)}{x}\ge\int_0^1
\liminf_{x\downarrow
0}\pibar^+(xy)\,\rmd y=\infty.
\]
Define $b_t$ again by \eqref{bdef}. Then the same working as in case~(i) gives
$t\pibar(b_t)\to0$ and $b_t$ regularly varying with index 1, so again
we get \eqref{prs2}.

\eqref{prs2} $\iff$ \eqref{prs1} is in Theorem 2.2 of Doney and
Maller \cite{doneymallera} in this case also; their theorem only
requires $\pibar(0+)>0$.

The remaining equivalences in case~(ii) follow exactly as in case~(i).
This completes the proof of Theorem \ref{PRS}.
\end{longlist}
\end{pf*}

The {\it domain of attraction of the normal distribution}, as
$t\downarrow0$,
appears in the next result, which is a corollary to Theorem \ref{PRS}.
We say $X\in D(N)$ at 0 if there are functions $a_t\in\myRR$,
$b_t>0$, such that $(X_t-a_t)/b_t\stackrel{\mathrm
{D}}{\longrightarrow}N(0,1)$ (a standard normal
random variable ) as $t\downarrow0$. If $a_t$ may be taken as 0, we write
$X\in D_0(N)$ (no centering required).
The following condition characterises the domain of attraction of the
normal at 0 (Doney and Maller \cite{doneymallera}, Theorem 2.5):
%
%e5.46 #&#
\begin{equation}
\label{cors2} \lim_{x\downarrow0} \frac{U(x)}{x^2\pibar(x)}=\infty;
\end{equation}
in fact, $D(N)$ (at 0) equals $D_0(N)$ (at 0) (Maller and Mason \cite
{MM10}, Theorem 2.4).
A characterisation for $D_0(N)$ at 0 (equivalent to \eqref{cors2}) is
%
%e5.47 #&#
\begin{equation}
\label{cors2a} \lim_{x\downarrow0} \frac{U(x)}{x\llvert A(x)\rrvert
+x^2\pibar(x)}=\infty.
\end{equation}

The following conditions are also equivalent to $X_t/b_t\stackrel
{\mathrm{D}}{\longrightarrow}N(0,1)$
(Kallenberg \cite{Kall}, Theorem~15.14): for all $x>0$,
%
%e5.48 #&#
\begin{equation}
\label{normcon} \lim_{t\downarrow0} t\pibar(xb_t)=0,\qquad \lim
_{t\downarrow0} \frac{tA(xb_t)}{b_t}=0,\qquad \lim_{t\downarrow0}
\frac{tU(xb_t)}{b_t^2}=1.
\end{equation}
Obvious modifications of these characterise the convergence
$X_t/b_t\stackrel{\mathrm{D}}{\longrightarrow}N(0,1)$ through a
subsequence $t_k\downarrow0$.

%
%co3 #&#
\begin{corollary}[(Corollary to Theorem \ref{PRS})]\label{corSPRS}
Assume $\pibar^+(0+)=\infty$. Then the following are equivalent:
%
%e5.49 #&#
\begin{eqnarray}
\label{cors1}
&&\mbox{there is a nonstochastic function } c_t>0
\mbox{ such that } \frac{V_t}{c_t} \buildrel P \over\to 1,\qquad\mbox
{as }t\downarrow0;
\\
%e5.50 #&#
\label{cors3} &&\frac{V_{t}}{\sup_{0<s\le t}\llvert {\Delta X}_s\rrvert
^2} \buildrel P \over\to
\infty,\qquad\mbox{as }t
\downarrow0;
\\
%
%%\label{cors4}
&& X \mbox{ is in the domain of attraction of the normal distribution, as }
t\downarrow0.\nonumber
\end{eqnarray}
\end{corollary}

\begin{pf} %{Proof of Corollary \ref{corSPRS}}
$V_t$ is a subordinator with drift $\rmd_V = \sigma^2$ and L\'evy
measure $\Pi_V$, where $\pibar_V(x)=\pibar(\sqrt{x})\bfeins_{\{
x>0\}}$.
Let the triplet of $V_t$ be $(\gamma_V,0,\Pi_V(\cdot))$. Then $\rmd
_V = \gamma_V +\int_0^1 y \Pi_V(\rmd y)$.
Thus, in obvious notation\vspace*{-2pt}
\begin{eqnarray*}
A_V(x)&=&\gamma_V+\pibar_V(1)-\int
_x^1 \pibar_V(y)\,\rmd y
= \rmd_V+\int_0^x
\pibar_V(y)\,\rmd y
\\[-3pt]
&=& \sigma^2+2\int_0^{\sqrt{x}} y
\pibar(y)\,\rmd y
= U(\sqrt{x}),\qquad x>0.
\end{eqnarray*}
Hence,\vspace*{-2pt}
\[
\frac{A_V(x)}{x\pibar_V(x)}= \frac{U(\sqrt{x})}{(\sqrt{x})^2\pibar
(\sqrt{x})}
\]
tends to $\infty$ iff \eqref{cors2} holds. By Theorem \ref{PRS}
these are equivalent to \eqref{cors1} and \eqref{cors3}, and \eqref
{cors2} characterises the domain of attraction of the normal, as noted.
\end{pf}

%
%re4 #&#
\begin{remark}\label{rem5}
(i) Another interesting kind of ``self-normalisation'' of a L\'evy process
is to divide $X_t$ by $\sqrt{V_t}$, possibly after removal of one or
the other kind of maximal jump. See, for example, Maller and Mason
\cite{MM08,MM13}. Our methods can be used to extend these
results in a variety of directions, but we omit further details here.

(ii) Relative stability of $X$ is directly related to the
stability of the ``one-sided'' and ``two-sided'' passage times over
power law boundaries defined by
\[
{\overline T}_b(r):= \inf\bigl\{t\ge0: X_t>rt^b
\bigr\},\qquad r\ge0,
\]
and
\[
T^*_b(r):= \inf\bigl\{t\ge0: \llvert X_t\rrvert
>rt^b\bigr\},\qquad r\ge0,
\]
when\footnote{Griffin and Maller \cite{GM2} show that relative stability
of ${\overline T}_b(r)$ or $T^*_b(r)$
cannot obtain when $b\ge1$.} $0\le b<1$. Griffin and Maller \cite
{GM2} show that, then,
${\overline T}_b(r)$ is relatively stable as $r\downarrow0$, in the
sense that
${\overline T}_b(r)/C(r)\buildrel P \over\to1$ as $r\downarrow0$
for a nonstochastic
function $C(r)>0$, iff $X\in \mathrm{PRS}$,
while
$T^*_b(r)$ is relatively stable as $r\downarrow0$, in the sense that
$T^*_b(r)/C(r)\buildrel P \over\to1$ as $r\downarrow0$ for a
nonstochastic function
$C(r)>0$, iff $X\in RS$.
Further connections made in Griffin and Maller \cite{GM2} are that
$X\in \mathrm{PRS}$ iff
${\overline X_t}:= \sup_{0<s\le t}X_s$ is relatively stable,
while
$X\in RS$ iff $X_t^*:= \sup_{0<s\le t}\llvert X_s\rrvert $ is
relatively stable.
Auxiliary results are (i)
there is a nonstochastic function $b_t^*>0$ and constants
$0<c_1<c_2<\infty$ such that $\lim_{t\downarrow0}P(c_1<\llvert
X_t\rrvert /b_t^*<c_2)=1$
iff $X\in RS$, and (ii) there is a nonstochastic function $b_t^\dagger
>0$ such that
each sequence $t_k\downarrow0$ contains a subsequence
$t_{k'}\downarrow0$ with
$\llvert X_{t_{k'}}\rrvert /b^\dagger_{t_{k'}}\buildrel P \over\to
c'$, where
$0<\llvert c'\rrvert <\infty$, iff
$X\in RS$.
See also Griffin and Maller \cite{GM1}.
\end{remark}

%s5.3 #&#
\subsection{Relative stability, attraction to normality and
dominance}\label{sub43}
The next theorems look at two-sided results, concerning stability and
dominance of $\llvert X\rrvert $.
Now the domain of attraction of the normal enters as an alternative to
relative stability.

%
%th8 #&#
\begin{theorem}\label{PRS3}
Assume $\pibar(0+)=\infty$. Then the following are equivalent:
%
%e5.51 #&#
\begin{eqnarray}
\label{mrs4} \frac{\llvert X_t\rrvert }{\llvert \widetilde{\Delta
X}_t^{(1)}\rrvert }&\buildrel P \over\to
&\infty,\qquad\mbox{as }t
\downarrow0;
\\
%
%%% next claim is incorrect?
%\be\label{mrs1}
%\frac{\left\vert X_t\right\vert}{(\Delta X^\pm)_t^{(1)}}\topr\infty, {
%\rm as} t\dto0;
%\ee
%
%e5.52 #&#
\label{mrs2} \lim_{x\downarrow0} \frac{x\llvert A(x)\rrvert
+U(x)}{x^2\pibar(x)} &=&\infty;
\\
%
%
%e5.53 #&#
\label{mrs3} \lim_{x \downarrow0} {U(x) \over x \llvert A(x)\rrvert
+x^2\pibar(x)} &=&+\infty,
\quad\mbox{or}\quad\lim_{x \downarrow0} {\llvert A(x)\rrvert \over
x\pibar(x)} =+
\infty;
\\
%
%e5.54 #&#
\label{mrs5} X&\in& D_0(N)\cup RS\qquad\mbox{at } 0.
\end{eqnarray}
\end{theorem}

\begin{pf} %{Proof of Theorem \ref{PRS3}} %\begin{proof}[of
%Theorem {\rm\ref{PRS3}}]
Assume $\pibar(0+)=\infty$. \eqref{mrs4} $\Longrightarrow$ \eqref
{mrs2}: Assume \eqref{mrs4}. This implies
\[
\frac{\llvert ^{(1)}\widetilde X_t\rrvert }{\llvert \widetilde{\Delta
X}_t^{(1)}\rrvert }\buildrel P \over\to\infty,\qquad\mbox{as
}t\downarrow0,
\]
so we have
%
%e5.55 #&#
\begin{equation}
\label{mrs11} \lim_{t\downarrow0}P \bigl(\bigl\llvert
{}^{(1)}\widetilde X_t\bigr\rrvert\le a\bigl\llvert
\widetilde{\Delta X}_t^{(1)}\bigr\rrvert\bigr)=0\qquad\mbox{for some
} a>0.
\end{equation}
Without loss of generality take $a\le1$.

We again abbreviate $\overline{\Pi}^{\leftarrow}(v)$ to $y_v$ throughout.
Then by \eqref{cor1a}, we can write
%
%e5.56 #&#
\begin{equation}
\label{mrs12} P \bigl(\bigl\llvert{}^{(1)}\widetilde X_t
\bigr\rrvert\le a\bigl\llvert\widetilde{\Delta X}_t^{(1)}
\bigr\rrvert\bigr) =\int_0^\infty P \bigl(\bigl
\llvert\widetilde{ X}_t^{v} + \widetilde{G}_t^{v}
\bigr\rrvert\le ay_v \bigr)P(\EEEE\in t\,\rmd v).
\end{equation}
By \eqref{ega}, we have
%
%e5.57 #&#
\begin{eqnarray}
\label{egam} \bigl\llvert E\widetilde{G}_t^{v}\bigr
\rrvert\le ty_v\pibar(y_v-)\le tU(y_v)/y_v,
\end{eqnarray}
and \eqref{vgt} remains true.
%\be\label{vgtm}
%\operatorname{Var}(\wt{G}_t^{v})\le ty_v^2 \pibar(y_v-).
%\ee
Also, as in \eqref{exa},
%
%e5.58 #&#
\begin{eqnarray}
\label{exam} \bigl\llvert E\widetilde{ X}_t^{v}\bigr
\rrvert&=& t\bigl\llvert A(y_v)-y_v \bigl(
\pibar^+(y_v-) -\pibar^-(y_v-) \bigr)\bigr\rrvert\nonumber
\\
&\le& t \bigl(\bigl\llvert A(y_v)\bigr\rrvert+y_v
\pibar(y_v-) \bigr)
\\
&
\le& t\bigl(\bigl\llvert A(y_v)\bigr\rrvert+U(y_v)/y_v
\bigr).\nonumber
\end{eqnarray}
Apply \eqref{egam} and \eqref{exam} to obtain from \eqref{mrs12}
\begin{eqnarray}
\label{mrs13} &&P \bigl(\bigl\llvert{}^{(1)}\widetilde
X_t\bigr\rrvert\le a\bigl\llvert\widetilde{\Delta
X}_t^{(1)}\bigr\rrvert\bigr)
\nonumber
\\
&&\quad\ge\int_0^\infty P \bigl(\bigl\llvert
\widetilde{ X}_t^{v}-E\widetilde{ X}_t^{v}
+ \widetilde{G}_t^{v}-E\widetilde{ G}_t^{v}
\bigr\rrvert\le ay_v-\bigl\llvert E\widetilde{ X}_t^{v}
\bigr\rrvert- \bigl\llvert E\widetilde{ G}_t^{v}\bigr
\rrvert\bigr)P(\EEEE\in t\,\rmd v)
\\
&&\quad\ge\int_0^\infty P \bigl(\bigl\llvert
\widetilde{ X}_t^{v}-E\widetilde{ X}_t^{v}
+ \widetilde{G}_t^{v}-E\widetilde{ G}_t^{v}
\bigr\rrvert\le ay_v-2t\bigl(\bigl\llvert A(y_v)\bigr
\rrvert+U(y_v)/y_v\bigr) \bigr)P(\EEEE\in t\,\rmd v).
\nonumber
\end{eqnarray}

For $t>0$, define
%
%e5.59 #&#
\begin{equation}
\label{mbdef} b_t:= \sup\biggl\{x>0: \frac{x\llvert A(x)\rrvert
+U(x)}{x^2} >
\frac
{a^2}{56t} \biggr\},
\end{equation}
with $b_0:=0$. Since $\pibar(0+)=\infty$, we have
$\lim_{x\downarrow0}(x\llvert A(x)\rrvert +U(x))/x^2=\infty$.
In addition, $\lim_{x\to\infty}(x\llvert A(x)\rrvert +U(x))/x^2=0$.
Then $0<b_t<\infty$, $b_t$ is strictly increasing, $b(t)\downarrow0$ as
$t\downarrow0$,
and
%
%e5.60 #&#
\begin{equation}
\label{mrs14} \frac{t(b_t\llvert A(b_t)\rrvert +U(b_t))}{b_t^2}= \frac
{a^2}{56},\qquad t>0.
\end{equation}
Now keep $v<\pibar(b_t)$. Then $y_v\ge b_t$, and so
\[
t\bigl(\bigl\llvert A(y_v)\bigr\rrvert+U(y_v)/y_v
\bigr)\le\frac{a^2y_v}{56}\le\frac{ay_v}{4},
\]
by definition of $b_t$.
Thus, by Chebyshev's inequality and \eqref{mrs13}
\begin{eqnarray*}
%%\label{mrs14a}
&& P \bigl(\bigl\llvert{}^{(1)}\widetilde X_t\bigr
\rrvert\le a\bigl\llvert\widetilde{\Delta X}_t^{(1)}\bigr
\rrvert\bigr)
\\
&&\quad  \ge \int_0^{\pibar(b_t)} P \bigl(\bigl
\llvert\widetilde{ X}_t^{v}-E\widetilde{
X}_t^{v} + \widetilde{G}_t^{v}-E
\widetilde{ G}_t^{v}\bigr\rrvert\le ay_v/2
\bigr)P(\EEEE\in t\,\rmd v)
\\
&&\quad \ge \int_0^{\pibar(b_t)} \biggl(1-
\frac{4 (\operatorname{Var}(\widetilde{X}_t^{v})+\operatorname
{Var}(\widetilde
{G}_t^{v}) )}{
a^2y_v^2} \biggr)P(\EEEE\in t\,\rmd v).
\end{eqnarray*}
Also, as in \eqref{prs15},
\begin{eqnarray*}
%%\label{mprs15}
\operatorname{Var}\bigl(\widetilde{X}_t^{v}\bigr)+\operatorname
{Var}\bigl(\widetilde{G}_t^{v}\bigr) %&\le&
%t\left(V(y_v)+y_v^2\pibar(y_v-)\right)
% \nonumber\\
% &\le&
% 2tU(y_v)
% \nonumber\\
&
\le& a^2y_v^2/8,
\end{eqnarray*}
giving
%
%e5.61 #&#
\begin{equation}
\label{sprs17} P \bigl(\bigl\llvert{}^{(1)}\widetilde X_t
\bigr\rrvert\le a\bigl\llvert\widetilde{\Delta X}_t^{(1)}
\bigr\rrvert\bigr) \ge t\int_0^{\pibar(b_t)}e^{-tv}\,
\rmd v/2 = \bigl(1- e^{-t\pibar(b_t)} \bigr)/2.
\end{equation}
Since the left-hand side tends to 0 as $t\downarrow0$ by \eqref
{mrs11} we
see that
%
%e5.62 #&#
\begin{equation}
\label{mrs16} t\pibar(b_t) =\frac{a^2b_t^2\pibar(b_t)}{56(b_t\llvert
A(b_t)\rrvert +U(b_t))}\to0,\qquad\mbox{as }t
\downarrow0.
\end{equation}
%
% thus
% \ben
%\frac{b_t\left\vert A(b_t)\right\vert+U(b_t)}
%{b_t^2\pibar(b_t)}\to\infty, {\rm as} t\dto0.
% \een

We need to replace $b_t$ by a continuous variable $x\downarrow0$ in this.
% It will be convenient to define the function
% \be\label{gdef}
% g(x):=
% \frac{x\left\vert A(x)\right\vert+U(x)}{x^2}, x>0.
% \ee
By \eqref{mrs14}, for $\lambda>1$ and $t>0$
\begin{eqnarray}
\label{bin} \frac{ b_{t\lambda}^2}{b_t^2} &=& \frac{56t\lambda
(b_{t\lambda}\llvert A(b_{t\lambda})\rrvert +U(b_{t\lambda}))}{a^2b_t^2}
\nonumber
\\
&=& \frac{56t\lambda(b_{t\lambda}\llvert A(b_t)\rrvert +U(b_t))}
{a^2b_t^2}+ \frac{56t\lambda b_{t\lambda}(\llvert A(b_{t\lambda
})\rrvert -\llvert A(b_t)\rrvert )} {a^2b_t^2}
\\
&&{} + \frac{56t\lambda
(U(b_{t\lambda})-U(b_t))} {a^2b_t^2}
\nonumber
\\
&\le& %%\frac{t\lambda(b_t\left\vert A(b_t)\right\vert+U(b_t))} {b_t^2}
\lambda+ \frac{56t\lambda(b_{t\lambda}-b_t)\llvert A(b_t)\rrvert }
{a^2b_t^2}+ \frac{56t\lambda b_{t\lambda}(b_{t\lambda}-b_t)\pibar(b_t)}
{a^2b_t^2}+
\frac{56t\lambda(b_{t\lambda}^2-b_t^2)\pibar(b_t)} {a^2b_t^2}.
\nonumber
\end{eqnarray}
Observe that $56t\lambda(b_{\lambda t}-b_t)\llvert A(b_t)\rrvert
/a^2b_t^2\le
\lambda(b_{\lambda t}-b_t)/b_t$. Since $t\pibar(b_t)=o(1)$, \eqref
{bin} implies
\[
\frac{ b_{t\lambda}^2}{b_t^2} \le\lambda+ \lambda\biggl( \frac{
b_{t\lambda}}{b_t}-1 \biggr) +
o \biggl( \frac{ b_{t\lambda}^2}{b_t^2} \biggr) \le\lambda+ \lambda
\frac{ b_{t\lambda}}{b_t}+ o
\biggl( \frac{
b_{t\lambda}^2}{b_t^2} \biggr).
\]
From this, we deduce that $\limsup_{t\downarrow0}b_{t\lambda
}/b_t<\infty$.
%This is ``dominated variation" of $b_t$ at 0, and it is well known
%then and easily established that there are constants $c>0$, $\beta>0$,
%such that
%$b_{t\lambda}/b_t\le c\lambda^\beta$ for all $\lambda>1$ and $t$
%small enough.

Now return to \eqref{mrs16} and take $x>0$. Choose $t=t(x)$ such
that $b_t\le x\le b_{\lambda t}$, $\lambda>1$.
It is shown in Klass and Wittmann \cite{KW}
that the function $x\llvert A(x)\rrvert +U(x)$ is
nondecreasing\footnote{Klass and
Wittmann prove this for versions of $A$ and $U$ defined for distribution
functions. But their proof is easily modified to apply to the present
$A$ and $U$.} in $x>0$. Thus,
\begin{eqnarray*}
\frac{x\llvert A(x)\rrvert +U(x)}{x^2\pibar(x)} &\ge& \frac
{b_t\llvert A(b_t)\rrvert +U(b_t)}{b_t^2\pibar(b_t)}\times\frac
{b_t^2}{b_{\lambda t}^2}.
\end{eqnarray*}
The first factor on the right tends to $\infty$ as $t\downarrow0$ by
\eqref
{mrs16}, and
$\liminf_{t\downarrow0}b_t/b_{t\lambda}>0$,
so we get \eqref{mrs2}.

\eqref{mrs2} $\iff$ \eqref{mrs3} is proved in Lemma 4 of Doney and
Maller \cite{doneymallerb}.

\eqref{mrs3} $\Longrightarrow$ \eqref{mrs5}:
Assume~\eqref{mrs3}. If $\sigma^2>0$ then by Lemma \ref{sigma2},
$X\in D_0(N)$ hence $X\in D_0(N)\cup RS$.
So suppose $\sigma^2=0$. Then the left-hand side of \eqref{mrs3} is
equivalent to
$X\in D_0(N)$ at 0 by \eqref{cors2a}, and
the right-hand side of \eqref{mrs3} is equivalent to $X_t\in RS$ at 0
by \eqref{prs2a}. Thus again, $X\in D_0(N)\cup RS$.

\eqref{mrs5} $\Longrightarrow$ \eqref{mrs4}:
Finally, if $X\in D_0(N)\cup RS$ then $X_t/b_t\stackrel
{D}{\longrightarrow}N(0,1)$ for some $b_t>0$
with $\widetilde{\Delta X}_t^{(1)}=o_P(b_t)$ or $X_t/c_t\buildrel P
\over\to \pm1$ for some
$c_t>0$ with $\widetilde{\Delta X}_t^{(1)}=o_P(c_t)$, and in either case
\eqref{mrs4} holds. This completes Theorem \ref{PRS3}.
\end{pf}
\subsection{Subsequential relative stability and dominance}\label{sub44}
We say that $X$ is \textit{subsequentially relatively stable} (SRS) at 0
if there are nonstochastic sequences $t_k\downarrow0$ and $b_k>0$ such that
%
%e5.63 #&#
\begin{equation}
\label{srs0} \frac{X_{t_k}}{b_k} \buildrel P \over\to\pm1,\qquad\mbox
{as }k\to
\infty.
\end{equation}
Define positive and negative subsequential relative stability (PSRS and NSRS)
in the obvious ways.

%
%th9 #&#
\begin{theorem}\label{SPRS}
Assume $\pibar^+(0+)=\infty$.
Then the following are equivalent:
there is a nonstochastic sequence $t_k\downarrow0$ such that
%
%e5.64 #&#
\begin{equation}
\label{sprs4} \frac{X_{t_k}}{\llvert \widetilde{\Delta
X}_{t_k}^{(1)}\rrvert }\buildrel P \over
\to \infty,\qquad\mbox{as }k\to
\infty;
\end{equation}
there is a nonstochastic sequence $t_k\downarrow0$ such that
%
%e5.65 #&#
\begin{eqnarray}
\label{sprs4a}
\frac{X_{t_k}}{(\Delta X^+)_{t_k}^{(1)}}&\buildrel P
\over\to
&\infty,\qquad\mbox{as }k\to
\infty;
\\
%
%e5.66 #&#
\label{sprs1} X&\in&\operatorname{PSRS} \mbox{at } 0;
\\
%
%e5.67 #&#
\label{sprs3} \limsup_{x\downarrow0} \frac{A(x)}{\sqrt{U(x)\pibar
(x)}}&=&\infty;
\\
%
%e5.68 #&#
\label{sprs3a} \limsup_{x\downarrow0} \frac{xA(x)}{U(x)}&=&\infty.
\end{eqnarray}
\end{theorem}

\begin{pf}%{Proof of Theorem \ref{SPRS}} %\begin{proof}[of
%Theorem {\rm\ref{SPRS}}]
Assume $\pibar^+(0+)=\infty$. %% and $\pibar^-(0+)>0$.
Each of \eqref{sprs4}--\eqref{sprs3a} implies $\sigma^2=0$; by Lemma
\ref{sigma2} in the case of \eqref{sprs4} and \eqref{sprs4a},
by Lemma \ref{sig1} in the case of \eqref{sprs3}, and by \eqref{xnu}
and $U(x)\ge\sigma^2$, in the case of \eqref{sprs3a}.
So we assume throughout that $\sigma^2=0$.

\eqref{sprs4} $\iff$ \eqref{sprs4a}: clearly, \eqref{sprs4}
implies \eqref{sprs4a}. Conversely, assume \eqref{sprs4a}. From
\eqref{Xdels}, we have that
$X_{t_k}/(\Delta X^-)_{t_k}^{(1)}\buildrel P \over\to\infty$, as
$k\to\infty$,
when $\lim_{k\to\infty}P(X_{t_k}> 0)=1$.
Together with \eqref{sprs4a} and
$\llvert \widetilde{\Delta X}_t^{(1)}\rrvert =\max((\Delta
X^+)_t^{(1)}, (\Delta
X^-)_t^{(1)})$,
this implies \eqref{sprs4}.

%%We first prove the equivalence of \eqref{sprs1} and \eqref{sprs3}.
\eqref{sprs3} $\iff$ \eqref{sprs3a}: Assume \eqref{sprs3}, so
there is a nonstochastic sequence $x_k\downarrow0$ such that
\[
%%\label{41}
\frac{A(x_k)}{\sqrt{U(x_k)\pibar(x_k)}}\to\infty,\qquad\mbox{as }k\to
\infty.
\]
Define
\[
%%\label{tdef}
t_k= \frac{1}{A(x_k)} \sqrt{\frac{U(x_k)}{\pibar(x_k)}}.
\]
Then
\[
%%\label{43x}
t_k\pibar(x_k) =\frac{\sqrt{U(x_k)\pibar(x_k)}}{A(x_k)} \to0
\]
and so, since $\pibar(0+)>0$,
$t_k\to0$. Also
\[
%%\label{43y}
\frac{U(x_k)}{t_kA^2(x_k)}=\frac{1}{A(x_k)}\sqrt{\pibar
(x_k)U(x_k)}\to0.
\]
Let $b_k=t_kA(x_k)$, then
\[
\frac{b_k}{x_k} =\frac{t_kA(x_k)}{x_k}= \sqrt{\frac{U(x_k)}{x_k^2\pibar(x_k)}}
\ge1.
\]
Now since $b_k \ge x_k$ we have
\begin{eqnarray*}
\frac{t_kU(b_k)}{b_k^2} &=& \frac{U(x_k)}{t_kA^2(x_k)} + \frac{2t_k\int
_{x_k}^{b_k} y\pibar(y)\,\rmd y}{b_k^2}
\nonumber
\\
&\le& o(1) +O \bigl(t_k\pibar(x_k) \bigr) =o(1).
\end{eqnarray*}
This implies $t_kU(xb_k)/b_k^2=o(1)$ for all $x\in(0,1]$, hence
%
%
%e5.69 #&#
\begin{equation}
\label{47} \lim_{k\to\infty} t_k
\pibar(xb_k)= 0\qquad\mbox{for all } x\in(0,1],
\end{equation}
because $U(x)\ge x^2\pibar(x)$.
But then since $\pibar$ is nonincreasing, \eqref{47} holds for all $x>0$.
Thus, also, for $x>1$,
%
%
%e5.70 #&#
\begin{equation}
\label{48} \frac{t_kU(xb_k)}{b_k^2} = \frac{t_kU(b_k)}{b_k^2} +O
\bigl(t_k\pibar(b_k) \bigr)=o(1).
\end{equation}
Again since $b_k \ge x_k$, we can write
%
%
%e5.71 #&#
\begin{equation}
\label{49} \frac{t_kA(b_k)}{b_k}=1+\frac{t_k\int_{x_k}^{b_k}
(\pibar^+(y)-\pibar^-(y))\, \rmd y}{b_k} =1+O
\bigl(t_k\pibar(x_k) \bigr) =1+o(1).
\end{equation}
\eqref{48} and \eqref{49}, hence \eqref{sprs3}, imply \eqref{sprs3a}.
Conversely, \eqref{sprs3a} implies \eqref{sprs3} because
$U(x)\ge x^2\pibar(x)$.

\eqref{sprs3} $\iff$ \eqref{sprs1}: \eqref{sprs3} implies \eqref
{47}--\eqref{49}, as just shown, and
these together imply \eqref{srs0} (with a ``$+$'' sign) by the
subsequential version of \eqref{degcon}. Thus, \eqref{sprs1} holds.
Conversely, assuming \eqref{sprs1}, we get \eqref{47}--\eqref{49} by
the subsequential version of
\eqref{degcon}.
But then \eqref{sprs3} holds because
\[
\frac{A(b_k)}{\sqrt{U(b_k)\pibar(b_k)}}= \frac{t_kA(b_k)}{b_k} \sqrt{
\biggl(
\frac{b_k^2}{t_kU(b_k)} \biggr) \biggl(\frac{1}{t_k\pibar(b_k)} \biggr)
} \to\infty.
\]
So we have proved the equivalence of \eqref{sprs1}--\eqref{sprs3a}.

%%It remains to deal with \eqref{sprs4} and \eqref{sprs4a}.

\eqref{sprs4} $\Longrightarrow$ \eqref{sprs3}: Assume \eqref{sprs4}.

\begin{longlist}
\item[\textit{Case} (i).] Suppose $\pibar^-(0+)>0$.
Then, using Theorem \ref{th33s}, we have $\lim_{k\to\infty
}P(X_{t_k}> 0)=1$, $\sigma^2=0$,
and \eqref{36as}. Since $\pibar^-(0+)>0$ and $U(x)\ge x^2\pibar
^-(x)$, \eqref{36as} implies $\limsup_{x\downarrow0}A(x)/ x=\infty$.
\eqref{sprs4} also implies
\[
\frac{^{(1)}\widetilde{X}_{t_k}}{\llvert \widetilde{\Delta
X}_{t_k}^{(1)}\rrvert }\buildrel P \over\to\infty,\qquad\mbox{as
}k\to\infty,
\]
so we have
\[
%%\label{sprs11}
\lim_{k\to\infty}P \bigl({}^{(1)}\widetilde
X_{t_k}\le a\bigl\llvert\widetilde{\Delta X}_{t_k}^{(1)}
\bigr\rrvert\bigr)=0\qquad\mbox{for some } a\in(0,1).
\]
%
%%(In fact, this holds for all $a>0$. But it will be enough to assume
%\eqref{sprs11}.)
Define $b_k$ similarly as in \eqref{mbdef}:
%
%e5.72 #&#
\begin{equation}
\label{bkdef} b_k:= \sup\biggl\{x>0: \frac{x\llvert A(x)\rrvert
+U(x)}{x^2} >
\frac
{a^2}{56t_k} \biggr\}.
\end{equation}
%
% Since $\limsup_{x\dto0}(x\left\vert A(x)\right\vert+U(x))/x^2=
%\infty$ and
% $\lim_{x\to\infty}(x\left\vert A(x)\left\vert+U(x))/x^2=0$,
% the $b_k$ are finite and positive, and satisfy
%\be\label{sprs5a}
%\frac{t_k(b_k\left\vert A(b_k)\left\vert+U(b_k))}{b_k^2} =
%\frac{a^2}{56}.
%\ee
%Again abbreviate $\pibarinv(v)$ to $y_v$ throughout.
%As in \eqref{ega}--\eqref{exa},
%\[
%\left\vert E\wt{G}_t^{v}\left\vert\le tU(y_v)/y_v\quad{\rm and}\quad
%\left\vert E\wt{X}_t^{v}\left\vert\le
%t\left(\left\vert A(y_v)\left\vert+U(y_v)/y_v\right).
%\]
%Keep $v\le\pibar(b_k)$, so $y_v\ge b_k$.
%By the definition of $b_k$, $y_v\ge b_k$ implies \ben
%t_k\left(\right\vert A(y_v)\right\vert+t_kU(y_v)/y_v\right)\le
%a^2y_v/56\le ay_v/4.
%\een
% Hence, as in \eqref{prs14a},
%\be\label{sprs13}
%P\left({}^{(1)}\wt X_{t_k}\le a\right\vert\wt{\Delta X}_{t_k}^{(1)}
%\right\vert\right)
% \ge t_k\int_0^{\pibar(b_k)}
% \left(1- \frac{4 \left(\operatorname{Var}(\wt{X}_{t_k}^{v})+
%\operatorname{Var}(
%\wt{G}_{t_k}^{v})\right)}
% {a^2y_v^2}\right)e^{-t_kv}\rmd v.
%\ee
% As in \eqref{prs15}
% \bea\label{prs155}
% \frac{\operatorname{Var}(\wt{X}_{t_k}^{v})+\operatorname{Var}(
%\wt{G}_{t_k}^{v})}{a^2y_v^2}
%&\le&
%\frac{ 2t_kU(y_v)}{a^2y_v^2}
% \nonumber\\
% &\le&
%\frac{ 2t_kU(b_k)}{a^2b_k^2}
% \quad({\rm since} y_v\ge b_k {\rm and} x^{-2}U(x) {\rm
%decreases})
% \nonumber\\
% &\le&
%\frac{1}{23}.
% \eea
%
Then by the same calculation as in \eqref{mbdef}--\eqref{sprs17}, we
find, for large $k$,
\[
P \bigl({}^{(1)}\widetilde X_{t_k}\le a\bigl\llvert
\widetilde{\Delta X}_{t_k}^{(1)}\bigr\rrvert\bigr) \ge
t_k\int_0^{\pibar(b_k)}e^{-t_kv}\,
\rmd v\big/2 = \bigl(1- e^{-t_k\pibar(b_k)} \bigr)/2.
\]
From this, we conclude that $t_k\pibar(b_k)\to0$.
Take a subsequence $k'\to\infty$ if necessary so that

%e5.73 #&#
\begin{equation}
\label{AB} \frac{t_{k'}A(b_{k'})}{b_{k'}} \to A \quad\mbox{and}\quad
\frac{t_{k'}U(b_{k'})}{b_{k'}^2}\to B,
\end{equation}
where $B\ge0$ and $\llvert A\rrvert +B=a^2/56$.

Now $A\le0$ is not possible in \eqref{AB}. To see this, take a
further subsequence of $k'$ if necessary so that, for some functions
$\lambar^\pm(x)$ and $B(x)$,
\[
\lim_{k'\to\infty}t_{k'}\pibar^\pm(xb_{k'})=
\lambar^\pm(x) \quad\mbox{and}\quad\lim_{k'\to\infty}
\frac{t_{k'}U(xb_{k'})}{b_{k'}^2}= B(x)
\]
at continuity points $x>0$ of these functions. Let $\Lambda$ be the measure
having positive and negative tails $\lambar^\pm$.
Then $\lambar(x):= \lambar^+(x)+\lambar^-(x)=0$ for all $x\ge1$.
Fatou's lemma gives
\[
\infty>B=\lim_{k'\to\infty}\frac{t_{k'}U(b_{k'})}{b_{k'}^2}=2\lim
_{k'\to\infty}\int_0^1yt_{k'}
\pibar(yb_{k'})\,\rmd y \ge2\int_0^1y
\lambar(y)\,\rmd y,
\]
and shows that the integral on the right is finite.
This means that $\Lambda$ is a L\'evy measure on $\myRR$ and by
Kallenberg (\cite{Kall}, Theorem 15.14), as $k'\to\infty$ we have
$(X_{t_{k'}}-t_{k'}\nu(b_{k'}))/b_{k'}\stackrel{\mathrm
{D}}{\longrightarrow}Y'$,
an infinitely divisible r.v. with canonical measure
$\Lambda$. Since $\lambar(x)=0$ for all $x\ge1$,
$Y'$ has finite variance. Further, since
$t_k\pibar(b_k)\to0$ we have
$\lim_{k'\to\infty}t_{k'}\nu(b_{k'})/b_{k'}=A$
(recall \eqref{Apart}).
The L\'evy--It\^o decomposition can equivalently be written as
%
%
%e5.74 #&#
\begin{equation}
\label{decomp1} X_t=t\nu(b)+\sigma Z_t+X_t^{(S,b)}+X_t^{(B,b)},\qquad
t\ge0,
\end{equation}
where $b>0$, $X_t^{(S,b)}$ is the compensated small jump component of
$X$, that is, having jumps less than or equal to $b$ in modulus, and
$X_t^{(B,b)}$ is the sum of jumps larger in modulus than $b$; see, for
example, Doney and Maller (\cite{doneymallera}, Lemma 6.1).
Choose $b=b_k$ in \eqref{decomp1}, and notice that the sum of jumps
larger in modulus than $b_k$ is $o(b_k)$ as $k \to\infty$ because
$t_k\pibar(b_k)\to0$. Also, $\sigma^2=0$. So we deduce
%
%e5.75 #&#
\begin{equation}
\label{stu} \frac{X_{t_{k'}}^{(S,b_{k'})}-t_{k'}\nu(b_{k'})}{b_{k'}}
=\frac{X_{t_{k'}}-t_{k'}\nu(b_{k'})}{b_{k'}}+o_P(1)
\stackrel{\mathrm{D}} {\longrightarrow}Y'.
\end{equation}
From the inequality,
\[
\frac{E(X_{t_{k'}}^{(S,b_{k'})})^2}{b^2_{k'}} \le\frac
{t_{k'}U(b_{k'})}{b_{k'}^2} \le\frac{a^2}{56}
\]
we see that $(X_{t_{k'}}^{(S,b_{k'})}/b_{k'})$
is uniformly integrable. Thus, we deduce from \eqref{stu} that
\[
\frac{E(X_{t_{k'}}^{(S,b_{k'})})}{b_{k'}}\to EY'+A.
\]
The expectation on the left equals 0, so this implies $EY'=-A$. Now
argue that
\[
\lim_{k'\to\infty}P( X_{t_{k'}}\le0) = \lim
_{k'\to\infty}P \biggl( \frac{X_{t_{k'}}-t_{k'}\nu(b_{k'})}{b_{k'}}\le
-\frac{t_{k'}\nu(b_{k'})}{b_{k'}}
\biggr)= P\bigl(Y'\le-A\bigr).
\]
But since $Y'+A$ has mean 0 and finite variance,
$P(Y'\le-A)=P(Y'+A\le0)>0$, in contradiction to \eqref{sprs4}.
Thus, $A\le0$ is not possible.

% $A=-\infty$ in \eqref{AB} also is not possible. If it occurred we
%would have
% $t_{k'}A(b_{k'})/b_{k'}\to-\infty$, hence, since $t_k\pibar(b_k)\to
%0$, also $t_{k'}\nu(b_{k'})/b_{k'}\to-\infty$.
%% Let $X_t^{(S,b)}$ for $b>0$ be the compensated small jump component
%of $X$, i.e., having jumps less than or equal to $b$ in modulus.
%Then by Chebyshev's inequality
%\bean
%P( X_{t_{k'}}\le0)&\ge&
%P\left(X_{t_{k'}}^{(S,b_{k'})} \le-t_{k'}\nu(b_{k'})\right) - t_{k'}
%\pibar(b_{k'})\cr
%&=&
%1-P\left(X_{t_{k'}}^{(S,b_{k'})} >t_{k'}\left\vert\nu(b_{k'})\right
%\vert\right) -o(1)\cr
%&\ge&
%1- \frac{U(b_{k'})}{t_{k'}\nu^2(b_{k'})}-o(1).
%\eean
%The last quantity tends to 1 as $k'\to\infty$ because, by
%\eqref{sprs5a},
%\ben
% \frac{t_{k'}U(b_{k'})}{b_{k'}^2}=\frac{a^2}{56}-
%\frac{t_{k'}\right\vert A(b_{k'})\right\vert}{b_{k'}}=\frac{a^2}{56}-
%\frac{t_{k'}\right\vert
%\nu(b_{k'})\right\vert}{b_{k'}}+o(1),
% \een
% giving
%\ben
% \frac{U(b_{k'})}{t_{k'}\nu^2(b_{k'})}=\left(\frac{a^2}{56}+o(1)\right)
% \frac{b_{k'}^2}{t_{k'}^2\nu^2(b_{k'})}-
% \frac{b_{k'}}{t_{k'}\right\vert\nu(b_{k'})\right\vert}\to0.
% \een
% Again we get a contradiction to \eqref{sprs4}.
We conclude that $A>0$ and $B<\infty$. It follows from \eqref{AB} that
\[
\frac{A(b_{k'})}{\sqrt{U(b_{k'})\pibar(b_{k'})}}\to\infty,
\]
which implies \eqref{sprs3}.
%% is implied by \eqref{sprs4} and $\pibar^-(0+)>0$.

\item[\textit{Case} (ii).] Still assuming \eqref{sprs4},
suppose $\pibar^-(0+)=0$. \eqref{sprs4} implies $P(X_{t_k}> 0)\to1$,
hence by Theorem \ref{th33s}, $X$ is a subordinator and $A(x)\ge0$
for all $x\ge0$. Then
\[
x^{-1}A(x)=x^{-1} \biggl(\rmd_X+\int
_0^x\pibar^+(y)\,\rmd y \biggr) \ge\int
_0^1\pibar^+(xy)\,\rmd y \to\infty,\qquad\mbox{as }x
\downarrow0,
\]
so we can define $b_k$ by \eqref{bkdef} and proceed as before to get
$t_k\pibar(b_k)\to0$, and hence \eqref{sprs3}.
%%Thus \eqref{sprs1} holds in this case too.

Conversely, in either cases (i) or~(ii), we know \eqref{sprs3}
$\Longrightarrow$ \eqref{sprs1}, and \eqref{sprs1} $\Longrightarrow
$ \eqref{sprs4}
follows easily from the subsequential version of \eqref{degcon}.\quad\qed
\end{longlist}\noqed
\end{pf}

%
% \eqref{sprs3a} implies \eqref{sprs3} because $U(x)\ge x^2\pibar(x)$.
% Conversely, \eqref{sprs1} implies \eqref{sprs3a} follows from the
%subsequential version of \eqref{degcon}.

%\eqref{sprs4} $\iff$ \eqref{sprs4a}: clearly \eqref{sprs4} implies
%\eqref{sprs4a}. Conversely, assume \eqref{sprs4a}. From \eqref{Xdels}
%we have that
%$X_{t_k}/(\Delta X^-)_{t_k}^{(1)}\topr\infty$, as $k\to\infty$,
%when $\lim_{k\to\infty}P(X_{t_k}> 0)=1$.
%Together with \eqref{sprs4a} and
%$\left\vert\wt{\Delta X}_t^{(1)}\right\vert=\max((\Delta
%X^+)_t^{(1)}, (\Delta
%X^-)_t^{(1)})$,
%this implies \eqref{sprs4}.
%%%\end{pf*}

The following corollary to Theorem \ref{SPRS} is also proved in
Theorem 4 of Maller \cite{PLMS09}.

%
%co4 #&#
\begin{corollary}\label{th4}
Assume $\pibar(0+)>0$.
The following are equivalent:

\begin{longlist}[(iii)]
\item[(i)] $X_t\in SRS$ at 0;

\item[(ii)] there are nonstochastic sequences $t_k \downarrow0$ and
$b_k>0$, such that, as $k \to\infty$,
%
%
%e5.76 #&#
\begin{equation}
\label{094a} \frac{\llvert X_{t_k }\rrvert }{b_k} \buildrel P \over
\to 1;
\end{equation}

\item[(iii)]
%
%e5.77 #&#
\begin{equation}
\label{095}
\sigma^2=0\quad\mbox{and}\quad \limsup_{x \downarrow0} \frac{
\vert A(x)\vert}{\sqrt{\pibar(x)U(x)}} =
\infty;
\end{equation}

\item[(iv)]
%
%e5.78 #&#
\begin{equation}
\label{096} \limsup_{x \downarrow0} \frac{x\llvert A(x)\rrvert
}{U(x)}=\infty.
\end{equation}
\end{longlist}
\end{corollary}

\begin{pf} %{Proof of Corollary \ref{th4}}
Assume $\pibar(0+)>0$.
First, %\begin{proof}[of Corollary {\rm\ref{th4}}]
$X_t\in SRS$ at 0 $\Longrightarrow$ \eqref{094a}
%%{\rm(i)} $\implies$ {\rm(ii)}
is obvious by definition.

%%{\rm(ii)} $\implies$ {\rm(iii)}:

\eqref{094a} $\Longrightarrow$ \eqref{095} and \eqref{096}: Let (\ref
{094a}) hold with $t_k\downarrow0$ and $b_k>0$.
Take a further subsequence $t_{k'}\downarrow0$
if necessary so that
$X_{t_{k'}}/b_{k'} \stackrel{\mathrm{D}}{\longrightarrow}Z'$.
$Z'$ is infinitely divisible by
Lemma 4.1 of Maller and Mason \cite{MM08}.
Then $\llvert Z'\rrvert =1$ a.s., thus, as a bounded infinitely
divisible random variable,
$Z'$ is degenerate at a constant which must be $\pm1$. When $Z = +1$,
$X \in \mathrm{PSRS}$. Apply Theorem \ref{SPRS} to get \eqref{095} and
\eqref{096}. If $Z = -1$, $-X \in \mathrm{PSRS}$. Then apply Theorem \ref
{SPRS} to $-X$ to get \eqref{095} and \eqref{096} again.

%%{\rm(iii)} $\implies$ {\rm(iv)}:

\eqref{095} or \eqref{096} $\Longrightarrow$ $X_t\in SRS$ at
0: Let
\eqref{095} or \eqref{096} hold.
Then there is a sequence $x_k\downarrow0$ as $k\to\infty$ such that
$\llvert A(x_k)\rrvert >0$. By taking a further subsequence, we may assume
that $A(x_k)>0$ for all $k$ or $A(x_k)<0$ for all $k$.
Suppose the former; then \eqref{sprs3} or \eqref{sprs3a} holds, so we
get $X\in \mathrm{PSRS}$ by Theorem \ref{SPRS}. If the latter, then by applying
Theorem \ref{SPRS} to $-X$, we get $X\in \mathrm{NSRS}$.
\end{pf}

%%{\rm(iv)}$\implies$ {\rm(i)}:
%\eqref{096} $\implies$ $X_t\in SRS$ at 0: %By the same argument as
%above, apply Theorem \ref{SPRS} to $X$ and
%$-X$ respectively according to the sign of $A(x)$. Thus, $X \in SRS$.
%%\end{pf*}

%s5.5 #&#
\subsection{Subsequential attraction to normality and dominance}\label{sub45}
We can also have subsequential convergence to normality, as
$t\downarrow0$.
The next theorem gives an ``uncentered'' version of this.
We describe \eqref{06p} as ``$X\in D_{P0}(N)$ at 0''.

%
%th10 #&#
\begin{theorem}\label{thsan}
Assume $\sigma^2>0$ or $\pibar(0+)=\infty$.
Then there are nonstochastic sequences
$t_k \downarrow0$ and $b_k\downarrow0$ such that, as $k\to\infty$,
%
%
%e5.79 #&#
\begin{equation}
\label{06p} \frac{X_{t_k}}{b_k} \stackrel{\mathrm{D}} {\longrightarrow}N(0,1);
\end{equation}
iff
%
%
%e5.80 #&#
\begin{equation}
\label{082} \limsup_{x \downarrow0} \frac{U(x)}{x^2 \pibar(x)+x\llvert
A(x)\rrvert } = \infty.
\end{equation}
\end{theorem}

\begin{pf}%{Proof of Theorem \ref{thsan}}
%\begin{proof}[of Theorem {\rm
%\ref{thsan}}]
Both conditions hold when $\sigma^2>0$, so
we can assume $\sigma^2=0$, thus, $\pibar(0+)=\infty$.
Let \eqref{082} hold and choose $x_k\downarrow0$ such that
%
%
%e5.81 #&#
\begin{equation}
\label{52} \frac{U(x_k)}{x_k^2\pibar(x_k)}\to\infty\quad\mbox{and}\quad
\frac{U(x_k)}{x_k\llvert A(x_k)\rrvert }\to\infty.
\end{equation}
Then define
%
%
%e5.82 #&#
\begin{eqnarray}
\label{53} t_k&=& \min\biggl\{ \sqrt{
\frac{x_k^2}{\pibar(x_k)U(x_k)}}, \sqrt{\frac{x_k^3}{\llvert
A(x_k)\rrvert U(x_k)}} \biggr\}.
\end{eqnarray}
(If $A(x_k)=0$ interpret the second component in
\eqref{53} as $+\infty$.)
Thus,
\[
t_k \pibar(x_k) \le\sqrt{
\frac{x_k^2\pibar(x_k)}{U(x_k)}}\to0,
\]
and since $\pibar(0+)>0$, we have $t_k\to0$ as $k\to\infty$.
Now let
\[
b_k^2=t_kU(x_k).
\]
Since $\sigma^2=0$, $U(x_k) = 2 \int_0^{x_k} y \pibar(y )\,\rmd y \to
0 $ as $k \to\infty$.
Then $b_k\to0$ as $k\to\infty$.
Also
\[
\frac{b_k^2}{x_k^2} =\min\biggl\{ \sqrt{\frac{U(x_k)}{x_k^2\pibar
(x_k)}}, \sqrt{
\frac{U(x_k)}{x_k\llvert A(x_k)\rrvert }} \biggr\} \to\infty\qquad
\bigl(\mbox{by \eqref{52}}\bigr).
\]
Given $x > 0$ choose $k$ so large that $x b_k \geq x_k$. Then
\[
%%\label{F--}
t_k \pibar(x b_k) \leq t_k
\pibar(x_k) \to0,
\]
and
%
%
%e5.83 #&#
\begin{eqnarray}
\label{U-} \frac{t_k U(x b_k)}{b_k^2} &=& 1 + \frac{t_k (U(x b_k)
-U(x_k) )}{b_k^2} = 1 +
\frac{2t_k \int_{x_k}^{x b_k} y \pibar(y) \,\rmd y}{b_k^2}
\nonumber\\[-8pt]\\[-8pt]\nonumber
&=& 1 + O\bigl(t_k \pibar(x_k)\bigr) = 1 + o(1).
\end{eqnarray}
Also
\[
\frac{t_k\llvert A(x_k)\rrvert }{x_k} \le\sqrt{ \frac{x_k\llvert
A(x_k)\rrvert }{U(x_k)} } \to0,
\]
while
%
%
%e5.84 #&#
\begin{eqnarray}
\label{A-} \frac{t_k\llvert A(b_k)\rrvert }{b_k} &\le& o \biggl(\frac
{t_k\llvert A(x_k)\rrvert }{x_k} \biggr) +
\frac{t_k\llvert \int_{x_k}^{b_k} (\pibar^+(y)-\pibar
^-(y) ) \,\rmd y\rrvert }{b_k}
\nonumber\\[-8pt]\\[-8pt]\nonumber
&\le& o(1)+t_k\pibar(x_k)=o(1).
\end{eqnarray}
It follows from \eqref{U-}, \eqref{A-} and the subsequential version
of \eqref{normcon} that $X_{t_k}/b_k \stackrel{\mathrm
{D}}{\longrightarrow}N(0,1)$.

Conversely, if there is a $t_k\downarrow0$ such that
$X_{t_k}/b_k \stackrel{\mathrm{D}}{\longrightarrow}N(0,1)$, then by
the subsequential version of
\eqref{normcon} we get \eqref{082}.
\end{pf}

Our final result in this section shows that a 2-sided version of \eqref
{sprs4} holds iff $X\in D_{P0}(N)$ at 0 or $X\in SRS$ at 0.

%
%th11 #&#
\begin{theorem}\label{SPRS3}
Assume $\pibar(0+)=\infty$. Then the following are equivalent:
%
%e5.85 #&#
\begin{eqnarray}
&& \mbox{there is a nonstochastic sequence } t_k
\downarrow0
\nonumber\\[-8pt]\label{smrs4}  \\[-8pt]\nonumber
&&\quad \mbox{such that } \frac{\llvert X_{t_k}\rrvert
}{\llvert \widetilde{\Delta X}_{t_k}^{(1)}\rrvert }\buildrel P
\over\to \infty,\qquad\mbox{as }k\to\infty;
\\
%
%e5.86 #&#
\label{smrs2} &&\limsup_{x\downarrow0} \frac{x\llvert A(x)\rrvert
+U(x)}{x^2\pibar(x)} =\infty;
\\
%
%
%e5.87 #&#
\label{smrs3} &&\mbox{\textup{(a)} }\limsup_{x \downarrow0}
{U(x) \over x \llvert A(x)\rrvert +x^2\pibar(x)} =+\infty,
\quad\mbox{or}\quad\mbox{\textup{(b)} }\limsup
_{x \downarrow0} {x\llvert A(x)\rrvert \over U(x)} =+\infty;
%e5.88 #&#
\\
\label{smrs5}
&& X\in D_{P0}(N)\cup SRS\qquad\mbox{at } 0.
\end{eqnarray}
\end{theorem}

%%\bigskip
\begin{pf} %{Proof of Theorem \ref{SPRS3}} %\begin{proof}[of Theorem {
%\rm\ref{SPRS3}}]
Assume $\pibar(0+)=\infty$.

\eqref{smrs4} $\Longrightarrow$ \eqref{smrs2}: Assume \eqref{smrs4}.
Then just as in the proof of Theorem \ref{PRS3}, we find
$t_k\pibar(b_{k})\to0$ as $k\to\infty$ where $b_k$ satisfies \eqref{mrs14}.
Thus, \eqref{smrs2} holds.

\eqref{smrs2} $\Longrightarrow$ \eqref{smrs3} follows from Theorem 3 of
Maller \cite{PLMS09}.

\eqref{smrs3} $\iff$ \eqref{smrs5}: follows from Theorem \ref
{thsan} and Corollary \ref{th4}.

\eqref{smrs5} $\Longrightarrow$ \eqref{smrs4}: \eqref{smrs5} implies
that there are $t_k\downarrow0$, $b_k\downarrow
0$ such
that $X_{t_k}/b_k\stackrel{\mathrm{D}}{\longrightarrow}N(0,1)$ or
$\llvert X_{t_k}\rrvert /b_k\buildrel P \over\to1$ as $k\to
\infty$.
Either of these implies $t_k\pibar(b_k)\to0$ as $k\to\infty$
and hence $\sup_{0<s\le t_k}\llvert \Delta X_s\rrvert /b_k\buildrel
P \over\to
0$ as $k\to\infty$.
Thus, \eqref{smrs4} holds.
%%That \eqref{smrs5} implies \eqref{smrs4} follows by the same argument.
%% asshowing that \eqref{mrs3} implies \eqref{mrs4} in the proof of
%Theorem \ref{PRS3}.
\end{pf}

%
%re5 #&#
\begin{remark}\label{rem6} {\rm(i) Theorems \ref{thsan} and \ref{SPRS3}
have deep connections to
generalised iterated logarithm laws for $X_t$ as $t\downarrow0$.
It is shown in Theorem 3 of Maller \cite{PLMS09} that \eqref{smrs2}
is equivalent to the existence of a nonstochastic function $B_t>0$ such that
\[
\limsup_{t\downarrow0} \frac{\llvert X_t\rrvert }{B_t} =1\qquad\mbox{a.s. }
\]
Maller \cite{PLMS09} also gives a.s. equivalences for \eqref{cors2} and
\eqref{smrs3}(a). We hope to consider a.s. results related to those in
Sections~\ref{s3}--\ref{s4} elsewhere.

(ii) We note that in many conditions such as \eqref{smrs2} and \eqref
{smrs3} we may replace the functions
$A(x)$ and $U(x)$ in \eqref{Adef} and \eqref{Udef} by the functions
$\nu(x)$ and $V(x)$ in \eqref{nudef}. This is because
\[
x\bigl\llvert A(x)-\nu(x)\bigr\rrvert\le x^2\pibar(x) \quad
\mbox{and}\quad0\le U(x)-V(x)=x^2\pibar(x),\qquad x>0.
\]
But there is some advantage to working with the continuous functions
$A(x)$ and $U(x)$, and sometimes it is essential, for example, in
Theorem \ref{th33}.
}
\end{remark}

%s6 #&#
\section{Related large time results}\label{s8}
Most of the small time results derived herein have exact or close
analogues for large times (i.e., allowing $t\to\infty$ rather than
$t\downarrow0$), some of them having been suggested by such analogies. In
fact, many of the identities hold generally, for all $t>0$; this is the
case for all results in Section~\ref{s2}, as well as
Lemmas \ref{lemtsum2} and \ref{lemtsum4}.
Some analogous large time results for L\'evy processes can be found in
Kevei and Mason \cite{KM}, and Maller and Mason \cite{MM09,MM13}, and we expect that others can be derived by straightforward
modification of our small time methods.
These would include compound Poisson processes as special cases.
%%, and analogous random walk versions.

%For a random walk $S_n=\sum_{i=1}^n\xi_i$ comprised of i.i.d.
%increments $\xi_i$, ``large time" means ``as $n\to\infty$", that is,
%also, ``large sample". Among known random walk results we mention
%Maller and Resnick \cite{MR} (which provided inspiration for
%Theorems \ref{tX}, \ref{pmtX} and \ref{+-}, in particular), and Mori
%\cite{mori}.
%Theorem \ref{randrep} can be used to transfer results of Arov and
%Bobrov \cite{ab} and Darling \cite{darl} concerning trimmed random
%walks and related order statistics from random walks to L\'evy
%processes.

\section*{Acknowledgements}
We are grateful to a referee for a very careful reading and for
suggesting substantial improvements to the original version of the paper.
R. Maller's research was partially supported by ARC grant DP1092502.

%\begin{appendix}
%\section{}
%\end{appendix}

%\section*{Acknowledgements}

%\begin{supplement}%[id=suppA]
%\slink[doi]{10.3150/00-BEJXXXXSUPP} %[doi,text={...}] - jei reikia
%suskaldyti doi
%\sfilename{BEJ000\_supp.pdf}
%\end{supplement}

% imsref loaded by linak, 2015-08-14 15:13:28
%
% imsref loaded by linak, 2015-08-17 16:15:57

\printhistory
\end{document}